\documentclass[11pt]{article}
\usepackage{amsfonts,amsmath}
\usepackage{amssymb,mathrsfs,latexsym,bm,cite}
\usepackage{arydshln}

\begin{document}
\newcommand{\qed}{\hphantom{.}\hfill $\Box$\medbreak}
\newcommand{\proof}{\noindent{\bf Proof \ }}
\newtheorem{Theorem}{Theorem}[section]
\newtheorem{Lemma}[Theorem]{Lemma}
\newtheorem{Corollary}[Theorem]{Corollary}
\newtheorem{Remark}[Theorem]{Remark}
\newtheorem{Example}[Theorem]{Example}
\newtheorem{Definition}[Theorem]{Definition}
\newtheorem{Construction}[Theorem]{Construction}

\thispagestyle{empty}
 \renewcommand{\thefootnote}{\fnsymbol{footnote}}

\begin{center}
{\Large\bf Semi-cyclic holey group divisible designs with block size three\footnote{Supported by the Fundamental Research Funds for the Central Universities under Grants $2011$JBM$298$, $2013$JB$Z005$ (T. Feng), $2011$JBZ$012$ (Y. Chang), and the NSFC under Grant $10901016$ (T. Feng), the NSFC under Grants $11126123$ and $11201252$ (X. Wang), the NSFC under Grant $11271042$ (Y. Chang).}}

\vskip12pt

Tao Feng$^1$, Xiaomiao Wang$^2$ and Yanxun Chang$^1$ \\[2ex] {\footnotesize {\footnotesize $^1$Institute of Mathematics, Beijing Jiaotong
University, Beijing 100044, P. R. China}
\\$^2$Department of Mathematics, Ningbo University, Ningbo 315211, P. R. China}\\  {\footnotesize
tfeng@bjtu.edu.cn, wangxiaomiao@nbu.edu.cn, yxchang@bjtu.edu.cn}
\vskip12pt

\end{center}

\vskip12pt

\noindent {\bf Abstract:} In this paper we discuss the existence problem for a semi-cyclic holey group divisible design of type $(n,m^t)$ with block size $3$, which is denoted by a $3$-SCHGDD of type $(n,m^t)$. When $n=3$, a $3$-SCHGDD of type $(3,m^t)$ is equivalent to a $(3,mt;m)$-cyclic holey difference matrix, denoted by a $(3,mt;m)$-CHDM.

It is shown that there is a $(3,mt;m)$-CHDM if and only if $(t-1)m\equiv 0\ ({\rm mod}\ 2)$ and $t\geq 3$ with the exception of $m\equiv 0\ ({\rm mod}\ 2)$ and $t=3$. When $n\geq 4$, the case of $t$ odd is considered. It is established that if $t\equiv 1\ ({\rm mod}\ 2)$ and $n\geq 4$, then there exists a $3$-SCHGDD of type $(n,m^t)$ if and only if $t\geq 3$ and $(t-1)n(n-1)m\equiv 0\ ({\rm mod}\ 6)$ with some possible exceptions of $n=6$ and $8$. The main results in this paper have been used to construct optimal two-dimensional optical orthogonal codes with weight $3$ and different auto- and cross-correlation constraints by the authors recently. \\

\noindent {\bf Keywords}: holey group divisible design; semi-cyclic; cyclic holey difference matrix; cyclic difference family; perfect difference family \\

\noindent {\bf Mathematics Subject Classification $(2000)$}: $05$B$30$ $\cdot$ $94$B$25$


\section{Introduction}

A holey group divisible design is an important combinatorial configuration. Its study has been motivated by many applications in constructing various types of combinatorial objects, for example \cite{a2,bcrw,cd,gww}. Throughout this paper we always assume that $I_u=\{0,1,\ldots,u-1\}$ and denote by $Z_v$ the additive group of integers modulo $v$.

Let $n,m$ and $t$ be positive integers. Let $K$ be a set of positive integers. A {\em holey group divisible design} (HGDD) $K$-HGDD is a quadruple $(X,{\cal G},{\cal H},{\cal B})$ which satisfies the following properties:
\begin{enumerate}
    \item[(1)] $X$ is a finite set of $nmt$ {\em points};
    \item[(2)] ${\cal G}$ is a partition of $X$ into $n$ subsets (called {\em groups}), each of size $mt$;
    \item[(3)] ${\cal H}$ is another partition of $X$ into $t$ subsets (called {\em holes}), each of size $nm$ such that $|H\cap G|=m$ for each $H\in {\cal H}$ and $G\in {\cal G}$;
    \item[(4)] $\cal B$ is a set of subsets (called {\em blocks}) of $X$, each of
    cardinality from $K$, such that no block contains two distinct points of any group or any hole, but any other pair of distinct points of $X$ occurs in exactly one block of $\cal B$.
\end{enumerate}
Such a design is denote by a $K$-HGDD of type $(n,m^t)$. When $m=1$, a $K$-HGDD of type $(n,1^t)$ is often said to be a {\em modified group divisible design}, denoted by a $K$-MGDD of type $t^n$. If $K=\{k\}$, we write a $\{k\}$-HGDD as a $k$-HGDD, and a $\{k\}$-MGDD as a $k$-MGDD.

Assaf \cite{a} first introduced the notion of MGDDs and settled the existence of $3$-MGDDs. The existence of $4$-MGDDs was investigated in \cite{a1,aw,lc}. There are also some results on $5$-MGDDs in \cite{aa,aa1}. Wei \cite{wei} first introduced the concept of HGDDs and gave a complete existence theorem for $3$-HGDDs. The existence of $4$-HGDDs has been completely settled in \cite{cww,gw}. We only quote the following result for the later use.

\begin{Theorem}\label{3-HGDD}{\rm \cite{wei}}
There exists a $3$-HGDD of type $(n,m^t)$ if and only if $n,t\geq 3$, $(t-1)(n-1)m\equiv 0\ ({\rm mod}\ 2)$ and $t(t-1)n(n-1)m^2\equiv 0\ ({\rm mod}\ 3)$.
\end{Theorem}

HGDDs can be seen as a special case of double group divisible designs (DGDDs), which are introduced in \cite{zhu} to simplify Stinson's proof \cite{s} on a recursive construction for group divisible designs. HGDDs can also be considered as a generalization of holey mutually orthogonal Latin squares (HMOLS), since a $k$-HGDD of type $(k,m^t)$ is equivalent to $k-2$ HMOLSs of type $m^t$ (see \cite{sz,wei} for details).

A way to construct $k$-HGDDs of type $(n,m^t)$ is the pure and mixed difference method. Let $S=\{0,t,\ldots,(m-1)t\}$ be a subgroup of order $m$ in $Z_{mt}$, and $S_l=S+l$ be a coset of $S$ in $Z_{mt}$, $0\leq l\leq t-1$. Let $X=I_n\times Z_{mt}$, ${\cal G}=\{\{i\}\times Z_{mt}:i\in I_n\}$, and ${\cal H}=\{I_n\times S_l:0\leq l\leq t-1\}$. Take a family ${\cal B}^*$ of some $k$-subsets (called {\em base blocks}) of $X$. For each $i,j\in I_n$ and $B\in{\cal B}^*$, define a multi-set $\Delta_{ij}(B)=\{x-y\ ({\rm mod}\ mt): (i,x),(j,y)\in B, (i,x)\neq(j,y)\}$, and a multi-set $\Delta_{ij}({\cal B}^*)=\bigcup_{B\in{\cal B}^*}\Delta_{ij}(B)$. If
$$\Delta_{ij}({\cal B}^*)=\left\{
\begin{array}{lll}
Z_{mt}\setminus S, & {\rm if\ } (i,j)\in I_n\times I_n {\rm \ and\ } i\neq j,\\
\emptyset, & {\rm otherwise},\\
\end{array}
\right.
$$
then a $k$-HGDD of type $(n,m^t)$ with the point set $X$, the group set $\cal G$ and the hole set $\cal H$ can be generated from ${\cal B}^*$. The required blocks are
obtained by developing all base blocks of ${\cal B}^*$ by successively
adding $1$ to the second component of each point of these base
blocks modulo $mt$. Usually a $k$-HGDD obtained by this
manner is said to be a {\em semi-cyclic $k$-HGDD} and denoted by a $k$-SCHGDD. $k$-SCHGDDs were first introduced in \cite{wsy} to construct optimal two-dimensional optical orthogonal codes.

\begin{Example}\label{3-1^t}
Here we construct a $3$-SCHGDD of type $(3,1^t)$ on $I_3\times Z_t$ for any odd integer $t\geq 3$. Let $S=\{0\}$ be a subgroup of order $1$ in $Z_t$, and $S_l=\{l\}$ be a coset of $S$ in $Z_t$, $0\leq l\leq t-1$. Take the group set ${\cal G}=\{\{i\}\times Z_t:i\in I_3\}$, and the hole set ${\cal H}=\{I_3\times S_l:0\leq l\leq t-1\}$. Let the base block set ${\cal B}^*=\{\{(0,0),(1,x),(2,2x)\}:1\leq x\leq t-1\}$. It is readily checked that for each $(i,j)\in I_3\times I_3$ and $i\neq j$, $\Delta_{ij}({\cal B}^*)=Z_t\setminus\{0\}$, and for each $(i,j)\in I_3\times I_3$ and $i=j$, $\Delta_{ij}({\cal B}^*)=\emptyset$.
\end{Example}

SCHGDDs are closely related to {\em cyclic holey difference matrices} (CHDMs). A $(k,wt;w)$-CHDM is a $k\times w(t-1)$ matrix $D=(d_{ij})$ with entries from $Z_{wt}$ such that for any two distinct rows $x$ and $y$, the difference list $L_{xy}=\{d_{xj}-d_{yj}:j\in I_{w(t-1)}\}$ contains each integer of $Z_{wt}\setminus S$ exactly once, where $S=\{0,t,\ldots,(w-1)t\}$ is a subgroup of order $w$ in $Z_{wt}$, while any integer of $S$ does not appear in $L_{xy}$. A $(k,wt;w)$-CHDM is also called a {\em $w$-regular $(wt,k,1)$-incomplete difference matrix} and denoted by a $w$-regular ICDM$(k;wt)$ in \cite{cm}.

\begin{Lemma}\label{relation SCH and CHDM} {\rm \cite{wy}}
A $k$-SCHGDD of type $(k,w^t)$ is equivalent to a $(k,wt;w)$-CHDM.
\end{Lemma}

\begin{Example}\label{3-t-1-CHDM} Example $\ref{3-1^t}$ presents a $3$-SCHGDD of type $(3,1^t)$ for any odd integer $t\geq 3$, which yields a $(3,t;1)$-CHDM $(A_1|A_2|\cdots|A_{t-1})$, where $A_i=(0,i,2i)^T$, $1\leq i\leq t-1$.
\end{Example}

In this paper we shall focus on the existence problem of $3$-SCHGDDs. It is easy to see that the number of base blocks in a $3$-SCHGDD of type $(n,m^t)$ is $(t-1)n(n-1)m/6$. Hence combining the result of Theorem \ref{3-HGDD}, we have the following necessary condition for the existence of $3$-SCHGDDs.

\begin{Lemma}\label{nece}
If there exists a $3$-SCHGDD of type $(n,m^t)$, then $n,t\geq 3$, $(t-1)(n-1)m\equiv 0\ ({\rm mod}\ 2)$ and $(t-1)n(n-1)m\equiv 0\ ({\rm mod}\ 6)$.
\end{Lemma}

As the main result, we are to prove the following theorems, which will be established at the end of Sections $5$ and $7$, respectively. We remark that these two theorems have been used to construct optimal two-dimensional optical orthogonal codes with weight $3$ and different auto- and cross-correlation constraints \cite{wcf}.

\begin{Theorem}\label{3-chdm}
There exists a $3$-SCHGDD of type $(3,m^t)$ $($i.e., a $(3,mt;m)$-CHDM$)$ if and only if $(t-1)m\equiv 0\ ({\rm mod}\ 2)$ and $t\geq 3$ with the exception of $m\equiv 0\ ({\rm mod}\ 2)$ and $t=3$.
\end{Theorem}

\begin{Theorem}\label{main theorem}
Assume that $t\equiv 1\ ({\rm mod}\ 2)$ and $n\geq 4$. There exists a $3$-SCHGDD of type $(n,m^t)$ if and only if $t\geq 3$ and $(t-1)n(n-1)m\equiv 0\ ({\rm mod}\ 6)$ except when $(n,m,t)=(6,1,3)$, and possibly when $(1)$ $n=6$, $m\equiv 1,5\ ({\rm mod}\ 6)$ and $t\equiv 3,15\ ({\rm mod}\ 18)$, $(2)$ $n=8$, $m\equiv 2,10\ ({\rm mod}\ 12)$ and $t\equiv 7\ ({\rm mod}\ 12)$.
\end{Theorem}

\section{Preliminaries}

Now we try to make the reader familiar with some basic concepts and terminologies in
combinatorial design theory adopted in this paper.

Let $K$ be a set of positive integers. A {\em group divisible
design} (GDD) $K$-GDD is a triple ($X, {\cal G},{\cal B}$)
satisfying the following properties: $(1)$ $X$ is a finite set of {\em points};
$(2)$ $\cal G$ is a partition of $X$ into subsets (called {\em groups});
$(3)$ $\cal B$ is a set of subsets (called {\em blocks}) of $X$, each of
cardinality from $K$, such that every $2$-subset of $X$ is either
contained in exactly one block or in exactly one group, but not in
both. If $\cal G$ contains $u_i$ groups of size $g_i$ for $1\leq
i\leq r$, then we call $g_1^{u_1}g_2^{u_2}\cdots g_r^{u_r}$ the {\em
group type} (or {\em type}) of the GDD. If $K=\{k\}$, we write a
$\{k\}$-GDD as a $k$-GDD. A $K$-$GDD$ of type $1^v$ is commonly called a
{\em pairwise balanced design}, denoted by a $(v,K,1)$-$PBD$. When $K=\{k\}$, a pairwise balanced design is called a {\em balanced incomplete block design}, denoted by a $(v,
k,1)$-BIBD.

\begin{Lemma}
\label{346pbd} {\rm \cite{abg}}
\begin{enumerate}
\item[$(1)$] There exists a $(v,\{3,4\},1)$-PBD for any integer $v\equiv 0,1\ ({\rm mod}\ 3)$ and $v\geq 3$ with the exception of $v=6$.
\item[$(2)$] There exists a $(v,\{3,4,5\},1)$-PBD for any integer $v\geq 3$ with the exception of $v=6,8$.
\item[$(3)$] There exists a $(v,\{4,6,7,9\},1)$-PBD for any integer $v\equiv 0,1\ ({\rm mod}\ 3)$, $v\geq 4$ and $v\not\in \{10,12,15,18,19,24,27\}$.
\end{enumerate}
\end{Lemma}

An {\em automorphism} $\pi$ of a GDD $(X,{\cal G},{\cal B})$ is a permutation on $X$ leaving ${\cal G}$, ${\cal B}$ invariant, respectively. Let $H$ be the cyclic group
generated by $\alpha$ under the compositions of permutations. Then all blocks of the GDD can be partitioned into some block orbits under $H$. Choose any fixed block from each block orbit and then call it a {\em base block} of this GDD under $H$. The number of the blocks contained in a block orbit is called the {\em length} of the block orbit.

A $K$-GDD of type $m^n$ is said to be {\em cyclic}, if it admits an automorphism consisting of a cycle of length $mn$. A cyclic $K$-GDD is denoted by a $K$-CGDD. For a $K$-CGDD $(X,{\cal G},{\cal B})$, we can always identify $X$ with $Z_{mn}$ and ${\cal G}$ with $\{\{in+j : 0\leq i\leq m-1\} : 0\leq j\leq n-1\}$. If the length of each block orbit in a $K$-CGDD of type $m^n$ is $mn$, then the $K$-GDD is called {\em strictly cyclic}.

\begin{Lemma}\label{3-sCGDD}{\rm \cite{wc}}
There exists a strictly cyclic $3$-GDD of type $m^n$ if and only if
\begin{enumerate}
\item[$(1)$] $m(n-1)\equiv 0\ ({\rm mod}\ 6)$ and $n\geq 4;$
\item[$(2)$] $n\not\equiv 2,3\ ({\rm mod}\ 4)$ when $m\equiv 2\ ({\rm mod}\ 4).$
\end{enumerate}
\end{Lemma}

A $K$-GDD of type $m^n$ is said to be {\em semi-cyclic}, if it admits an automorphism which permutes the elements of each group $G\in{\cal G}$ in an $m$ cycle. Such a GDD is denoted by a $K$-SCGDD of type $m^n$. For a $K$-SCGDD $(X,{\cal G},{\cal B})$, we can always identify $X$ with $I_n\times Z_m$ and ${\cal G}$ with $\{\{i\}\times Z_m:i\in I_n\}$. In this case the automorphism can be taken as $(i,x)\longmapsto(i,x+1)$ (mod $(-,m)$), $i\in I_n$ and $x\in Z_m$. Assume that ${\cal B}^*$ is the set of base blocks of a $K$-SCGDD of type $m^n$. It is easy to verify that
$$\Delta_{ij}({\cal B}^*)=\left\{
\begin{array}{lll}
Z_m, & {\rm if\ } (i,j)\in I_n\times I_n {\rm \ and\ } i\neq j,\\
\emptyset, & {\rm otherwise}.\\
\end{array}
\right.
$$
Note that no element of $Z_m$ occurs more than once in $\Delta_{ij}({\cal B}^*)$, and the length of each block orbit is $m$ in a $K$-SCGDD.

\begin{Lemma}{\rm\cite{jiang}}\label{scgdd}
There is a $3$-SCGDD of type $m^n$ if and only if $n\geq 3$ and
\begin{enumerate}
\item[$(1)$] $n\equiv 1,3\ ({\rm mod}\ 6)$ when $m\equiv 1,5\ ({\rm
mod}\ 6);$
\item[$(2)$] $n\equiv 1\ ({\rm mod}\ 2)$ when $m\equiv 3\ ({\rm mod}\ 6);$
\item[$(3)$] $n\equiv 0,1,4,9\ ({\rm mod}\ 12)$ when $m\equiv 2,10\ ({\rm mod}\ 12);$
\item[$(4)$] $n\equiv 0,1\ ({\rm mod}\ 3)$, $n\neq 3$ when $m\equiv 4,8\ ({\rm mod}\ 12);$
\item[$(5)$] $n\neq 3$ when $m\equiv 0\ ({\rm mod}\ 12);$
\item[$(6)$] $n\equiv 0,1\ ({\rm mod}\ 4)$ when $m\equiv 6\ ({\rm mod}\ 12)$.
\end{enumerate}
\end{Lemma}

\begin{Lemma}{\rm\cite{wy,ww}}\label{4-scgdd}
There is a $4$-SCGDD of type $m^n$ if and only if $n\geq 4$, $m(n-1)\equiv 0\ ({\rm mod}\ 3)$ and $mn(n-1)\equiv 0\ ({\rm mod}\ 12)$, except when $n=4$ or $(m,n)\in\{(2,10),(4,7),(6,5)\}$, and possibly when $(1)$ $n=5$, $m\equiv \pm6\ ({\rm mod}\ 36)$ and $m\geq 30$, $(2)$ $n=7$, $m\equiv \pm4\ ({\rm mod}\ 24)$ and $m\geq 20$, $(3)$ $n=10$, $m\equiv \pm2\ ({\rm mod}\ 12)$ and $m\geq 10$.
\end{Lemma}

\begin{Lemma}{\rm\cite{wy}}\label{SCGDD from SCHGDD}
Suppose that a $k$-SCHGDD of type $(n,m^t)$ and a $k$-SCGDD of type $m^n$ exist. Then a $k$-SCGDD of type $(mt)^n$ exists.
\end{Lemma}

We remark that $k$-SCGDDs are closely related to {\em cyclic difference matrices} (CDMs). A $(k,m)$-CDM is a $k\times m$ matrix $D=(d_{ij})$ with entries from $Z_m$ such that for any two distinct rows $x$ and $y$, the difference list $\{d_{xj}-d_{yj}:j\in I_m\}$ contains each integer of $Z_m$ exactly once.

\begin{Lemma}\label{relation SC and CDM} {\rm \cite{yin2002}}
A $k$-SCGDD of type $m^k$ is equivalent to a $(k,m)$-CDM.
\end{Lemma}

Let $g,t$ be positive integers and $K$ be a set of positive integers. A {\em $(gt,g,K,1)$-cyclic difference family} (briefly $(gt,g,K,1)$-CDF) is a family $\mathscr{F}$ of subsets
(called {\em base blocks}) of $Z_{gt}$ such that $(1)$ if $B\in
\mathscr{F}$, then $|B|\in K$; $(2)$
$\Delta\mathscr{F}=\bigcup_{B\in \mathscr{F}}\Delta B$ covers
each element of $Z_{gt}\setminus\{0,t,2t,\ldots,(g-1)t\}$ exactly once, where
$\Delta B=\{x-y:\, x,y\in B, x\neq y\}$. It is clear that all the base blocks of a $(gt,g,K,1)$-CDF constitute the base block set of a strictly cyclic $K$-GDD of type $g^t$.

Let $h$ be a positive integer. Assume that ${\cal F}=\{\{0,b_{1i},b_{2i},\ldots,b_{k_i-1,i}\}:i=1,2,\ldots,r\}$ is the family of base blocks of an $(hgt,hg,K,1)$-CDF, where $k_i\in K$. Define $ele({\cal F})=\bigcup_{i=1}^r\{b_{1i},b_{2i},\ldots,\linebreak b_{k_i-1,i}\}.$ The $(hgt,hg,K,1)$-CDF is said to be {\em $h$-perfect}, denoted by an $(hgt,hg,K,1)$-$h$-PDF, if $ele(\mathcal{F})\subseteq\{a+bgt:\ 0\leq a\leq
\lfloor gt/2\rfloor,a\not\equiv 0\ ({\rm mod}\ t),
b=0,1,\ldots,h-1\}$. As usual, when $h=1$, a $(gt,g,K,1)$-$1$-PDF is abbreviated to a $(gt,g,K,1)$-PDF. Furthermore, if $g=1$, a $(t,1,K,1)$-PDF is simply written as a $(t,K,1)$-PDF.

To illustrate these definitions, we give the following examples, which are useful later.

\begin{Example}\label{example-CDF} There is a $(256,32,\{3,5\},1)$-CDF as follows:
\begin{center}\tabcolsep 0.07in
\begin{tabular}{lllllll}
$\{0,1,45\},$ &$\{0,7,58\},$ &$\{0,17,172\},$ &$\{0,33,181\},$ &$\{0,38,179\},$&$\{0,25,191\},$\\
$\{0,2,49\},$ &$\{0,9,130\},$ &$\{0,18,132\},$ &$\{0,35,178\},$ &$\{0,62,165\},$&$\{0,26,164\},$\\
$\{0,3,46\},$ &$\{0,10,177\},$ &$\{0,20,189\},$ &$\{0,36,183\},$ &$\{0,27,186\},$&$\{0,42,55\},$\\
$\{0,4,54\},$ &$\{0,11,182\},$ &$\{0,22,173\},$ &$\{0,37,187\},$ &$\{0,28,185\},$&$\{0,61,163\},$\\
$\{0,5,57\},$ &$\{0,12,174\},$ &$\{0,23,134\},$ &$\{0,30,170\},$ &$\{0,39,60\},$&$\{0,15,34,161,190\},$\\
$\{0,6,59\},$ &$\{0,14,133\},$ &$\{0,63,131\},$ &$\{0,31,180\},$ &$\{0,41,158\}.$&\\
\end{tabular}
\end{center}
\end{Example}

\begin{Example}\label{example-PDF} For $(t,K)\in\{(13,\{4\})$, $(33,\{3,5\})$, $(59,\{3,4,5\})$, $(71,\{3,5\})\}$, there is a $(t,K,1)$-PDF as follows:
\begin{center}\tabcolsep 0.035in
\begin{tabular}{llllll}
$(t,K)=(13,\{4\}):$&
$\{0,1,4,6\}$.\\
$(t,K)=(33,\{3,5\}):$&
$\{0,1,6,14,16\}$,&
$\{0,3,12\}$,&
$\{0,4,11\}$.\\
$(t,K)=(59,\{3,4,5\}):$&
$\{0,1,6,23,26\}$,&
$\{0,7,15,19,28\}$,&
$\{0,11,27,29\}$,&
$\{0,14,24\}$.\\
$(t,K)=(71,\{3,5\}):$&
$\{0,1,3,10,27\}$,&
$\{0,4,15,29,35\}$,&
$\{0,5,28\}$,&
$\{0,8,30\}$,\\&
$\{0,12,33\}$,&
$\{0,13,32\}$,&
$\{0,16,34\}$.
\end{tabular}
\end{center}
\end{Example}

\begin{Example}\label{example-PRDF} For $(g,t)\in\{(8,8),(16,6)\}$, there is a $(gt,g,\{3,5\},1)$-PDF as follows:
\begin{center}\tabcolsep 0.035in
\begin{tabular}{llllll}
$(g,t)=(8,8):$&
$\{0,1,3,7,26\}$,& $\{0,5,27\}$, &$\{0,9,21\}$, &$\{0,10,28\}$,& $\{0,11,31\}$,\\& $\{0,13,30\}$,& $\{0,14,29\}$.\\
$(g,t)=(16,6):$&
$\{0,1,4,11,26\}$,& $\{0,8,40\}$,& $\{0,13,44\}$,& $\{0,16,45\}$, &$\{0,19,47\}$,\\&$\{0,2,35\}$,&
$\{0,5,39\}$,&$\{0,9,46\}$,& $\{0,14,41\}$,&$\{0,17,38\}$,\\&$\{0,20,43\}$.&\\
\end{tabular}
\end{center}
\end{Example}

\begin{Example}\label{example-h-PDF} For $(g,t)\in\{(8,8),(16,6)\}$, there is a $(2gt,2g,\{3,5\},1)$-$2$-PDF as follows:
\begin{center}\tabcolsep 0.01in
\begin{tabular}{llllll}
$(g,t)=(8,8):$ & $\{0,1,4,21,30\}$,& $\{0,2,7,25,67\}$,&$\{0,11,90\}$, &$\{0,13,94\}$,&$\{0,15,91\}$,\\&$\{0,22,66\}$,&$\{0,6,75\}$,&$\{0,10,87\}$,&$\{0,12,95\}$,&$\{0,14,92\}$,\\&$\{0,19,93\}$,&
$\{0,27,73\}$,&
$\{0,28,71\}$,& $\{0,31,70\}$.\\
$(g,t)=(16,6):$ &
$\{0,1,14,47,135\}$,&$\{0,15,41,97,142\}$,&$\{0,7,131\}$,&$\{0,19,128\}$,&$\{0,25,137\}$,\\&
$\{0,2,11,31,118\}$,&$\{0,17,40,103,140\}$,&$\{0,10,44\}$,& $\{0,21,143\}$,&$\{0,28,141\}$,\\&
$\{0,3,8,35,133\}$,&$\{0,4,43\}$,&$\{0,16,115\}$,&$\{0,22,139\}$,&$\{0,38,111\}$.\\
\end{tabular}
\end{center}
\end{Example}

\begin{Example}\label{example-h-PDF-h=4} There is a $(96,16,\{3,5\},1)$-$4$-PDF as follows:
\begin{center}
\begin{tabular}{llllll}
$\{0,1,4,11,26\}$,&$\{0,2,29\}$,&
$\{0,5,33\}$,&
$\{0,8,59\}$,&
$\{0,9,49\}$,&
$\{0,31,83\}$,\\
$\{0,32,82\}$,&
$\{0,34,53\}$,&
$\{0,35,55\}$,&
$\{0,57,73\}$,&
$\{0,58,75\}$.
\end{tabular}
\end{center}
\end{Example}

$h$-Perfect difference families were first introduced by Chang and Miao \cite{cm} as a generalization of perfect difference families to establish a recursive construction for optimal optical orthogonal codes. They used a different terminology called $h$-perfect $g$-regular cyclic packing. For more information on perfect difference families, the interested reader may refer to \cite{ab,gms,wc2}.

\section{Construction methods}

In this section we shall establish some recursive constructions for $k$-SCHGDDs.

\begin{Construction}\label{SCHGDD from SCHGDD}
If there exist a $k$-SCHGDD of type $(n,(gw)^t)$ and a $k$-SCHGDD of type
$(n,g^w)$, then there exists a $k$-SCHGDD of type $(n,g^{wt})$.
\end{Construction}

\proof Let $S=\{0,t,\ldots,(gw-1)t\}$ be a subgroup of order $gw$ in $Z_{gwt}$, and $S_l=S+l$ be a coset of $S$ in $Z_{gwt}$, $0\leq l\leq t-1$. By assumption we can construct a $k$-SCHGDD of type $(n,(gw)^t)$ on $I_n\times Z_{gwt}$ with the group set $\{\{i\}\times Z_{gwt}:i\in I_n\}$ and the hole set $\{I_n\times S_l:0\leq l\leq t-1\}$. Denote the set of its base blocks by ${\cal B}^*$.

Let $S'=\{0,wt,\ldots,(g-1)wt\}$ be a subgroup of order $g$ in $S$, and $S'_r=S'+rt$ be a coset of $S'$ in $S$, $0\leq r\leq w-1$. By assumption we can construct a $k$-SCHGDD of type $(n,g^w)$ on $I_n\times S$ with the group set $\{\{i\}\times S:i\in I_n\}$ and the hole set $\{I_n\times S'_r:0\leq r\leq w-1\}$. Denote the set of all its base blocks by ${\cal A}^*$.

Let $S''_j=S'+j$ be a coset of $S'$ in $Z_{gwt}$, $0\leq j\leq wt-1$. Now we construct the required $k$-SCHGDD of type $(n,g^{wt})$ on $I_n\times Z_{gwt}$ with the group set $\{\{i\}\times Z_{gwt}:i\in I_n\}$ and the hole set $\{I_n\times S''_j:0\leq j\leq wt-1\}$. It is readily checked that ${\cal B}^*\cup{\cal A}^*$ constitutes all base blocks of the required design. \qed

\begin{Construction} \label{SCHGDD-recur}
Suppose that there exist a $K$-SCGDD of type $g^n$ and an $l$-SCHGDD of type $(k,w^t)$ for each $k\in K$. Then there exists an $l$-SCHGDD of type $(n,(gw)^t)$.
\end{Construction}

\proof Let ${\cal B}^*$ be the set of base blocks of the given $K$-SCGDD of type $g^n$, which is constructed on $I_n\times Z_g$ with the group set $\{\{i\}\times Z_g:i\in I_n\}$.

Let $S=\{0,t,\ldots,(w-1)t\}$ be a subgroup of order $w$ in $Z_{wt}$, and $S_r=S+r$ be a coset of $S$ in $Z_{wt}$, $0\leq r\leq t-1$. For each $B\in {\cal B}^*$, construct an $l$-SCHGDD of type $(|B|,w^t)$ on $B\times Z_{wt}$ with the group set $\{\{j\}\times Z_{wt}:j\in B\}$ and the hole set $\{B\times S_r:0\leq r\leq t-1\}$. Denote the set of its base blocks by ${\cal A}^*_B$.

Let $S'=\{0,t,\ldots,(gw-1)t\}$ be a subgroup of order $gw$ in $Z_{gwt}$, and $S'_h=S'+h$ be a coset of $S'$ in $Z_{gwt}$, $0\leq h\leq t-1$. Now we construct the required $l$-SCHGDD of type $(n,(gw)^t)$ on $I_n\times Z_{gwt}$ with the group set $\{\{i\}\times Z_{gwt}:i\in I_n\}$ and the hole set $\{I_n\times S'_h:0\leq h\leq t-1\}$. For each $A=\{(a_1,x_1,y_1),(a_2,x_2,y_2),\ldots,(a_l,x_l,y_l)\}\in {\cal A}^*_B$, define
$$\overline{A}=\{(a_1,y_1+wtx_1),(a_2,y_2+wtx_2),\ldots,(a_l,y_l+wtx_l)\}.$$
It is readily checked that $\bigcup_{B\in{\cal B}^*}\{\overline{A}:A\in {\cal A}^*_B\}$ forms the set of base blocks of the required design.  \qed

\begin{Construction} \label{SCHGDD-from strictly CGDD}
Suppose that there exist a strictly cyclic $K$-GDD of type $w^t$ and an $l$-MGDD of type $k^n$ for each $k\in K$. Then there exists an $l$-SCHGDD of type $(n,w^t)$.
\end{Construction}

\proof Let $S=\{0,t,\ldots,(w-1)t\}$ be a subgroup of order $w$ in $Z_{wt}$, and $S_r=S+r$ be a coset of $S$ in $Z_{wt}$, $0\leq r\leq t-1$. Let ${\cal B}^*$ be the set of base blocks of the given strictly cyclic $K$-GDD of type $w^t$, which is constructed on $Z_{wt}$ with the group set $\{S_r:0\leq r\leq t-1\}$.

For each $B\in{\cal B}^*$, we construct an $l$-MGDD of type $|B|^n$ on $I_n\times B$ with the group set $\{\{i\}\times B:i\in I_n\}$ and the hole set $\{I_n\times \{j\}: j\in B\}$. Denote the set of its blocks by ${\cal A}_B$.

Now we construct the required $l$-SCHGDD of type $(n,w^t)$ on $I_n\times Z_{wt}$ with the group set $\{\{i\}\times Z_{wt}:i\in I_n\}$ and the hole set $\{I_n\times S_r:0\leq r\leq t-1\}$. It is readily checked that $\bigcup_{B\in {\cal B}^*}{\cal A}_B$ forms the set of base blocks of the required design. \qed

\begin{Construction}\label{SCHGDD from CDM}
If there exist a $k$-SCHGDD of type
$(n,w^t)$ and a $(k,v)$-CDM, then there exists a $k$-SCHGDD of type $(n,(wv)^t)$.
\end{Construction}

\proof Let ${\cal B}^*$ be the set of base blocks of the given $k$-SCHGDD of type
$(n,w^t)$. Let $D=(d_{ij})$ be the given $(k,v)$-CDM. For each
$B=\{(a_1,x_1),(a_2,x_2),\ldots,(a_k,x_k)\}\in {\cal B}^*$, we construct a family, ${\cal A}^*_B$, consisting of the following $v$ base blocks
$$\{(a_1,x_1+wtd_{1j}),(a_2,x_2+wtd_{2j}),\ldots,(a_k,x_k+wtd_{kj})\},$$
where $0\leq j\leq v-1$. Let $S=\{0,t,\ldots,(wv-1)t\}$ be a subgroup of order $wv$ in $Z_{vwt}$, and $S_l=S+l$ be a coset of $S$ in $Z_{vwt}$, $0\leq l\leq t-1$. Now we construct the required $k$-SCHGDD of type $(n,(wv)^t)$ on $I_n\times Z_{vwt}$ with the group set $\{\{i\}\times Z_{vwt}:i\in I_n\}$ and the hole set $\{I_n\times S_l:0\leq l\leq t-1\}$. It is readily checked that $\bigcup_{B\in {\cal B}^*}{\cal A}^*_B$ forms the set of base blocks of the required design. \qed

It is well known that a $(k,v)$-CDM does not exist for any $v\equiv 0\ ({\rm mod}\ 2)$ (see \cite{c}). Thus when $v$ is even, Construction \ref{SCHGDD from CDM} does not work. We shall present Constructions \ref{SCHGDD from CHDM} and \ref{SCHGDD from perfect} to deal with this problem.

\begin{Construction}\label{SCHGDD from CHDM}
Let gcd$(wt,hg)=1$. Suppose that there exists a $k$-SCHGDD of type $(n,w^t)$. If there exist a $k$-SCHGDD of type $(n,(hw)^t)$ and a $(k,hg;h)$-CHDM, then there exists a $k$-SCHGDD of type $(n,(hgw)^{t})$.
\end{Construction}

\proof Let $S=\{0,t,\ldots,(w-1)t\}$ be a subgroup of order $w$ in $Z_{wt}$, and $S_l=S+l$ be a coset of $S$ in $Z_{wt}$, $0\leq l\leq t-1$. By assumption we can construct a $k$-SCHGDD of type $(n,w^t)$ on $I_n\times Z_{wt}$ with the group set $\{\{i\}\times Z_{wt}:i\in I_n\}$ and the hole set $\{I_n\times S_l:0\leq l\leq t-1\}$. Denote the set of its base blocks by ${\cal B}^*$.

Let $D=(d_{ij})$ be the given $(k,hg;h)$-CHDM. Then for its two distinct rows $x$ and $y$, the difference list $L_{xy}=\{d_{xj}-d_{yj}:j\in I_{h(g-1)}\}$ contains each integer of $Z_{hg}\setminus H$ exactly once, where $H=\{0,g,\ldots,(h-1)g\}$ is a subgroup of order $h$ in $Z_{hg}$.

Let $X=I_n\times (Z_{wt}\times Z_{hg})$. Since gcd$(wt,hg)=1$, $X$ is isomorphic to $I_n\times Z_{wthg}$. We shall construct the required $k$-SCHGDD of type $(n,(hgw)^{t})$ on $X$ with the group set $\{\{i\}\times (Z_{wt}\times Z_{hg}):i\in I_n\}$ and the hole set $\{I_n\times (S_l\times Z_{hg}):0\leq l\leq t-1\}$ as follows (note that due to gcd$(w,hg)=1$, $S_0\times Z_{hg}$ is isomorphic to $Z_{whg}$):

First for each base block $B=\{(a_1,x_1),(a_2,x_2),\ldots,(a_k,x_k)\}\in{\cal B}^*$, we construct a family, ${\cal A}^*_B$, consisting of the following $h(g-1)$ base blocks
$$\{(a_1,x_1,d_{1j}),(a_2,x_2,d_{2j}),\ldots,(a_k,x_k,d_{kj})\},$$
where $0\leq j\leq h(g-1)-1$.

Next since gcd$(wt,hg)=1$, by assumption we can construct a $k$-SCHGDD of type $(n,(hw)^t)$ on $I_n\times (Z_{wt}\times H)$ with the group set $\{\{i\}\times (Z_{wt}\times H):i\in I_n\}$ and the hole set $\{I_n\times (S_l\times H):0\leq l\leq t-1\}$. Denote the set of its base blocks by ${\cal C}$. It is readily checked that $(\bigcup_{B\in {\cal B}^*}{\cal A}^*_B)\cup{\cal C}$ forms the set of base blocks of the required design. \qed

Construction \ref{SCHGDD from CHDM} has a requirement that gcd$(wt,hg)=1$. To relax this condition, we introduce a special kind of SCHGDD called $h$-perfect SCHGDD as follows, which is an analogy of an $h$-perfect difference family.

Suppose that ${\cal B}^*=\{B_i:i=1,2,\ldots,r\}$ is the family of base blocks of a $k$-SCHGDD of type $(n,(hg)^t)$ on $I_n\times Z_{hgt}$, where $B_i=\{(a_{1i},b_{1i}),(a_{2i},b_{2i}),\ldots,(a_{ki},b_{ki})\}$. Take $\delta_i={\rm min}\{b_{1i},b_{2i},\ldots,b_{ki}\}$ for each possible $i$. Then we can use $B_i-\delta=\{(a_{1i},b_{1i}-\delta),(a_{2i},b_{2i}-\delta),\ldots,(a_{ki},b_{ki}-\delta)\}$ instead of $B_i$. So without loss of generality we can assume that each base block of ${\cal B}^*$ is of the form $\{(a_{1i},b_{1i}),(a_{2i},b_{2i}),\ldots,(a_{l_i-1,i},b_{l_i-1,i}),(a_{l_i,i},0),
(a_{l_i+1,i},b_{l_i+1,i}),\ldots,$ $
(a_{ki},b_{ki})\},$
where $l_i$ is some value from $1$ to $k$ and $b_{ji}\not\equiv 0\ ({\rm mod}\ t)$ for each $1\leq i\leq r$, $1\leq j\leq k$ and $j\neq l_i$. Define a set
$$ele({\cal B}^*)=\bigcup_{i=1}^r \{b_{1i},b_{2i},\ldots,b_{l_i-1,i},b_{l_i+1,i},\ldots,b_{ki}\}.$$
The $k$-SCHGDD of type $(n,(hg)^t)$ is said to be {\em $h$-perfect} if
$$ele({\cal B}^*)\subseteq\{x+ygt:\ 0\leq x\leq
\lfloor gt/2\rfloor,x\not\equiv 0\ ({\rm mod}\ t),
0\leq y\leq h-1\}.$$
When $h=1$, a $1$-perfect $k$-SCHGDD is simply called a {\em perfect} $k$-SCHGDD.

\begin{Construction}\label{SCHGDD from perfect}
Suppose that there exists a perfect $k$-SCHGDD of type $(n,w^t)$. If there exist an $h$-perfect $k$-SCHGDD of type $(n,(hw)^t)$ and a $(k,hg;h)$-CHDM, then there exists an $hg$-perfect $k$-SCHGDD of type $(n,(hgw)^{t})$.
\end{Construction}

\proof Suppose that ${\cal A}^*=\{\{(a_{1i},b_{1i}),(a_{2i},b_{2i}),\ldots,(a_{l_i,i},0),\ldots,(a_{ki},b_{ki})\}:1\leq i\leq r\}$ is the set of base blocks of the given perfect $k$-SCHGDD of type $(n,w^t)$, where $l_i$ is some value from $1$ to $k$, $0\leq b_{ji}\leq \lfloor wt/2\rfloor$ and $b_{ji}\not\equiv 0\ ({\rm mod}\ t)$ for each $1\leq i\leq r$, $1\leq j\leq k$ and $j\neq l_i$.

Let ${\cal B}^*=\{\{(c_{1i},x_{1i}+y_{1i}wt),(c_{2i},x_{2i}+y_{2i}wt),\ldots,(c_{l'_i,i},0),\ldots,
(c_{ki},x_{ki}+y_{ki}wt)\}:1\leq i\leq s\}$ be the set of base blocks of the given $h$-perfect $k$-SCHGDD of type $(n,(hw)^t)$, where $l'_i$ is some value from $1$ to $k$, $0\leq x_{ji}\leq
\lfloor wt/2\rfloor$, $x_{ji}\not\equiv 0\ ({\rm mod}\ t)$ and $0\leq y_{ji}\leq h-1$ for each $1\leq i\leq s$, $1\leq j\leq k$ and $j\neq l'_i$.

Let $D=(d_{ij})$ be the given $(k,hg;h)$-CHDM, where $d_{ij}\in Z_{hg}$ for $1\leq i\leq k$ and $1\leq j\leq h(g-1)$ such that for any two distinct rows $\alpha$ and $\beta$, $\{d_{\alpha j}-d_{\beta j}:1\leq j\leq h(g-1)\}$ contains each integer of $Z_{hg}\setminus \{0,g,\ldots,(h-1)g\}$.

Now we construct the desired $hg$-perfect $k$-SCHGDD of type $(n,(hgw)^{t})$ on $I_n\times Z_{hgwt}$, whose base blocks consists of the following two parts:

\vspace{0.1cm}
$(1)$ For each $A_i=\{(a_{1i},b_{1i}),(a_{2i},b_{2i}),\ldots,(a_{l_i,i},0),\ldots,(a_{ki},b_{ki})\}\in {\cal A}^*$, we construct a family, ${\cal C}^*_{A_i}$, consisting of the following $h(g-1)$ base blocks
$$\{(a_{1i},b_{1i}+(d_{2\gamma}-d_{1\gamma})wt),(a_{2i},b_{2i}+(d_{3\gamma}-
d_{1\gamma})wt),\ldots,(a_{l_i,i},0),\ldots,(a_{ki},b_{ki}+(d_{k\gamma}-d_{1\gamma})wt)\},$$
where $1\leq \gamma\leq h(g-1)$ and the second coordinates are reduced modulo $hgwt$. Let ${\cal C}^*=\bigcup_{i=1}^r {\cal C}^*_{A_i}$. For $1\leq \rho\neq\theta\leq n$, define the notation $\Delta_{\rho\theta}({\cal C}^*)=\bigcup_{C\in{\cal C}^*}\Delta_{\rho\theta}(C)$, where $\Delta_{\rho\theta}(C)=\{x-y\ ({\rm mod}\ hgwt): (\rho,x),(\theta,y)\in C\}$. Then, by noting that $-\lfloor wt/2\rfloor\leq b_{ei}-b_{fi}\leq \lfloor wt/2\rfloor$ for each $1\leq i\leq r$ and each admissible $1\leq e\neq f\leq k$, it is readily checked that for any $1\leq \rho\neq\theta\leq n$,
$$\Delta_{\rho\theta}({\cal C}^*)=\pm\{p+qwt:0\leq p\leq \lfloor wt/2\rfloor,p\not\equiv 0\ ({\rm mod}\ t),q\in Z_{hg}\setminus \{0,g,\ldots,(h-1)g\}\}.$$

$(2)$ For each $B_i=\{(c_{1i},x_{1i}+y_{1i}wt),(c_{2i},x_{2i}+y_{2i}wt),\ldots,(c_{l'_i,i},0),\ldots,
(c_{ki},x_{ki}+y_{ki}wt)\}\in{\cal B}^*$, we take a base block
$$\bar{B}_i=\{(c_{1i},x_{1i}+y_{1i}wtg),(c_{2i},x_{2i}+y_{2i}wtg),\ldots,(c_{l'_i,i},0),\ldots,
(c_{ki},x_{ki}+y_{ki}wtg)\}.$$
Let ${\cal D}^*=\{\bar{B}_i:1\leq i\leq s\}$. Then similarly it can be checked that for any $1\leq \rho\neq\theta\leq n$,
$$\Delta_{\rho\theta}({\cal D}^*)=\pm\{p+qwt:0\leq p\leq \lfloor wt/2\rfloor,p\not\equiv 0\ ({\rm mod}\ t),q\in \{0,g,\ldots,(h-1)g\}\}.$$
Thus ${\cal C}^*\cup{\cal D}^*$ forms the set of base blocks of the required $hg$-perfect $k$-SCHGDD of type $(n,(hgw)^{t})$. \qed

To use Construction \ref{SCHGDD from perfect}, we need to find some perfect $k$-SCHGDDs of type $(n,w^t)$. So we establish the following construction, whose proof is similar to that of Construction \ref{SCHGDD-from strictly CGDD} (here only a routine check for the perfect property is needed).

\begin{Construction} \label{perfect SCHGDD-from perfect DF}
Suppose that there exist a $(hwt,hw,K,1)$-$h$-PDF and an $l$-MGDD of type $k^n$ for each $k\in K$. Then there exists a $h$-perfect $l$-SCHGDD of type $(n,(hw)^t)$.
\end{Construction}

\section{Some cyclic or $h$-perfect difference families}

To apply Constructions \ref{SCHGDD-from strictly CGDD} and \ref{perfect SCHGDD-from perfect DF}, we need some cyclic or $h$-perfect difference families. Actually cyclic difference families have their own rich design-theoretical content (see \cite{ab}) and many application backgrounds (for example \cite{gfm,glm}). In this section, we always assume that $[a,b]$ denotes the set of integers $n$ such that $a\leq n\leq b$, and $[a,b]_{e}$ denotes the set of even integers in $[a,b]$.

Chen \cite{czl} investigated the existences of a $(v,\{3,4^*\},1)$-PDF with exactly one block of size $4$ and a $(v,\{3,4^*,5^*\},1)$-PDF with exactly one block of size $4$ and exactly one block of size $5$. The main method he used to obtain these PDFs is by means of the following Lemma \ref{recursive for PDF}, which is proved by using extended Skolem sequences.

\begin{Lemma}{\rm \cite{czl}}\label{recursive for PDF} Suppose that there is a $(v,K,1)$-PDF. Let $u\geq v$ be a positive integer. Then there exists a $(6u+v,K\cup\{3\},1)$-PDF for any $u\equiv 0,1\ ({\rm mod}\ 4)$ and $v\equiv 1\ ({\rm mod}\ 4)$, or $u\equiv 0,3\ ({\rm mod}\ 4)$ and $v\equiv 3\ ({\rm mod}\ 4)$.
\end{Lemma}

\begin{Lemma}
{\rm \cite{czl}}\label{PDF-24-19}
\begin{enumerate}
\item[$(1)$] There exists a $(v,\{3,4^*\},1)$-PDF for $v\equiv 1\ ({\rm mod}\ 6)$ and $v\geq 19$.
\item[$(2)$] There exists a $(v,\{3,4^*,5^*\},1)$-PDF for $v\equiv 3\ ({\rm mod}\ 6)$ and $v\geq 39$.
\end{enumerate}
\end{Lemma}

\begin{Lemma}\label{PDF-24-5} There exists a $(v,\{3,4,5\},1)$-PDF for $v\equiv 5\ ({\rm mod}\ 6)$ and $v\geq 59$.
\end{Lemma}

\proof When $v\in\{59,71\}$, the conclusion follows from Example \ref{example-PDF}. When $v\in\{65,77,89$, $101,113,125,137,149\}$, a $(v,\{4,5\},1)$-PDF can be found in Appendix B of \cite{glm}. When $v\equiv 5\ ({\rm mod}\ 6)$, $77< v<455 $ and $v\not\in \{89,101,113,125,137,149,413,419,437,443\}$, we have given a direct construction for a $(v,\{3,4,5\},1)$-PDF by computer, which are shown in Appendix A.

When $v\equiv 5,11\ ({\rm mod}\ 24)$ and $v\geq 413$, assume that $v=6u+59$ with $u\equiv 0,3\ ({\rm mod}\ 4)$ and $u\geq 59$. Start from a $(59,\{3,4,5\},1)$-PDF, which exists by the first paragraph of this proof. Then apply Lemma \ref{recursive for PDF} to obtain the required $(6u+59,\{3,4,5\},1)$-PDF.

When $v\equiv 17,23\ ({\rm mod}\ 24)$ and $v\geq 455$, assume that $v=6u+65$ with $u\equiv 0,1\ ({\rm mod}\ 4)$ and $u\geq 65$. Start from a $(65,\{4,5\},1)$-PDF, which exists by the first paragraph of this proof. Then apply Lemma \ref{recursive for PDF} to obtain the required $(6u+65,\{3,4,5\},1)$-PDF. \qed

\begin{Lemma}\label{relation PDF and PDF} If there exists a $(2t-1,K,1)$-PDF, then there exists a $(2t,2,K,1)$-PDF.
\end{Lemma}
\proof Suppose there is a $(2t-1,K,1)$-PDF ${\cal A}$. Then $\Delta^+{\cal A}=\bigcup_{A\in{\cal A}}\{x-y:x,y\in A,x>y\}=[1,t-1]$. Obviously $\cal A$ is also a $(2t,2,K,1)$-PDF. \qed

\begin{Lemma}
\label{2t,P} There exists a $(2t,2,\{3,4,5\},1)$-PDF for each integer $t\geq 7$ and $t\not\in\{8,9,11$, $12,14,15,18,21,24,27\}$.
\end{Lemma}

\proof By Lemma \ref{relation PDF and PDF}, it suffices to construct the corresponding $(2t-1,\{3,4,5\},1)$-PDFs. By Example \ref{example-PDF} and Lemma \ref{PDF-24-19}$(1)$, when $t\equiv 1\ ({\rm mod}\ 3)$ and $t\geq 7$, there is a $(2t-1,\{3,4\},1)$-PDF. By Example \ref{example-PDF} and Lemma \ref{PDF-24-19}$(2)$, when $t\equiv 2\ ({\rm mod}\ 3)$ and $t\geq 17$, there is a $(2t-1,\{3,4,5\},1)$-PDF. By Lemma \ref{PDF-24-5}, when $t\equiv 0\ ({\rm mod}\ 3)$ and $t\geq 30$, a $(2t-1,\{3,4,5\},1)$-PDF exists.  \qed

\begin{Lemma}\label{langford}{\rm\cite{linek}} If $u\geq2d-1$ and $u\not\in[2d+2,8d-5]$, then the set $[1,2u+1]\setminus\{c\}$
can be partitioned into some pairs $\{x_i,y_i\}$, $1\leq i\leq u$, such that $\{y_i-x_i:1\leq i\leq u\}=[d,d+u-1]$ whenever $(u,c)\equiv (0,1),(1,d),(2,0),(3,d+1)\ ({\rm mod}\ (4,2))$.
\end{Lemma}

\begin{Lemma}\label{holey pdf} If there exists a $(v,K,1)$-PDF, then there exists a $(v+6u+3,4,K\cup\{3\},1)$-PDF for any
 $u\geq v$, $u\not\in[v+3,4v-1]$, $2u-v\equiv 3\ ({\rm mod}\ 4)$, and $(u,(2u-v+5)/4)\equiv (0,1),(1,(v+1)/2),(2,0),(3,(v+3)/2)\ ({\rm mod}\ (4,2))$.
\end{Lemma}

\proof  Suppose there is a $(v,K,1)$-PDF ${\cal A}$, then $\Delta^+{\cal A}=\bigcup_{A\in{\cal A}}\{x-y:x,y\in A,x>y\}=[1,(v-1)/2]$. Now apply Lemma \ref{langford} with $c=(2u-v+5)/4$ and $d=(v+1)/2$ (note that $(v+1)/2$ is an integer since the existence of a $(v,K,1)$-PDF implies $v$ is odd). We have that the set $[1,2u+1]\setminus\{(2u-v+5)/4\}$ can be partitioned into some pairs $\{x_i,y_i\}$, $1\leq i\leq u$, such that $\{y_i-x_i:1\leq i\leq u\}=[(v+1)/2,(v-1)/2+u]$.

Let $B_i=\{0,x_i+(v-1)/2+u,y_i+(v-1)/2+u\}$ and ${\cal B}=\bigcup_{1\leq i\leq u}B_i$. It is readily checked that $\Delta^+{\cal B}=\bigcup_{B\in{\cal B}}\{x-y:x,y\in B,x>y\}=[(v+1)/2,(v+1)/2+3u]\setminus\{(v+6u+3)/4\}$. So ${\cal A}\cup{\cal B}$ forms the desired $(v+6u+3,4,K\cup\{3\},1)$-PDF. \qed

\begin{Remark}\label{remark gap} Lemma $\ref{holey pdf}$ does not work when $u\in[v+3,4v-1]$. However, from the proof of Lemma $\ref{holey pdf}$, we know that if one can remove the condition $u\not\in[2d+2,8d-5]$ in Lemma $\ref{langford}$, then Lemma $\ref{holey pdf}$ also holds for $u\in[v+3,4v-1]$.
\end{Remark}

\begin{Lemma}\label{4t,4-3,5-1} There exists a $(4t,4,\{3,5\},1)$-PDF for $t\equiv 0\ ({\rm mod}\ 2)$, $t\geq 4$ and $t\neq8$.
\end{Lemma}

\proof $(1)$ When $t\equiv 4\ ({\rm mod}\ 6)$ and $t\geq 4$, the conclusion follows from Lemma $3.5$ in \cite{wc}.

$(2)$ When $t\equiv 0\ ({\rm mod}\ 6)$, $t\geq 210$ or $t=60$, start from a $(33,\{3,5\},1)$-PDF, which exists by Example \ref{example-PDF}. Apply Lemma \ref{holey pdf} to obtain a $(36+6u,4,\{3,5\},1)$-PDF with $u\equiv 2\ ({\rm mod}\ 4)$, $u\geq 33$ and  $u\not\in[36,131]$. Write $t=9+3u/2$. Then a $(4t,4,\{3,5\},1)$-PDF exists.

For $t\equiv 0\ ({\rm mod}\ 6)$, $6\leq t\leq 54$, we have given a direct construction for a $(4t,4,\{3,5\},1)$-PDF by computer, which are shown in Appendix B. For $t\equiv 0\ ({\rm mod}\ 6)$ and $66\leq t\leq204$, by Remark \ref{remark gap}, if we can show that for each $u\in[36,131]$ and $u\equiv 2\ ({\rm mod}\ 4)$, the set $[1,2u+1]\setminus\{(u-14)/2\}$
can be partitioned into some pairs $\{x_i,y_i\}$, $1\leq i\leq u$, such that $\{y_i-x_i:1\leq i\leq u\}=[17,16+u]$, then take $t=9+3u/2$ and a $(4t,4,\{3,5\},1)$-PDF can be obtained. By computer, we have found all these required partitions, which are shown in Appendix D.

$(3)$ Similarly, when $t\equiv 2\ ({\rm mod}\ 6)$, $t\geq 410$ or $t=116$, start from a $(65,\{3,5\},1)$-PDF, which exists by Lemma \ref{PDF-24-5}. Apply Lemma \ref{holey pdf} to obtain a $(68+6u,4,\{3,5\},1)$-PDF with $u\equiv 2\ ({\rm mod}\ 4)$, $u\geq 65$ and  $u\not\in[68,259]$. Write $t=17+3u/2$. Then a $(4t,4,\{3,5\},1)$-PDF exists.

For $t\equiv 2\ ({\rm mod}\ 6)$, $14\leq t\leq 110$, we have given a direct construction for a $(4t,4,\{3,5\},1)$-PDF by computer, which are shown in Appendix C. For $t\equiv 2\ ({\rm mod}\ 6)$ and $122\leq t\leq 404$, by Remark \ref{remark gap}, if we can show that for each $u\in[68,259]$ and $u\equiv 2\ ({\rm mod}\ 4)$, the set $[1,2u+1]\setminus\{(u-30)/2\}$
can be partitioned into some pairs $\{x_i,y_i\}$, $1\leq i\leq u$, such that $\{y_i-x_i:1\leq i\leq u\}=[33,32+u]$, then take $t=17+3u/2$ and a $(4t,4,\{3,5\},1)$-PDF can be obtained. By computer, we have found all these required partitions, which are shown in Appendix E. \qed

\begin{Corollary}
\label{4t,4,PSCHGDD} There exists a perfect $3$-SCHGDD of type $(3,4^t)$ for $t\equiv 0\ ({\rm mod}\ 2)$, $t\geq 4$ and $t\neq8$.
\end{Corollary}

\proof By Lemma \ref{4t,4-3,5-1}, there is a $(4t,4,\{3,5\},1)$-PDF for $t\equiv 0\ ({\rm mod}\ 2)$, $t\geq 4$ and $t\neq8$. Then apply Construction \ref{perfect SCHGDD-from perfect DF} with $h=1$ and $w=4$ to obtain a perfect $3$-SCHGDD of type $(3,4^t)$, where the needed $3$-MGDDs of type $k^3$, $k\in\{3,5\}$, are from Theorem \ref{3-HGDD}. \qed

The following two lemmas are from the use of extended Langford sequences.

\begin{Lemma}
\label{se1}{\rm \cite[Lemma 1.7]{zhang} } Let $(u,d)\equiv
(0,1),(1,1),(0,0),(3,0)\ ({\rm mod}\ (4,2))$ such that $u\geq 2d-1.$
Then $[d,d+3u-1]$ can be partitioned into triples $\{a_i,b_i,c_i\},
1\leq i\leq u,$ such that $a_i+b_i=c_i.$
\end{Lemma}

\begin{Lemma}
\label{se}{\rm \cite[Lemma 1.10]{zhang} } For $1\leq d\leq4,$ if
$(u,k)\equiv (0,1),(1,d),(2,0),(3,d+1)\ ({\rm mod}\ (4,\linebreak2))$ such
that $u\geq 2d-3$ and $(u/2)(2d-1-u)+1\leq k\leq (u/2)(u-2d+5)+1,$
then $[d,d+3u]\backslash\{k+d+u-1\}$ can be partitioned into triples
$\{a_i,b_i,c_i\}, 1\leq i\leq u,$ such that $a_i+b_i=c_i.$
\end{Lemma}

\begin{Lemma}
\label{8t5} There exists a $(8t,8,\{3,5\},1)$-$2$-PDF for $t\equiv 0\ ({\rm mod}\ 2)$ and $t\geq 8$.
\end{Lemma}

\proof $(1)$ When $t\equiv 4\ ({\rm mod}\ 6)$ and $t\geq 10$, the conclusion follows immediately from Lemma $3.6$ in \cite{wc}.

$(2)$ When $t=12$, a $(8t,8,\{3,5\},1)$-$2$-PDF is listed below:
\begin{center}
\begin{tabular}{lllllll}
$\{0,1,20\},$ & $\{0,3,62\},$ & $\{0,5,69\},$ & $\{0,17,70\},$ & $\{0,13,22,52,68\},$\\
$\{0,2,10\},$ & $\{0,4,71\},$ & $\{0,6,21\},$ & $\{0,23,58\},$ & $\{0,11,18,51,65\}.$\\
\end{tabular}
\end{center}
When $t=18$, a $(8t,8,\{3,5\},1)$-$2$-PDF is listed below:
\begin{center}
\begin{tabular}{lllllll}
$\{0,1,21\},$ & $\{0,4,27\},$ & $\{0,7,88\},$ & $\{0,10,103\},$ & $\{0,14,105\},$ & $\{0,19,34,76,104\},$\\
$\{0,2,24\},$ & $\{0,5,82\},$ & $\{0,8,100\},$ & $\{0,11,89\},$ & $\{0,16,95\},$ & $\{0,17,30,75,101\},$\\
$\{0,3,32\},$ & $\{0,6,31\},$ & $\{0,9,107\},$ & $\{0,12,106\},$ & $\{0,33,80\},$ & $\{0,35,83\}.$\\
\end{tabular}
\end{center}
When $t\equiv 0\ ({\rm mod}\ 6)$ and $t\geq 24$, take $A_1=\{0,4t+4,t+1,2t-2,6t-4\}$ and $A_2=\{0,4t+3,t-1,6t-7,2t-6\}$ as two base blocks with block size five. Let $S=([1,2t-1]\cup[4t+1,6t-1])\setminus \{t,5t\}$ and $T=\{t-5,t-3,t-1,t+1,2t-10,2t-8,2t-6,2t-2,4t+1,4t+2,4t+3,4t+4,5t-6,5t-5,5t-4,5t-3,6t-9,6t-7,
6t-6,6t-4\}$.

By Lemma \ref{se} with $(d,u,k)=(3,t/3-3,t/6+1)$, the set $[3,t-6]\setminus \{t/2\}$ can be partitioned into triples $\{a'_i,b'_i,c'_i\},$ $1\leq i\leq t/3-3,$ such that $a'_i+b'_i=c'_i.$ Taking $a_i=2a'_i,$ $b_i=2b'_i$ and $c_i=2c'_i$ for $1\leq i\leq t/3-3,$ we have that $[6,2t-12]_{e}\setminus \{t\}$ can be partitioned into triples $\{a_i,b_i,c_i\},$ $1\leq i\leq t/3-3,$ such that $a_i+b_i=c_i.$

Furthermore, $(S\setminus T)\setminus ([6,2t-12]_{e}\setminus \{t\})$ can be partitioned into triples $\{a_i,b_i,c_i\},$ $t/3-2\leq
i\leq 4t/3-8$, such that $a_i+b_i\equiv c_i\ ({\rm mod}\ 8t)$ as follows:

\begin{center}
\begin{tabular}{ll}
$\{13+2j,5t-9-j,5t+4+j\},$ & $j\in[0,t/2-11];$\\
$\{t+3+2j,9t/2-3-j,11t/2+j\},$ & $j\in[0,t/2-10];$
\end{tabular}

\begin{tabular}{lll}
$\{2,9,11\},$ & $\{t-7,5t+2,6t-5\},$ & $\{9t/2-2,11t/2-5,2t-7\},$\\
$\{4,2t-15,2t-11\},$ & $\{5,11t/2-6,11t/2-1\},$ & $\{9t/2-1,11t/2-4,2t-5\},$\\
$\{3,5t-2,5t+1\},$ & $\{2t-13,4t+5,6t-8\},$ & $\{9t/2,11t/2-3,2t-3\},$\\
$\{7,5t-8,5t-1\},$ & $\{2t-9,4t+6,6t-3\},$ & $\{9t/2+1,11t/2-2,2t-1\},$\\
$\{1,6t-2,6t-1\},$ & $\{5t-7,5t+3,2t-4\}.$ &  \\
\end{tabular}
\end{center}
Then $\{\{0,a_i,c_i\}:1\leq i\leq 4t/3-8\}\cup\{A_1,A_2\}$ forms a $(8t,8,\{3,5\},1)$-$2$-PDF.

$(3)$ When $t=8$, a $(8t,8,\{3,5\},1)$-$2$-PDF is listed below:
\begin{center}{\tabcolsep 0.05in
\begin{tabular}{lllllll}
$\{0,1,3,7,44\},$ & $\{0,5,14\},$ & $\{0,10,45\},$  & $\{0,11,47\},$  & $\{0,12,34\},$  & $\{0,13,38\},$  & $\{0,15,33\}.$\\
\end{tabular}}
\end{center}
When $t=14$, a $(8t,8,\{3,5\},1)$-$2$-PDF is listed below:
\begin{center}
\begin{tabular}{lllllll}
 $\{0,1,10\},$ & $\{0,4,12\},$ &  $\{0,7,81\},$ & $\{0,21,82\},$ & $\{0,15,26,58,83\},$ &
 $\{0,2,18\},$ \\ $\{0,5,71\},$ & $\{0,13,78\},$ & $\{0,22,59\},$ & $\{0,24,60\},$& $\{0,3,20\},$ & $\{0,6,79\},$
  \\ $\{0,19,64\},$ & $\{0,23,63\},$ & $\{0,27,62\}.$\\
\end{tabular}
\end{center}
\vspace{0.1cm}
When $t\equiv 2\ ({\rm mod}\ 6)$ and $t\geq20$, take $A=\{0,4t+2,t+1,2t-2,6t-1\}$ as the base block with block size five. Let $S=([1,2t-1]\cup[4t+1,6t-1])\setminus \{t,5t\}$ and $T=\{t-3,t+1,2t-3,2t-2,4t+1,4t+2,5t-2,5t-1,6t-4,6t-1\}$.

By Lemma \ref{se} with $(d,u,k)=(3,(t-5)/3,(t-2)/6)$, the set $[3,t-2]\setminus \{t/2\}$ can be partitioned into triples $\{a'_i,b'_i,c'_i\},$ $1\leq i\leq (t-5)/3,$ such that $a'_i+b'_i=c'_i.$ Taking $a_i=2a'_i,$ $b_i=2b'_i$ and $c_i=2c'_i$ for $1\leq i\leq (t-5)/3$, we have that $[6,2t-4]_{e}\setminus \{t\}$ can be partitioned into triples $\{a_i,b_i,c_i\},$ $1\leq i\leq (t-5)/3,$ such that $a_i+b_i=c_i.$

Furthermore, $(S\setminus T)\setminus ([6,2t-4]_{e}\setminus \{t\})$ can be partitioned into triples $\{a_i,b_i,c_i\},$ $(t-2)/3\leq
i\leq (4t-14)/3$, such that $a_i+b_i\equiv c_i\ ({\rm mod}\ 8t)$ as follows:
\begin{center}
\begin{tabular}{ll}
$\{11+2j,5t-6-j,5t+5+j\},$ & $j\in[0,t/2-9];$\\
$\{t+3+2j,9t/2-3-j,11t/2+j\},$ & $j\in[0,t/2-6];$\\
\end{tabular}

\begin{tabular}{lll}
$\{4,5,9\},$ & $\{t-5,5t+3,6t-2\},$ & $\{5t-3,5t-4,2t-7\},$\\
$\{3,5t+1,5t+4\},$ & $\{t-1,9t/2-1,11t/2-2\},$ & $\{9t/2-2,11t/2-3,2t-5\},$\\
$\{7,5t-5,5t+2\},$ & $\{1,9t/2+1,9t/2+2\},$ & $\{9t/2,11t/2-1,2t-1\},$\\
$\{2,6t-5,6t-3\}.$ &  & \\
\end{tabular}
\end{center}
Then $\{\{0,a_i,c_i\}:1\leq i\leq (4t-14)/3\}\cup\{A\}$ forms a $(8t,8,\{3,5\},1)$-$2$-PDF. \qed

\begin{Corollary}
\label{8t,8,PSCHGDD} There exists a $2$-perfect $3$-SCHGDD of type $(3,8^t)$ for $t\equiv 0\ ({\rm mod}\ 2)$ and $t\geq 8$.
\end{Corollary}

\proof By Lemma \ref{8t5}, there is a $(8t,8,\{3,5\},1)$-$2$-PDF for $t\equiv 0\ ({\rm mod}\ 2)$ and $t\geq 8$. Then apply Construction \ref{perfect SCHGDD-from perfect DF} with $h=2$ and $w=4$ to obtain a $2$-perfect $3$-SCHGDD of type $(3,8^t)$, where the needed $3$-MGDDs of type $k^3$, $k\in\{3,5\}$, are from Theorem \ref{3-HGDD}. \qed

With similar methods to those in Lemma \ref{8t5}, making use of Lemmas \ref{se1} and \ref{se} thoroughly, we can construct a $(4t,4,\{3,5\},1)$-CDF for any $t\equiv 1\ ({\rm mod}\ 2)$, $t\geq 7$ and $t\neq11$, and a $(16t,16,\{3,5\},1)$-CDF for any $t\equiv 0\ ({\rm mod}\ 2)$ and $t\geq 4$. We have moved these proofs to Appendices F and G, respectively.

\begin{Lemma}
\label{4t,2p} There exists a $(4t,4,\{3,5\},1)$-CDF for $t\equiv 1\ ({\rm mod}\ 2)$, $t\geq 7$ and $t\neq11$.
\end{Lemma}

\begin{Lemma}
\label{16t,CDF} There exists a $(16t,16,\{3,5\},1)$-CDF for $t\equiv 0\ ({\rm mod}\ 2)$ and $t\geq 4$.
\end{Lemma}

\section{Cyclic holey difference matrices}

In this section we shall establish the necessary and sufficient condition for the existence of a $(3,wt;w)$-CHDM, which is equivalent to a $3$-SCHGDD of type $(3,w^t)$ by Lemma \ref{relation SCH and CHDM}.

\begin{Lemma}{\rm \cite{yin2005}}\label{yin-CHDM}
Let $q=2^{\alpha}3^{\beta}p_1^{\alpha_1}p_2^{\alpha_2}\cdots p_s^{\alpha_s}$ be the prime factorization of $q\geq 4$. If $(\alpha,\beta)\neq (1,0),(0,1)$, then there is a $(4,2q;2)$-CHDM.
\end{Lemma}

Let $d$ and $v$ be integers with $1\leq d\leq 2v+1$. A {\em $d$-extended Skolem sequence} of order $v$ is a sequence $(a_1,a_2,\ldots,a_v)$ of $v$ integers satisfying $\bigcup_{i=1}^v\{a_i,a_i-i\}=\{1,2,\ldots,2v+1\}\setminus\{d\}.$
When $d=2v+1$, it is simply called a {\em Skolem sequence} of order $v$. For more details on Skolem sequences, the reader may refer to \cite{shalaby}.

\begin{Lemma}{\rm \cite{baker}}\label{extended-skolem}
A $d$-extended Skolem sequence of order $v$ exists if and only if
\emph{(i)} $d$ is odd and $v\equiv 0,1\ ({\rm mod}\ 4)$, or
\emph{(ii)} $d$ is even and $v\equiv 2,3\ ({\rm mod}\ 4)$.
\end{Lemma}

\begin{Lemma}\label{3-chdm-2} There exists a $(3,2t;2)$-CHDM for any integer $t\geq 4$.
\end{Lemma}

\proof When $t\equiv 0,1\ ({\rm mod}\ 4)$ and $t\geq 4$, a $(3,2t;2)$-CHDM is constructed as follows:
$$(A_1|A_2|\cdots|A_{t-1}|B_1|B_2|\cdots|B_{t-1}),$$
where $A_i=(0,i,a_i)^{\rm T}$ and $B_i=(0,-i,a_i-i)^{\rm T}$, $1\leq i\leq t-1$, such that $(a_1,a_2,\ldots,a_{t-1})$ is an $t$-extended Skolem sequence of
order $t-1$, which exists by Lemma \ref{extended-skolem}.

When $t\in\{6,11,18,27\}$, by Lemma \ref{yin-CHDM} there is a $(4,2t;2)$-CHDM, and deleting any row of this CHDM yields a $(3,2t;2)$-CHDM. When $t\in\{14,15\}$, we list a $(3,2t;2)$-CHDM as follows:
\begin{center}
{\footnotesize\tabcolsep 0.043in
\noindent {\rm $\left(
         \begin{tabular}{cccccccccccccccccccccccccc}
           0 & 0 & 0 & 0 & 0 & 0 & 0 & 0 & 0 & 0 & 0 & 0 & 0 & 0 & 0 & 0 & 0 & 0 & 0 & 0 & 0 & 0 & 0 & 0 & 0 & 0 \\
           1 & 2 & 3 & 4 & 5 & 6 & 7 & 8 & 9 & 10 & 11 & 12 & 13 & 15 & 16 & 17 & 18 & 19 & 20 & 21 & 22 & 23 & 24 & 25 & 26 & 27 \\
           16 & 1 & 4 &2 & 26 & 3 & 10 & 13 & 15 & 20 & 22 & 24 & 8 & 23 & 25 & 7 & 9 & 21 & 5 & 17 & 11 & 27 & 12 & 19 & 18 & 6 \\
         \end{tabular}
       \right);
$}}
\end{center}

\begin{center}
{\footnotesize\tabcolsep 0.043in
\noindent {\rm $\left(
         \begin{tabular}{cccccccccccccccccccccccccccc}
           0 & 0 & 0 & 0 & 0 & 0 & 0 & 0 & 0 & 0 & 0 & 0 & 0 & 0 & 0 & 0 & 0 & 0 & 0 & 0 & 0 & 0 & 0 & 0 & 0 & 0 & 0 & 0\\
           1 & 2 & 3 & 4 & 5 & 6 & 7 & 8 & 9 & 10 & 11 & 12 & 13 & 14 & 16 & 17 & 18 & 19 & 20 & 21 & 22 & 23 & 24 & 25 & 26 & 27 & 28 & 29 \\
           2 & 4 & 6 & 8 & 10 & 12 & 14 & 16 & 18 & 20 & 23 & 28 & 27 & 1 & 29 & 5 & 7 & 9 & 11 & 13 & 3 & 17 & 19 & 21 & 24 & 26 & 25 & 22 \\
         \end{tabular}
       \right).
$}}
\end{center}
When $t\equiv 2,3\ ({\rm mod}\ 4)$, $t\geq 7$ and $t\not\in\{11,14,15,18,27\}$, by Lemma \ref{2t,P}, there is a $(2t,2$, $\{3,4,5\},1)$-PDF, which implies a strictly cyclic $\{3,4,5\}$-GDD of type $2^t$. Start from this GDD. Apply Construction \ref{SCHGDD-from strictly CGDD} with a $3$-MGDD of type $k^3$, $k\in\{3,4,5\}$, which exists by Theorem \ref{3-HGDD}, to get a $3$-SCHGDD of type $(3,2^t)$. It is equivalent to a $(3,2t;2)$-CHDM.  \qed

\begin{Lemma}\label{3-chdm-4} There exists a $(3,4t;4)$-CHDM for any integer $t\geq 4$.
\end{Lemma}

\proof When $t\in\{5,8,11\}$, we list a $(3,4t;4)$-CHDM as follows:

\begin{center}
{\footnotesize\tabcolsep 0.043in
\noindent {\rm $\left(
         \begin{tabular}{cccccccccccccccccccccccccc}
0&0&0&0&0&0&0&0&0&0&0&0&0&0&0&0\\
1&2&3&4&6&7&8&9&11&12&13&14&16&17&18&19\\
2&4&6&8&12&14&17&1&19&3&7&11&9&13&16&18\\
         \end{tabular}
       \right);
$}}
\end{center}

\begin{center}
{\footnotesize\tabcolsep 0.043in
\noindent {\rm $\left(
         \begin{tabular}{ccccccccccccccccccccccccccccccccccccccccccccc}
0&0&0&0&0&0&0&0&0&0&0&0&0&0&0&0&0&0&0&0&0&0&0&0&0&0&0&0\\
1&2&3&4&5&6&7&9&10&11&12&13&14&15&17&18&19&20&21&22&23&25&26&27&28&29&30&31\\
2&4&6&9&11&10&14&18&20&22&25&27&29&1&5&3&31&7&12&15&13&19&21&23&17&26&28&30\\
         \end{tabular}
       \right);
$}}
\end{center}

\begin{center}
{\footnotesize\tabcolsep 0.043in
\noindent {\rm $\left(
         \begin{tabular}{ccccccccccccccccccccccccccccccccccccccccccccccccccccccccccccccccc}
0&0&0&0&0&0&0&0&0&0&0&0&0&0&0&0&0&0&0&0\\
1&2&3&4&5&6&7&8&9&10&12&13&14&15&16&17&18&19&20&21\\
2&4&6&8&10&12&14&16&18&24&25&28&26&31&41&43&35&38&40&39\\
\hline
0&0&0&0&0&0&0&0&0&0&0&0&0&0&0&0&0&0&0&0\\
23&24&25&26&27&28&29&30&31&32&34&35&36&37&38&39&40&41&42&43\\
7&17&5&3&15&13&23&9&21&27&20&1&19&34&30&37&36&32&29&42\\

         \end{tabular}
       \right).
$}}
\end{center}

\noindent When $t\equiv 0\ ({\rm mod}\ 2)$, $t\geq 4$ and $t\neq8$, by Corollary \ref{4t,4,PSCHGDD}, there exists a perfect $3$-SCHGDD of type $(3,4^t)$, which is also a $(3,4t;4)$-CHDM by Lemma \ref{relation SCH and CHDM}. When $t\equiv 1\ ({\rm mod}\ 2)$, $t\geq 7$ and $t\neq11$, by Lemma \ref{4t,2p}, there exists a $(4t,4,\{3,5\},1)$-CDF. Then apply Construction \ref{SCHGDD-from strictly CGDD} with a $3$-MGDD of type $k^3$, $k\in\{3,5\}$, which exists by Theorem \ref{3-HGDD}, to obtain a $(3,4t;4)$-CHDM. \qed

\begin{Lemma}\label{3-chdm-2^x} There exists a $(3,2^xt;2^x)$-CHDM for any integers $x\geq 3$ and $t\geq 5$.
\end{Lemma}

\proof For any odd integer $t\geq 5$ and $x\geq 3$, start from a $3$-SCHGDD of type $(3,1^t)$, which is equivalent to a $(3,t;1)$-CHDM and exists by Example \ref{3-t-1-CHDM}. Then apply Construction \ref{SCHGDD from CHDM} with $w=1$, $h=2$ and $g=2^{x-1}$ to obtain a $(3,2^xt;2^x)$-CHDM, where the needed $3$-SCHGDD of type $(3,2^t)$ (i.e., a $(3,2t;2)$-CHDM) and the needed  $(3,2^x;2)$-CHDM are both from Lemma \ref{3-chdm-2}.

For $(x,t)=(3,6)$, we list a $(3,8\times6;8)$-CHDM as follows:
\begin{center}
{\footnotesize\tabcolsep 0.043in
\noindent {\rm $\left(
         \begin{tabular}{ccccccccccccccccccccccccccccccccccccccccccccccccccccccccccccccccc}
0&0&0&0&0&0&0&0&0&0&0&0&0&0&0&0&0&0&0&0\\
1&2&3&4&5&7&8&9&10&11&13&14&15&16&17&19&20&21&22&23\\
2&4&7&9&8&14&16&19&21&20&26&28&31&33&32&38&40&43&47&46\\

\hline
0&0&0&0&0&0&0&0&0&0&0&0&0&0&0&0&0&0&0&0\\
25&26&27&28&29&31&32&33&34&35&37&38&39&40&41&43&44&45&46&47\\
3&5&10&1&13&11&17&22&15&25&23&29&34&27&37&35&41&44&39&45\\

         \end{tabular}
       \right).
$}}
\end{center}

\noindent For any even integer $t\geq 6$, $x\in\{3,4\}$ and $(x,t)\neq(3,6)$, by Lemmas \ref{8t5} and \ref{16t,CDF}, there is a $(2^xt,2^x,\{3,5\},1)$-CDF. Then apply Construction \ref{SCHGDD-from strictly CGDD} with a $3$-MGDD of type $k^3$, $k\in\{3,5\}$, which exists by Theorem \ref{3-HGDD}, to obtain a $3$-SCHGDD of type $(3,(2^x)^t)$. It is also a $(3,2^xt;2^x)$-CHDM.

For $t=6$ and $x=5$, by Example \ref{example-h-PDF} there is a $(32\times6,32,\{3,5\},1)$-$2$-PDF. Then applying Construction \ref{perfect SCHGDD-from perfect DF} with a $3$-MGDD of type $k^3$, $k\in\{3,5\}$, we have a $2$-perfect $3$-SCHGDD of type $(3,32^6)$. It is also a $(3,32\times6;32)$-CHDM.

For $t=6$ and $x=6$, by Example \ref{example-h-PDF-h=4}, there is a $(16\times6,16,\{3,5\},1)$-$4$-PDF. Then applying Construction \ref{perfect SCHGDD-from perfect DF} with a $3$-MGDD of type $k^3$, $k\in\{3,5\}$, we have a $4$-perfect $3$-SCHGDD of type $(3,16^6)$. Start from a perfect $3$-SCHGDD of type $(3,4^6)$, which exists by Corollary \ref{4t,4,PSCHGDD}. Apply Construction \ref{SCHGDD from perfect} with $w=h=g=4$, we obtain a $(3,64\times6;64)$-CHDM, where the needed $(3,16;4)$-CHDM is from Lemma \ref{3-chdm-4}.

For $t=6$ and $x\geq 7$, by Example \ref{example-PRDF}, there is a $(16\times6,16,\{3,5\},1)$-PDF. Then applying Construction \ref{perfect SCHGDD-from perfect DF} with a $3$-MGDD of type $k^3$, $k\in\{3,5\}$, we have a perfect $3$-SCHGDD of type $(3,16^6)$. Now using Construction \ref{SCHGDD from perfect} with $w=16$, $h=2$ and $g=2^{x-5}$, we obtain a $(3,2^{x}\times6;2^x)$-CHDM, where the needed $(3,2^{x-4};2)$-CHDM is from Lemma \ref{3-chdm-2}, and the needed $2$-perfect $3$-SCHGDD of type $(3,32^6)$ exists by the first paragraph of this proof.

For $t=8$ and $x=5$, by Example \ref{example-CDF}, there exists a $(32\times8,32,\{3,5\},1)$-CDF, which implies a strictly cyclic $\{3,5\}$-GDD of type $32^8$. Start from this GDD. Apply Construction \ref{SCHGDD-from strictly CGDD} with a $3$-MGDD of type $k^3$, $k\in\{3,5\}$, to obtain a $(3,32\times8;32)$-CHDM.

For $t=8$ and $x\geq 6$, by Example \ref{example-PRDF}, there is a $(8\times8,8,\{3,5\},1)$-PDF. Then applying Construction \ref{perfect SCHGDD-from perfect DF} with a $3$-MGDD of type $k^3$, $k\in\{3,5\}$, we have a perfect $3$-SCHGDD of type $(3,8^8)$. Similarly, by Example \ref{example-h-PDF} there is a $(128,16,\{3,5\},1)$-$2$-PDF, which yields a $2$-perfect $3$-SCHGDD of type $(3,16^8)$. Now using Construction \ref{SCHGDD from perfect} with $w=8$, $h=2$ and $g=2^{x-4}$, we obtain a $(3,2^{x}\times8;2^x)$-CHDM, where the needed $(3,2^{x-3};2)$-CHDM is from Lemma \ref{3-chdm-2}.

For any even integer $t\geq 10$ and $x\geq 5$, start from a perfect $3$-SCHGDD of type $(3,4^t)$, which exists by Corollary \ref{4t,4,PSCHGDD}. Then apply Construction \ref{SCHGDD from perfect} with $w=4$, $h=2$ and $g=2^{x-3}$ to obtain a $(3,2^{x}t;2^x)$-CHDM, where the needed $2$-perfect $3$-SCHGDD of type $(3,8^t)$ is from Corollary \ref{8t,8,PSCHGDD} and the needed $(3,2^{x-2};2)$-CHDM is from Lemma \ref{3-chdm-2}. \qed

\begin{Lemma}\label{3-chdm-no-1}
There is no $(3,3m;m)$-CHDM for any positive integer $m\equiv 0\ ({\rm mod}\ 2)$.
\end{Lemma}

\proof Suppose that there were a $(3,3m;m)$-CHDM $D=(d_{ij})_{3\times 2m}$. For each $1\leq j\leq 2m$, let $d'_{1j}=0$, $d'_{2j}=d_{2j}-d_{1j}$ and $d'_{3j}=d_{3j}-d_{1j}$. Obviously $D'=(d'_{ij})$ is also a $(3,3m;m)$-CHDM. Thus without loss of generality, we can always assume that
\begin{center}
$D=\left(
         \begin{array}{ccccc;{2pt/2pt}ccccc}
           0 & 0 & 0 & \cdots & 0 & 0 & 0 & 0 & \cdots & 0 \\
           1 & 4 & 7 & \cdots & 3m-2 & 2 & 5 & 8 & \cdots & 3m-1\\
           d_{31} & d_{32} & d_{33} & \cdots & d_{3m} & d_{3,m+1} & d_{3,m+2} & d_{3,m+3} & \cdots & d_{3,2m} \\
         \end{array}
       \right),
$\end{center}
\noindent where $\{d_{31},d_{32},\ldots,d_{3m}\}=\{3l+2:0\leq l\leq m-1\}$ and $\{d_{3,m+1},d_{3,m+2},\ldots,d_{3,2m}\}=\{3l+1:0\leq l\leq m-1\}$. Hence,
\begin{align}
& \sum_{j=1}^m d_{2j}=\frac{m}{2}(3m-1)\equiv -\frac{m}{2}\ ({\rm mod}\ 3m),  \\
& \sum_{j=1}^m d_{3j}=\frac{m}{2}(3m+1)\equiv \frac{m}{2}\ ({\rm mod}\ 3m).
\end{align}
Counting $(2)-(1)$, we have $\sum_{j=1}^m (d_{3j}-d_{2j})\equiv m\ ({\rm mod}\ 3m).$
However, $\{d_{3j}-d_{2j}\ ({\rm mod}\ 3m):1\leq j\leq m\}=\{3l+1:0\leq l\leq m-1\}$. Thus
$$ \sum_{j=1}^m (d_{3j}-d_{2j})=\frac{m}{2}(3m-1)\equiv -\frac{m}{2}\ ({\rm mod}\ 3m).
$$
A contradiction occurs. \qed

\begin{Lemma}\label{3-chdm-no-2}
There is no $(3,mt;m)$-CHDM for any positive integer $m\equiv 1\ ({\rm mod}\ 2)$ and $t\equiv 0\ ({\rm mod}\ 2)$.
\end{Lemma}

\proof Suppose that there were a $(3,mt;m)$-CHDM for $m\equiv 1\ ({\rm mod}\ 2)$ and $t\equiv 0\ ({\rm mod}\ 2)$. Applying Construction \ref{SCHGDD from SCHGDD} and Lemma \ref{relation SCH and CHDM} with a $(3,m;1)$-CHDM, which exists by Example \ref{3-t-1-CHDM}, we have a $(3,mt;1)$-CHDM. However, Lemma \ref{SCGDD from SCHGDD} shows that a $(3,mt;1)$-CHDM implies a $3$-SCGDD of type $(mt)^3$, which does not exist by Lemma \ref{scgdd}. A contradiction occurs. \qed

\vspace{0.3cm}
\noindent \textbf{Proof of Theorem \ref{3-chdm}} When $m\equiv 1\ ({\rm mod}\ 2)$ and $t\equiv 0\ ({\rm mod}\ 2)$, or $m\equiv 0\ ({\rm mod}\ 2)$ and $t=3$, there is no $(3,mt;m)$-CHDM by Lemmas \ref{3-chdm-no-1} and \ref{3-chdm-no-2}. Then the necessity follows from Lemma \ref{nece}. It suffices to consider the sufficiency.

For $m\equiv 1\ ({\rm mod}\ 2)$ and $t\equiv 1\ ({\rm mod}\ 2)$, start from a $(3,t;1)$-CHDM, which exists by Example \ref{3-t-1-CHDM}. By Lemmas \ref{scgdd} and \ref{relation SC and CDM}, there exists a $(3,m)$-CDM. Then apply Construction \ref{SCHGDD from CDM} and Lemma \ref{relation SCH and CHDM} to obtain a $(3,mt;m)$-CHDM.

For $m\equiv 2\ ({\rm mod}\ 4)$ and $t\geq 4$, start from a $(3,2t;2)$-CHDM, which exists by Lemma \ref{3-chdm-2}. By Lemmas \ref{scgdd} and \ref{relation SC and CDM}, there exists a $(3,m/2)$-CDM. Then apply Construction \ref{SCHGDD from CDM} and Lemma \ref{relation SCH and CHDM} to obtain a $(3,mt;m)$-CHDM.

For $m\equiv 0\ ({\rm mod}\ 12)$ and $t\geq 4$, or $m\equiv 4,8\ ({\rm mod}\ 12)$, $t\equiv 1\ ({\rm mod}\ 3)$ and $t\geq 4$, start with a strictly cyclic $3$-GDD of type $m^t$ from Lemma \ref{3-sCGDD}. Apply Construction \ref{SCHGDD-from strictly CGDD} with a $3$-MGDD of type $3^3$, which exists by Theorem \ref{3-HGDD}. Then we have a $3$-SCHGDD of type $(3,m^t)$. It yields a $(3,mt;m)$-CHDM.

For $m\equiv 4,8\ ({\rm mod}\ 12)$, $t\equiv 0,2\ ({\rm mod}\ 3)$ and $t\geq 5$, write $m=2^x u$, where $x\geq 2$ and $u\equiv 1\ ({\rm mod}\ 2)$. Start from a $(3,2^xt;2^x)$-CHDM, which exists by Lemmas \ref{3-chdm-4} and \ref{3-chdm-2^x}. Then apply Construction \ref{SCHGDD from CDM} with a $(3,u)$-CDM  to obtain a $(3,mt;m)$-CHDM. \qed

\section{The cases of $n=4,5,6,8$}

In this section we shall establish the existence of $3$-SCHGDDs of type $(n,m^t)$ for $n=4,5,6,8$ and $t$ odd.

\begin{Lemma}\label{4-1^t}{\rm \cite{wsy}}
\item{$(1)$} There exists a $3$-SCHGDD of type $(4,1^t)$ for any odd integer $t\geq 3$.
\item{$(2)$} There exists a $3$-SCHGDD of type $(6,1^t)$ for any integer $t\equiv 1,5\ ({\rm mod}\ 6)$ and $t\geq 5$.
\item{$(3)$} There exists a $3$-SCHGDD of type $(6,1^9)$.
\end{Lemma}

\begin{Lemma}\label{4-SCHGDD}
There exists a $3$-SCHGDD of type $(4,m^t)$ for any positive integer $m$ and any odd integer $t\geq 3$.
\end{Lemma}

\proof For $m\equiv 0\ ({\rm mod}\ 2)$, start from a $3$-SCGDD of type $m^4$, which exists by Lemma \ref{scgdd}. By Example \ref{3-1^t}, there exists a $3$-SCHGDD of type $(3,1^t)$ for any odd integer $t\geq 3$. Then apply Construction \ref{SCHGDD-recur} to obtain a $3$-SCHGDD of type $(4,m^t)$.

For $m\equiv 1\ ({\rm mod}\ 2)$, start from a $3$-SCHGDD of type $(4,1^t)$, which exists for any odd integer $t\geq 3$ by Lemma \ref{4-1^t}. By Lemmas \ref{scgdd} and \ref{relation SC and CDM}, there exists a $(3,m)$-CDM. Then apply Construction \ref{SCHGDD from CDM} to obtain a $3$-SCHGDD of type $(4,m^t)$. \qed

\begin{Lemma}\label{5-SCHGDD}
Let $(t-1)m\equiv 0\ ({\rm mod}\ 6)$ and $t\geq 3$ be an odd integer. There exists a $3$-SCHGDD of type $(5,m^t)$.
\end{Lemma}

\proof For $m\equiv 0\ ({\rm mod}\ 3)$ and $t\equiv 1\ ({\rm mod}\ 2)$, start from a $3$-SCGDD of type $m^5$, which exists by Lemma \ref{scgdd}. By Example \ref{3-1^t}, there exists a $3$-SCHGDD of type $(3,1^t)$ for any odd integer $t\geq 3$. Then apply Construction \ref{SCHGDD-recur} to obtain a $3$-SCHGDD of type $(5,m^t)$.

For $m=2$ and $t\equiv 1\ ({\rm mod}\ 6)$, by Example \ref{example-PDF} and Lemma \ref{PDF-24-19}, there is a $(2t-1,\{3,4\},1)$-PDF, which yields a $(2t,2,\{3,4\},1)$-PDF from Lemma \ref{relation PDF and PDF}. Thus we have a strictly cyclic $\{3,4\}$-GDD of type $2^t$. Start from this GDD and apply Construction \ref{SCHGDD-from strictly CGDD} with a $3$-MGDD of type $k^5$, $k\in\{3,4\}$, which exists by Theorem \ref{3-HGDD}, to obtain a $3$-SCHGDD of type $(5,2^t)$. For $m\equiv 2,10\ ({\rm mod}\ 12)$ and $t\equiv 1\ ({\rm mod}\ 6)$, start from the resulting $3$-SCHGDD of type $(5,2^t)$, and apply Construction \ref{SCHGDD from CDM} with a $(3,m/2)$-CDM to obtain a $3$-SCHGDD of type $(5,m^t)$.

For $m\equiv 1,4,5,7,8,11\ ({\rm mod}\ 12)$ and $t\equiv 1\ ({\rm mod}\ 6)$, by Lemma \ref{3-sCGDD}, we have a strictly cyclic $3$-GDD of type $m^t$. Then apply Construction \ref{SCHGDD-from strictly CGDD} with a $3$-MGDD of type $3^5$ to obtain a $3$-SCHGDD of type $(5,m^t)$. \qed

\begin{Lemma}\label{6-SCHGDD-1^3}
There is no $3$-SCHGDD of type $(6,1^3)$.
\end{Lemma}

\proof Suppose that there exists a $3$-SCHGDD of type $(6,1^3)$ on $I_6\times Z_3$ with the group set $\{\{(i,0),(i,1),(i,2)\}:0\leq i\leq 5\}$ and the hole set $\{\{(i,j):0\leq i\leq 5\}:0\leq j\leq 2\}$. It has $10$ base blocks. Denote the set of its base blocks by ${\cal B}^*$. Let ${\cal B}^*=\{\{(a_1^{(l)},x_1^{(l)}),(a_2^{(l)},x_2^{(l)}),(a_3^{(l)},x_3^{(l)})\}:0\leq l\leq 9\}$. Given $b,c\in I_6$ and $b\neq c$, consider the number of base blocks containing the pairs of the form $\{(b,*),(c,*)\}$. The number is $2$. Now write ${\cal A}=\{\{a_1^{(l)},a_2^{(l)},a_3^{(l)}\}:0\leq l\leq 9\}$. We have that ${\cal A}$ forms a $(6,3,2)$-BIBD on $I_6$. Up to isomorphism, there is only one $(6,3,2)$-BIBD \cite{mr}. Hence, without loss of generality the $10$ base blocks can be assumed as follows.
\begin{center}\begin{tabular}{lll}
$\{(0,0),(1,x_2^{(0)}),(2,x_3^{(0)})\}$,&
$\{(1,0),(2,x_2^{(1)}),(5,x_3^{(1)})\}$,&
$\{(0,0),(2,x_2^{(2)}),(4,x_3^{(2)})\}$,\\
$\{(0,0),(4,x_2^{(3)}),(5,x_3^{(3)})\}$,&
$\{(1,0),(4,x_2^{(4)}),(5,x_3^{(4)})\}$,&
$\{(0,0),(1,x_2^{(5)}),(3,x_3^{(5)})\}$,\\
$\{(0,0),(3,x_2^{(6)}),(5,x_3^{(6)})\}$,&
$\{(1,0),(3,x_2^{(7)}),(4,x_3^{(7)})\}$,&
$\{(2,0),(3,x_2^{(8)}),(4,x_3^{(8)})\}$,\\
$\{(2,0),(3,x_2^{(9)}),(5,x_3^{(9)})\}$.
\end{tabular}
\end{center}
Note that $\{x_2^{(l)},x_3^{(l)}\}=\{1,2\}$ and $x_2^{(l)}\neq x_3^{(l)}$ for each $0\leq l\leq 9$. If $x_2^{(0)}=1$, then one can finish the first $5$ base blocks above as follows.
\begin{center}\begin{tabular}{lll}
$\{(0,0),(1,1),(2,2)\}$,&
$\{(1,0),(2,2),(5,1)\}$,&
$\{(0,0),(2,1),(4,2)\}$,\\
$\{(0,0),(4,1),(5,2)\}$,&
$\{(1,0),(4,1),(5,2)\}$.
\end{tabular}
\end{center}
So the pair $\{(4,1),(5,2)\}$ appears in two base blocks, a contradiction. Similar argument holds for $x_2^{(0)}=2$. \qed

\begin{Lemma}\label{6-SCHGDD-m^3}
There exists a $3$-SCHGDD of type $(6,m^3)$ for any integer $m\equiv 3\ ({\rm mod}\ 6)$.
\end{Lemma}

\proof When $m=3$, let $I=\{0,1,2,3,4,\infty\}$. We here construct a $3$-SCHGDD of type $(6,3^3)$ on $I\times Z_9$ with the group set $\{\{i\}\times Z_9:i\in I\}$ and the hole set $\{I\times \{j,3+j,6+j\}:0\leq j\leq 2\}$. Only initial base blocks are listed below, and all other base blocks are obtained by developing these base blocks by $(+1,-)$ modulo $(5,-)$, where $\infty+1=\infty$.

\begin{center}\begin{tabular}{lll}
$\{(0,0),(1,1),(3,2)\}$,&
$\{(0,0),(1,2),(3,4)\}$,&
$\{(0,0),(1,4),(\infty,8)\}$,\\
$\{(0,0),(1,5),(4,1)\}$,&
$\{(0,0),(1,7),(\infty,5)\}$,&
$\{(0,0),(2,8),(\infty,1)\}$.
\end{tabular}
\end{center}

\noindent When $m\geq 9$, start from the resulting $3$-SCHGDD of type $(6,3^3)$. Then apply Construction \ref{SCHGDD from CDM} with a $(3,m/3)$-CDM, which exists by Lemmas \ref{scgdd} and \ref{relation SC and CDM}, to obtain a $3$-SCHGDD of type $(6,m^3)$. \qed

\begin{Lemma}\label{6-SCHGDD-1^t}
There exists a $3$-SCHGDD of type $(6,1^t)$ for any odd integer $t\not\equiv 3,15\ ({\rm mod}\ 18)$.
\end{Lemma}

\proof When $t\equiv 1,5\ ({\rm mod}\ 6)$ and $t\geq 5$, the conclusion follows from Lemma \ref{4-1^t}. When $t\equiv 9\ ({\rm mod}\ 18)$, write $t=3^r u$, where $u\equiv 1,5\ ({\rm mod}\ 6)$ and $r\geq 2$.

If $u=1$, we use induction on $r$. When $r=2$, the conclusion follows from Lemma \ref{4-1^t}. When $r\geq 3$, assume that there exists a $3$-SCHGDD of type $(6,1^{3^{r-1}})$. By Lemma \ref{6-SCHGDD-m^3} we have a $3$-SCHGDD of type $(6,(3^{r-1})^3)$. Then apply Construction \ref{SCHGDD from SCHGDD} with the given $3$-SCHGDD of type $(6,1^{3^{r-1}})$, we have the required $3$-SCHGDD of type $(6,1^{3^r})$.

If $u\geq 2$, start from the resulting $3$-SCHGDD of type $(6,1^{3^r})$. Apply Construction \ref{SCHGDD from CDM} with a $(3,u)$-CDM, which exists by Lemmas \ref{scgdd} and \ref{relation SC and CDM}, to obtain a $3$-SCHGDD of type $(6,u^{3^r})$. Then making use of Construction \ref{SCHGDD from SCHGDD} with a $3$-SCHGDD of type $(6,1^u)$, which exists by Lemma \ref{4-1^t}, we have a $3$-SCHGDD of type $(6,1^{3^r u})$. \qed

Cyclotomic cosets play an important role in direct constructions for SCHGDDs. Let $p\equiv 1$ (mod $n$) be a prime and $\omega\in Z_p$
be a primitive element. Let $C_{0}^{n}$ denote the multiplicative
subgroup $\{{w}^{in}: 0\le i<(p-1)/n\}$ of the $n$-th powers in
$Z_p$ and $C_{j}^{n}$ denote the coset of $C_{0}^{n}$ in
$Z_{p}\setminus \{0\}$, i.e. $C_{j}^{n}={w}^{j}\cdot C_{0}^{n}$, $0\leq j\leq
n-1$. The following theorem is a variation of famous Weil's theorem on character sum, which can be also seen as a corollary of Theorem $2.2$ in \cite{mp}.

\begin{Theorem}\label{weil} {\rm \cite{cj}}
Let $p\equiv 1\ ({\rm mod}\ q)$ be a prime satisfying the inequality
$$
p-[\sum_{i=0}^{s-2} {{s}\choose{i}}(s-i-1)(q-1)^{s-i}]\sqrt{p}-sq^{s-1}>0.
$$
Then for any given $s$-tuple $(j_1,j_2,\ldots,j_s)\in
\{0,1,\ldots,q-1\}^s$ and any given $s$-tuple $(c_1,c_2$, $\ldots,c_s)$
of pairwise distinct elements of $Z_p$, there exists an element
$x\in Z_p$ such that $x+c_i\in C_{j_i}^q$ for each $i$.
\end{Theorem}

\begin{Lemma}\label{apply Weil Theorem} Let $p\geq 5$ be a prime. There exists an element $x\in Z_p$ such that $x\in C_1^2$, $x+1\in C_1^2$ and $x-1\in C_0^2$.
\end{Lemma}

\proof Take $q=2$ and $s=3$ in Theorem \ref{weil}. We have that if $p$ satisfies $p-5\sqrt{p}-12>0$, which yields $p>45$, then there exists an element $x\in Z_p$ such that $x\in C_1^2$, $x+1\in C_1^2$ and $x-1\in C_0^2$. When $5\leq p\leq 43$, it is readily checked that we may take $x$ as $(p,x)=(5,2),(7,5),(11,6),(13,5),(17,5),(19,2),(23,10),(29,2),(31,11),(37,5),(41,6),(43,2)$. \qed

\begin{Lemma}\label{4-SCHGDD-prime} There exists a $3$-SCHGDD of type $(6,2^p)$ for any prime $p\geq 3$.
\end{Lemma}

\proof Since $(p,2)=1$, $Z_{2p}$ is isomorphic to $Z_2\times Z_p$ under the mapping $\tau:\ x ({\rm mod}\ 2p)\rightarrow(x ({\rm mod}\ 2),x ({\rm mod}\ p))$. We construct the required $3$-SCHGDD of type $(6,2^p)$ on $I_6\times Z_2\times Z_p$. Let $S=\{(0,0),(1,0)\}$ be a subgroup of order $2$ in $Z_2\times Z_p$, and $S_l=S+(l,l)=\{(l,l),(l+1,l)\}$ (mod $(2,p)$) be a coset of $S$ in $Z_2\times Z_p$, $0\leq l\leq p-1$. Take the group set ${\cal G}=\{\{j\}\times Z_2\times Z_p:j\in I_6\}$, and the hole set ${\cal H}=\{I_6\times S_l:0\leq l\leq p-1\}$.

When $p=3$, all base blocks of a $3$-SCHGDD of type $(6,2^3)$ are listed as follows.

\begin{center}
\begin{tabular}{lll}
$\{(0,0,0),(1,1,1),(2,0,2)\}$,&
$\{(0,0,0),(2,1,1),(1,0,2)\}$,&
$\{(0,0,0),(3,1,1),(4,0,2)\}$,\\
$\{(0,0,0),(4,1,1),(3,0,2)\}$,&
$\{(0,0,0),(5,1,1),(1,1,2)\}$,&
$\{(0,0,0),(5,0,2),(1,0,1)\}$,\\
$\{(0,0,0),(2,0,1),(3,1,2)\}$,&
$\{(0,0,0),(3,0,1),(2,1,2)\}$,&
$\{(0,0,0),(4,0,1),(5,1,2)\}$,\\
$\{(0,0,0),(5,0,1),(4,1,2)\}$,&
$\{(1,0,0),(3,1,1),(4,1,2)\}$,&
$\{(1,0,0),(4,1,1),(2,0,2)\}$,\\
$\{(1,0,0),(5,1,1),(3,1,2)\}$,&
$\{(1,0,0),(3,0,2),(4,0,1)\}$,&
$\{(1,0,0),(4,0,2),(2,0,1)\}$,\\
$\{(1,0,0),(3,0,1),(5,1,2)\}$,&
$\{(2,0,0),(4,1,1),(5,1,2)\}$,&
$\{(2,0,0),(5,1,1),(3,0,2)\}$,\\
$\{(2,0,0),(4,0,2),(5,0,1)\}$,&
$\{(2,0,0),(5,0,2),(3,0,1)\}$.
\end{tabular}
\end{center}

When $p\geq 5$ is a prime, take $x\in Z_p$ such that $x\in C_1^2$, $x+1\in C_1^2$ and $x-1\in C_0^2$. This can be done by Lemma \ref{apply Weil Theorem}. Let $\omega$ be a primitive element of $Z_p$. All $10(p-1)$ base blocks of the required $3$-SCHGDD of type $(6,2^p)$ are generated by multiplying each of the following $20$ base blocks by $\omega^{2r}$, $0\leq r< (p-1)/2$.

$$\begin{array}{ll}
{\rm if}\ p\equiv 1\ ({\rm mod}\ 4):\\
\{(0,0,0),(1,0,1),(2,1,x)\},&\{(0,0,0),(1,1,1),(2,0,x+1)\},\\
\{(0,0,0),(1,0,x),(3,1,1)\},& \{(0,0,0),(1,1,x+1),(3,0,1)\},\\
\{(0,0,0),(2,0,1),(4,1,x)\},& \{(0,0,0),(2,1,1),(4,0,x+1)\},\\
\{(0,0,0),(3,0,x),(5,0,1)\},& \{(0,0,0),(3,1,x+1),(5,1,1)\},\\
\{(0,0,0),(4,0,1),(5,0,x)\},& \{(0,0,0),(4,1,1),(5,1,x+1)\},\\
\{(1,0,0),(2,0,1),(5,0,x+1)\},& \{(1,0,0),(2,0,x+1),(5,1,1)\},\\
\{(1,0,0),(3,0,1),(4,0,x)\},& \{(1,0,0),(3,0,x),(4,1,1)\},\\
\{(1,0,0),(4,0,1),(5,1,x+1)\},& \{(1,1,0),(4,0,x),(5,1,1)\},\\
\{(2,0,1),(3,0,0),(4,0,x)\},& \{(2,1,1),(3,0,0),(4,1,x+1)\},\\
\{(2,0,0),(3,0,x),(5,1,1)\},& \{(2,1,0),(3,0,x+1),(5,1,1)\}.
\end{array}$$

$$\begin{array}{ll}
{\rm if}\ p\equiv 3\ ({\rm mod}\ 4):\\
\{(0,0,0),(1,0,1),(2,1,x)\},& \{(0,0,0),(1,1,1),(2,0,x+1)\},\\
\{(0,0,0),(1,0,x),(3,1,1)\},& \{(0,0,0),(1,1,x),(3,0,-1)\},\\
\{(0,0,0),(2,0,1),(4,1,x)\},& \{(0,0,0),(2,1,1),(4,0,x+1)\},\\
\{(0,0,0),(3,0,1),(5,0,x)\},& \{(0,0,0),(3,1,x),(5,1,1)\},\\
\{(0,0,0),(4,0,-x),(5,0,1)\},& \{(0,0,0),(4,1,1),(5,1,x)\},\\
\{(1,0,0),(2,0,1),(5,0,x)\},& \{(1,0,0),(2,0,x),(5,1,1)\},\\
\{(1,0,0),(3,0,1),(4,0,x)\},& \{(1,0,0),(3,0,x),(4,1,1)\},\\
\{(1,0,0),(4,0,1),(5,1,x)\},& \{(1,1,0),(4,0,x),(5,1,1)\},\\
\{(2,0,1),(3,0,0),(4,0,x)\},& \{(2,1,x),(3,0,0),(4,1,1)\},\\
\{(2,0,0),(3,0,-x),(5,1,1)\},& \{(2,1,0),(3,0,x),(5,1,x+1)\}.
\end{array}$$

\noindent Let ${\cal B}^*$ be the set of all base block of the resulting $3$-SCHGDD of type $(6,2^p)$. For $(i,j)\in I_6\times I_6$, write $\Delta_{ij}({\cal B}^*)=\bigcup_{B\in{\cal B}^*}\Delta_{ij}(B)=\bigcup_{B\in{\cal B}^*}\{(y_1-y_2,z_1-z_2)\ ({\rm mod}\ (2,p)):(i,y_1,z_1),(j,y_2,z_2)\in {\cal B}^*\}$. It is readily checked that
$$\Delta_{ij}({\cal B}^*)=\left\{
\begin{array}{lll}
(Z_2\times Z_p)\setminus S, & {\rm if\ } (i,j)\in I_6\times I_6 {\rm \ and\ } i\neq j,\\
\emptyset, & {\rm otherwise}.\\
\end{array}
\right.\eqno{\Box}
$$

\begin{Lemma}\label{6-SCHGDD-2^t}
There exists a $3$-SCHGDD of type $(6,2^t)$ for any odd integer $t\geq 3$.
\end{Lemma}

\proof Let $t=p_1^{a_1}p_2^{a_2}\cdots p_t^{a_t}$ be its prime factorization, where $p_i\geq 3$ is a prime, $1\leq i\leq t$. Start from a $3$-SCHGDD of type $(6,2^{p_1})$, which exists by Lemma \ref{4-SCHGDD-prime}. Apply Construction \ref{SCHGDD from CDM} with a $(3,q)$-CDM, $q\in\{p_1,p_2,\ldots,p_t\}$, which exists by Lemmas \ref{scgdd} and \ref{relation SC and CDM}, to obtain a $3$-SCHGDD of type $(6,(2q)^{p_1})$. Then apply Construction \ref{SCHGDD from SCHGDD} with a $3$-SCHGDD of type $(6,2^q)$, which exists by Lemma \ref{4-SCHGDD-prime}, to obtain a $3$-SCHGDD of type $(6,2^{p_1 q})$. Repeating the above process will produce the required $3$-SCHGDD of type $(6,2^t)$ for any odd integer $t\geq 3$. \qed

\begin{Lemma}\label{6-SCHGDD}
There exists a $3$-SCHGDD of type $(6,m^t)$  for any positive integer $m$ and any odd integer $t\geq 3$ except possibly when $m\equiv 1,5\ ({\rm mod}\ 6)$ and $t\equiv 3,15\ ({\rm mod}\ 18)$.
\end{Lemma}

\proof When $m\equiv 2\ ({\rm mod}\ 4)$, start from a $3$-SCHGDD of type $(6,2^t)$, which exists for any odd integer $t\geq 3$ by Lemma \ref{6-SCHGDD-2^t}. By Lemmas \ref{scgdd} and \ref{relation SC and CDM}, there exists a $(3,m/2)$-CDM. Then apply Construction \ref{SCHGDD from CDM} to obtain a $3$-SCHGDD of type $(6,m^t)$.

When $m\equiv 0\ ({\rm mod}\ 4)$, start from a $3$-SCGDD of type $m^6$, which exists by Lemma \ref{scgdd}. By Example \ref{3-1^t}, there exists a $3$-SCHGDD of type $(3,1^t)$ for any odd integer $t\geq 3$. Then apply Construction \ref{SCHGDD-recur} to obtain a $3$-SCHGDD of type $(6,m^t)$.

When $m\equiv 1\ ({\rm mod}\ 2)$, $t\equiv 1\ ({\rm mod}\ 2)$ and $t\not\equiv 3,15\ ({\rm mod}\ 18)$, start from a $3$-SCHGDD of type $(6,1^t)$, which exists by Lemma \ref{6-SCHGDD-1^t}. Then apply Construction \ref{SCHGDD from CDM} with a $(3,m)$-CDM to obtain a $3$-SCHGDD of type $(6,m^t)$.

When $m\equiv 3\ ({\rm mod}\ 6)$ and $t\equiv 3,15\ ({\rm mod}\ 18)$, start from a $3$-SCHGDD of type $(6,t^3)$, which exists by Lemma \ref{6-SCHGDD-m^3}. Then applying Construction \ref{SCHGDD from SCHGDD} with a $3$-SCHGDD of type $(6,3^{t/3})$, which exists by the previous paragraph, we have a $3$-SCHGDD of type $(6,3^t)$. Now start from the resulting $3$-SCHGDD of type $(6,3^t)$, and apply Construction \ref{SCHGDD from CDM} with a $(3,m/3)$-CDM to obtain a $3$-SCHGDD of type $(6,m^t)$. \qed

\begin{Lemma}\label{8-SCHGDD}
Let $(t-1)m\equiv 0\ ({\rm mod}\ 6)$ and $t\geq 3$ be an odd integer. There exists a $3$-SCHGDD of type $(8,m^t)$ except possibly when $m\equiv 2,10\ ({\rm mod}\ 12)$ and $t\equiv 7\ ({\rm mod}\ 12)$.
\end{Lemma}

\proof When $m\equiv 0\ ({\rm mod}\ 6)$ and $t\equiv 1\ ({\rm mod}\ 2)$, start from a $3$-SCGDD of type $m^8$, which exists by Lemma \ref{scgdd}. By Example \ref{3-1^t}, there exists a $3$-SCHGDD of type $(3,1^t)$ for any odd integer $t\geq 3$. Then apply Construction \ref{SCHGDD-recur} to obtain a $3$-SCHGDD of type $(8,m^t)$.

When $m\equiv 3\ ({\rm mod}\ 6)$ and $t\equiv 1\ ({\rm mod}\ 2)$, start from a $4$-SCGDD of type $3^8$, which exists by Lemma \ref{4-scgdd}. By Lemma \ref{4-SCHGDD}, there exists a $3$-SCHGDD of type $(4,(m/3)^t)$. Then apply Construction \ref{SCHGDD-recur} to obtain a $3$-SCHGDD of type $(8,m^t)$.

When $m\equiv 1,4,5,7,8,11\ ({\rm mod}\ 12)$ and $t\equiv 1\ ({\rm mod}\ 6)$, or $m\equiv 2,10\ ({\rm mod}\ 12)$ and $t\equiv 1\ ({\rm mod}\ 12)$, take a strictly cyclic $3$-GDD of type $m^t$ from Lemma \ref{3-sCGDD}. Then apply Construction \ref{SCHGDD-from strictly CGDD} with a $3$-MGDD of type $3^8$, which exists by Theorem \ref{3-HGDD}, to obtain a $3$-SCHGDD of type $(8,m^t)$. \qed

\section{Proof of Theorem \ref{main theorem}}

\begin{Lemma}
\label{10,2^3} There exists a $3$-SCHGDD of type $(10,2^3)$.
\end{Lemma}

\proof We here give a construction of a $3$-HGDD of type $(10,2^3)$ on
$I_{60}$ with the group set $\{\{10i+j: 0\leq i\leq 5\} :$ $0\leq
j\leq 9\}$ and the hole set $\{\{i+10j,30+i+10j: 0\leq i\leq 9\} :0\leq
j\leq 2\}$. Let $\alpha=(0\ 10\ 20\ 30\ 40\ 50)(1\ 11\ 21\ 31\ 41\ 51)\cdots(9\ 19\ 29\ 39 \ 49\ 59)$ and $\beta=(0\ 2\ 4\ 6\ 8)(1\ 3\ 5\ 7\ 9)(10\ 12\ 14\ 16\ 18)(11\ 13\ 15\ 17\ 19) \cdots (50\ 52\ 54\ 56\ 58)(51\ 53 \ 55\ 57\ 59)$ be two permutations on $I_{60}$ and $G$ be the group generated by $\alpha$ and $\beta$. Only base blocks are listed below. All
other blocks are obtained by developing these base blocks under the
action of $G$. Obviously this design is isomorphic to a
$3$-SCHGDD of type $(10,2^3)$.

\begin{center}
{ \tabcolsep 0.05in
\begin{tabular}{lllllll}
$\{0,11,22\}$,&
$\{0,12,21\}$,&
$\{0,13,25\}$,&
$\{0,14,51\}$,&
$\{0,15,28\}$,&
$\{0,16,24\}$,&
$\{0,17,53\}$,\\
$\{0,23,41\}$,&
$\{0,26,43\}$,&
$\{0,29,45\}$,&
$\{0,49,55\}$,&
$\{1,15,23\}$.&
\end{tabular}}
\end{center}

\begin{Lemma}
\label{15,2^3} There exists a $3$-SCHGDD of type $(15,2^3)$.
\end{Lemma}

\proof We here give a construction of a $3$-HGDD of type $(15,2^3)$ on
$I_{90}$ with the group set $\{\{15i+j: 0\leq i\leq 5\} :$ $0\leq
j\leq 14\}$ and the hole set $\{\{i+15j,45+i+15j: 0\leq i\leq 14\} :0\leq
j\leq 2\}$. Let $\alpha=(0\ 15\ 30\ 45\ 60\ 75)(1\ 16\ 31\ 46\ 61\ 76)\cdots(14\ 31\ 44\ 59 \ 74\ 89)$ and $\beta=(0\ 3\ 6\ 9\ 12)(1\ 4\ 7\ 10\ 13)(2\ 5\ 8\ 11\ 14)(15\ 18\ 21\ 24\ 27)(16\ 19\ 22\ 25\ 28)(17\ 20\ 23\ 26\ 29) \cdots (75\ 78\ \\81\ 84\ 87)(76\ 79\ 82\ 85\ 88)(77\ 80\ 83\ 86\ 89)$ be two permutations on $I_{90}$ and $G$ be the group generated by $\alpha$ and $\beta$. Only base blocks are listed below. All
other blocks are obtained by developing these base blocks under the
action of $G$. Obviously this design is isomorphic to a
$3$-SCHGDD of type $(15,2^3)$.

\begin{center}
{ \tabcolsep 0.05in
\begin{tabular}{lllllll}
$\{0,16,32\}$,&
$\{0,17,31\}$,&
$\{0,18,37\}$,&
$\{0,20,33\}$,&
$\{0,21,34\}$,&
$\{0,22,35\}$,&
$\{0,23,39\}$,\\
$\{0,24,36\}$,&
$\{0,25,42\}$,&
$\{0,26,38\}$,&
$\{0,29,40\}$,&
$\{0,41,61\}$,&
$\{0,43,62\}$,&
$\{0,44,65\}$,\\
$\{0,64,76\}$,&
$\{0,67,86\}$,&
$\{0,68,85\}$,&
$\{0,70,79\}$,&
$\{0,71,80\}$,&
$\{0,73,83\}$,&
$\{0,74,82\}$,\\
$\{1,19,40\}$,&
$\{1,23,43\}$,&
$\{1,32,65\}$,&
$\{1,34,62\}$,&
$\{1,37,74\}$,&
$\{1,41,68\}$,&
$\{2,20,41\}$.
\end{tabular}}
\end{center}

\begin{Lemma}
\label{18,2^3} There exists a $3$-SCHGDD of type $(18,2^3)$.
\end{Lemma}

\proof We here give a construction of a $3$-HGDD of type $(18,2^3)$ on
$I_{108}$ with the group set $\{\{18i+j: 0\leq i\leq 5\} :$ $0\leq
j\leq 17\}$ and the hole set $\{\{i+18j,54+i+18j: 0\leq i\leq 17\} :0\leq
j\leq 2\}$. Let $\alpha=(0\ 18\ 36\ 54\ 72\ 90)(1\ 19\ 37\ 55\ 73\ 91)\cdots(17\ 35\ 53\ 71\ 89\ 107)$ and $\beta=(0\ 6\ 12)(1\ 7\ 13)(2\ 8\ 14)(3\ 9\ 15)(4\ 10\ 16)(5\ 11\ 17)(18\ 24\ 32)(19\ 25\ 33)(20\ 26\ 34)(21\ 27\ 35)(22\ \\28\ 36)(23\ 29\ 37) \cdots (90\ 96\ 102)(91\ 97\ 103)(92\ 98\ 104)(93\ 99\ 105)(94\ 100\ 106)(95\ 101\ 107)$ be two permutations on $I_{108}$ and $G$ be the group generated by $\alpha$ and $\beta$. Only base blocks are listed below. All
other blocks are obtained by developing these base blocks under the
action of $G$. Obviously this design is isomorphic to a
$3$-SCHGDD of type $(18,2^3)$.

\begin{center}
{ \tabcolsep 0.05in
\begin{tabular}{lllllll}
$\{0,32,52\}$,&$\{0,49,73\}$,&$\{0,50,74\}$,&$\{0,51,86\}$,&
$\{0,76,97\}$,&$\{0,77,91\}$,&$\{0,79,101\}$,\\
$\{0,80,103\}$,&$\{0,82,94\}$,&$\{0,83,104\}$,&$\{0,85,98\}$,&
$\{0,88,37\}$,&$\{0,19,38\}$,&$\{0,20,39\}$,\\
$\{0,21,40\}$,&$\{0,22,41\}$,&$\{0,23,42\}$,&$\{0,24,53\}$,&
$\{0,25,45\}$,&$\{0,26,43\}$,&$\{0,27,44\}$,\\
$\{0,28,48\}$,&$\{0,30,46\}$,&$\{0,31,47\}$,&$\{0,33,92\}$,&
$\{0,35,93\}$,&$\{0,75,95\}$,&$\{0,81,105\}$,\\
$\{0,87,99\}$,&$\{0,89,100\}$,&$\{1,22,38\}$,&$\{1,25,39\}$,&
$\{1,26,40\}$,&$\{1,27,41\}$,&$\{1,28,45\}$,\\
$\{1,29,44\}$,&$\{1,31,46\}$,&$\{1,47,74\}$,&$\{1,49,80\}$,&
$\{1,51,75\}$,&$\{1,53,76\}$,&$\{1,77,93\}$,\\
$\{1,81,94\}$,&$\{1,82,105\}$,&$\{1,83,100\}$,&$\{1,86,107\}$,&
$\{1,87,98\}$,&$\{1,88,99\}$,&$\{1,89,101\}$,\\
$\{2,26,39\}$,&$\{2,27,50\}$,&$\{2,28,41\}$,&$\{2,29,52\}$,&
$\{2,32,47\}$,&$\{2,40,82\}$,&$\{2,45,88\}$,\\
$\{2,46,89\}$,&$\{2,53,87\}$,&$\{2,76,100\}$,&$\{2,77,101\}$,&
$\{2,81,106\}$,&$\{3,29,51\}$,&$\{3,40,88\}$,\\
$\{3,41,83\}$,&$\{3,47,82\}$,&$\{3,52,89\}$,&$\{3,76,101\}$,&
$\{4,53,83\}$.
\end{tabular}}
\end{center}

\begin{Lemma}
\label{27,2^3} There exists a $3$-SCHGDD of type $(27,2^3)$.
\end{Lemma}

\proof We here give a construction of a $3$-HGDD of type $(27,2^3)$ on
$I_{162}$ with the group set $\{\{27i+j: 0\leq i\leq 5\} :$ $0\leq
j\leq 26\}$ and the hole set $\{\{i+27j,81+i+27j: 0\leq i\leq 26\} :0\leq
j\leq 2\}$. Let $\alpha=(0\ 27\ 54\ 81\ 108\ 135)(1\ 28\ 55\ 82\ 109\ 136)\cdots(26\ 53\ 80\ 107\ 134\ 161)$ and $\beta=(0\ 3\ 6\ 9\ 12\ 15\ 18\ 21\ 24)(1\ 4\ 7\ 10\ 13\ 16\ 19\ 22\ 25)(2\ 5\ 8\ 11\ 14\ 17\ 20\ 22\ 26)(27\ 30\ 33\ 36\ 39\ 42\ 45\ 48\ \\51)(28\ 31\ 34\ 37\ 40\ 43\ 46\ 49\ 52)(29\ 32\ 35\ 38\ 41\ 44\ 47\ 50\ 53) \cdots (135\ 138\ 141\ 144\ 147\ 150\ 153\ \\156\ 159)(136\ 139\ 142\ 145\ 148\ 151\ 154\ 157\ 160)(137\ 140\ 143\ 146\ 149\ 152\ 155\ 158\ 161)$ be two permutations on $I_{162}$ and $G$ be the group generated by $\alpha$ and $\beta$. Only base blocks are listed below. All
other blocks are obtained by developing these base blocks under the
action of $G$. Obviously this design is isomorphic to a
$3$-SCHGDD of type $(27,2^3)$.

\begin{center}
{ \tabcolsep 0.05in
\begin{tabular}{llllll}
$\{0,28,56\}$,&$\{0,29,55\}$,&$\{0,30,61\}$,&
$\{0,32,57\}$,&$\{0,33,58\}$,&$\{0,34,59\}$,\\
$\{0,35,63\}$,&$\{0,36,60\}$,&$\{0,37,66\}$,&
$\{0,38,62\}$,&$\{0,39,65\}$,&$\{0,40,64\}$,\\
$\{0,41,70\}$,&$\{0,42,136\}$,&$\{0,43,73\}$,&
$\{0,44,67\}$,&$\{0,45,68\}$,&$\{0,46,69\}$,\\
$\{0,47,77\}$,&$\{0,48,140\}$,&$\{0,49,71\}$,&
$\{0,72,124\}$,&$\{0,74,109\}$,&$\{0,75,128\}$,\\
$\{0,76,111\}$,&$\{0,110,142\}$,&$\{0,112,143\}$,&
$\{0,113,146\}$,&$\{0,115,148\}$,&$\{0,116,154\}$,\\
$\{0,118,157\}$,&$\{0,122,158\}$,&$\{0,125,145\}$,&
$\{0,130,151\}$,&$\{0,131,149\}$,&$\{0,133,152\}$,\\
$\{0,134,155\}$,&$\{1,35,70\}$,&$\{1,37,74\}$,&
$\{1,41,56\}$,&$\{1,43,59\}$,&$\{1,46,149\}$,\\
$\{1,58,110\}$,&$\{1,61,122\}$,&$\{1,64,112\}$,&
$\{1,65,125\}$,&$\{1,67,116\}$,&$\{1,68,118\}$,\\
$\{1,71,119\}$,&$\{1,134,146\}$,&$\{2,59,122\}$,&
$\{2,68,113\}$.
\end{tabular}}
\end{center}

\begin{Lemma}\label{SCHGDD-n-2^3}
Let $n\equiv 0,1\ ({\rm mod}\ 3)$ and $n\geq 4$. There exists a $3$-SCHGDD of type $(n,2^3)$.
\end{Lemma}

\proof When $n\in\{4,6\}$, the conclusion follows from Lemmas \ref{4-SCHGDD} and \ref{6-SCHGDD}. When $n\in\{10,15,18,27\}$, the conclusion follows from Lemmas \ref{10,2^3}-\ref{27,2^3}. When $n\in\{9,12,24\}$, start from a $3$-SCGDD of type $2^n$, which exists by Lemma \ref{scgdd}. By Example \ref{3-1^t}, there exists a $3$-SCHGDD of type $(3,1^3)$. Then apply Construction \ref{SCHGDD-recur} to obtain a $3$-SCHGDD of type $(n,2^3)$. When $n\in\{7,19\}$, start from a $4$-SCGDD of type $2^n$, which exists by Lemma \ref{4-scgdd}. By Lemma \ref{4-SCHGDD}, there exists a $3$-SCHGDD of type $(4,1^3)$. Then apply Construction \ref{SCHGDD-recur} to obtain a $3$-SCHGDD of type $(n,2^3)$.

When $n\equiv 0,1\ ({\rm mod}\ 3)$, $n\geq 4$ and $n\not\in \{10,12,15,18,19,24,27\}$, start from a $\{4,6,7,9\}$-SCGDD of type $1^n$, which is also a $(n,\{4,6,7,9\},1)$-PBD and exists by Lemma \ref{346pbd}. Then apply Construction \ref{SCHGDD-recur} with a $3$-SCHGDD of type $(k,2^3)$, $k\in \{4,6,7,9\}$, to obtain a $3$-SCHGDD of type $(n,2^3)$.  \qed

\begin{Lemma}\label{SCHGDD-n-0-1-no 3}
Let $n\equiv 0,1\ ({\rm mod}\ 3)$ and $n\geq 7$. There exists a $3$-SCHGDD of type $(n,m^t)$ for any odd integer $t\geq 3$ and any positive integer $m$.
\end{Lemma}

\proof For any odd integer $t\geq 5$ and any positive integer $m$, or $t=3$ and $m\equiv 1\ ({\rm mod}\ 2)$, start from a $\{3,4\}$-SCGDD of type $1^n$, which is also a $(n,\{3,4\},1)$-PBD and exists by Lemma \ref{346pbd}. Take a $3$-SCHGDD of type $(k,m^t)$ for $k\in \{3,4\}$, which exists from Theorem \ref{3-chdm} and Lemma \ref{4-SCHGDD}, and then apply Construction \ref{SCHGDD-recur} to obtain a $3$-SCHGDD of type $(n,m^t)$.

For $t=3$ and $m\equiv 0\ ({\rm mod}\ 4)$, start from a $3$-SCGDD of type $m^n$, which exists by Lemma \ref{scgdd}. By Example \ref{3-1^t}, there exists a $3$-SCHGDD of type $(3,1^3)$. Then apply Construction \ref{SCHGDD-recur} to obtain a $3$-SCHGDD of type $(n,m^3)$.

For $t=3$ and $m\equiv 2\ ({\rm mod}\ 4)$, start from a $3$-SCHGDD of type $(n,2^3)$, which exists by Lemma \ref{SCHGDD-n-2^3}. By Lemmas \ref{scgdd} and \ref{relation SC and CDM}, there exists a $(3,m/2)$-CDM. Then apply Construction \ref{SCHGDD from CDM} to obtain a $3$-SCHGDD of type $(n,m^3)$. \qed

\begin{Lemma}\label{SCHGDD-n-2-(1)}
Let $n\equiv 2\ ({\rm mod}\ 3)$ and $n\geq 11$. Let $(t-1)m\equiv 0\ ({\rm mod}\ 6)$ and $t\geq 3$ be an odd integer. There exists a $3$-SCHGDD of type $(n,m^t)$.
\end{Lemma}

\proof When $(t-1)m\equiv 0\ ({\rm mod}\ 6)$, $t\geq 4$ is an odd integer and $m$ is a positive integer, or $t=3$ and $m\equiv 3\ ({\rm mod}\ 6)$, start from a $\{3,4,5\}$-SCGDD of type $1^n$, which is also a $(n,\{3,4,5\},1)$-PBD and exists by Lemma \ref{346pbd}. Take a $3$-SCHGDD of type $(k,m^t)$ for $k\in \{3,4,5\}$ from Theorem \ref{3-chdm}, Lemmas \ref{4-SCHGDD} and \ref{5-SCHGDD}. Then apply Construction \ref{SCHGDD-recur} to obtain a $3$-SCHGDD of type $(n,m^t)$.

When $t=3$ and $m\equiv 0\ ({\rm mod}\ 12)$, start from a $3$-SCGDD of type $m^n$, which exists by Lemma \ref{scgdd}. By Example \ref{3-1^t}, there exists a $3$-SCHGDD of type $(3,1^3)$. Then apply Construction \ref{SCHGDD-recur} to obtain a $3$-SCHGDD of type $(n,m^3)$.

When $t=3$ and $m=6$, start from a $4$-SCGDD of type $6^n$, which exists by Lemma \ref{4-scgdd}. By Lemma \ref{4-1^t}, there exists a $3$-SCHGDD of type $(4,1^3)$. Then apply Construction \ref{SCHGDD-recur} to obtain a $3$-SCHGDD of type $(n,6^3)$.
When $t=3$ and $m\equiv 6\ ({\rm mod}\ 12)$, start from the resulting $3$-SCHGDD of type $(n,6^3)$. By Lemmas \ref{scgdd} and \ref{relation SC and CDM}, there exists a $(3,m/6)$-CDM. Then apply Construction \ref{SCHGDD from CDM} to obtain a $3$-SCHGDD of type $(n,m^3)$.   \qed

\vspace{0.3cm}
\noindent \textbf{Proof of Theorem \ref{main theorem}} Combining the results of Lemmas \ref{nece}, \ref{4-SCHGDD}, \ref{5-SCHGDD}, \ref{6-SCHGDD-1^3}, \ref{6-SCHGDD},  \ref{8-SCHGDD}, \ref{SCHGDD-n-0-1-no 3} and \ref{SCHGDD-n-2-(1)}, we complete the proof. \qed

\newpage

\begin{center}
{\Large\bf Semi-cyclic holey group divisible designs with block size three: Appendix}

\vskip12pt

Tao Feng$^1$, Xiaomiao Wang$^2$ and Yanxun Chang$^1$ \\[2ex] {\footnotesize {\footnotesize $^1$Institute of Mathematics, Beijing Jiaotong
University, Beijing 100044, P. R. China}
\\$^2$Department of Mathematics, Ningbo University, Ningbo 315211, P. R. China}\\  {\footnotesize
tfeng@bjtu.edu.cn, wangxiaomiao@nbu.edu.cn, yxchang@bjtu.edu.cn}
\vskip12pt

\end{center}

This appendix is a supplement to Lemmas $4.3$, $4.9$, $4.15$ and $4.16$ of the paper {\it Semi-cyclic holey group divisible designs with block size three}, where Section A is for Lemma $4.3$, Sections B-E are for Lemma $4.9$, Section F is for Lemma $4.15$ and Section G is for Lemma $4.16$.

\appendix{}
\section{Small $(v,\{3,4,5\},1)$-PDFs for $v\equiv 5\ ({\rm mod}\ 6)$}

Here we list all the base blocks of a $(v,\{3,4,5\},1)$-PDF for $77<v<455$ and $v\not\in \{89,101$, $113,125,137,149,413,419,437,443\}$.

{\scriptsize
$$
\begin{array}{lllllll}

\bullet\ \ v=83: \\
  \{0,1,3,18,38\} & \{0,4,25,33\} & \{0,7,30,39\} & \{0,13,24,40\} & \{0,6,34\} &
  \{0,5,19,31,41\}\\

\bullet\ \ v=95: \\
  \{0,2,9,21,46\} & \{0,1,31,36\} & \{0,6,17,45\} & \{0,8,32,42\} & \{0,18,40\} &
  \{0,4,20,33,47\} & \{0,3,26,41\} \\

\bullet\ \ v=107: \\
  \{0,1,11,34,51\} & \{0,4,19,49\} & \{0,7,28,53\} & \{0,16,29,47\} & \{0,5,37\} &
  \{0,2,14,38,41\} & \{0,6,26,48\} \\ \{0,8,43,52\}  \\

\bullet\ \ v=119: \\
\end{array}$$

$$\begin{array}{lllllllll}
\bullet\ \ v=203: \\
  \{0,8,31,67,100\} & \{0,20,80,97\} & \{0,4,49,90\} & \{0,9,57,85\} & \{0,40,96\} &
  \{0,14,82,98\} & \{0,42,89\} \\ \{0,5,70,83\} & \{0,11,55,74\} & \{0,34,61,87\} &
  \{0,24,88,95\} & \{0,1,39,51\} & \{0,6,81,99\} & \{0,3,46\} \\ \{0,10,35,72,101\} &
  \{0,21,79,94\} & \{0,2,32,54\}  \\

\bullet\ \ v=209: \\
  \{0,8,37,81,103\} & \{0,3,49,101\} & \{0,9,68,99\} & \{0,25,78,96\} & \{0,4,58,91\} &
  \{0,21,76\} & \{0,11,26\} \\ \{0,13,70,102\} & \{0,28,84,100\} & \{0,12,39\} &
  \{0,1,35,83\} & \{0,5,47,97\} & \{0,14,65,88\} & \{0,6,67\} \\ \{0,10,40,85,104\} &
  \{0,2,38,79\} & \{0,7,69,93\} & \{0,20,63,80\} \\
\end{array}$$

$$\begin{array}{lllllll}

\bullet\ \ v=293: \\\{0,11,33,78,146\} & \{0,21,114,141\} & \{0,16,87,142\} & \{0,20,110,139\} & \{0,3,63,77\} &
  \{0,4,88,132\} &\{0,17,125\} \\ \{0,19,89,124\} & \{0,34,131,137\} & \{0,15,79,133\}  &
  \{0,49,101,144\} & \{0,5,96,134\} & \{0,32,115\} & \{0,7,73\} \\  \{0,23,51,92,145\} &
  \{0,57,107,143\} & \{0,13,61,117\} & \{0,18,100,130\} & \{0,1,47,59\} & \{0,8,80,106\}  &
  \{0,10,75\} \\ \{0,25,127,136\} & \{0,39,76,138\} & \{0,2,42,123\} & \{0,31,116,140\}  \\

\bullet\ \ v=299: \\
  \{0,10,41,81,148\} & \{0,25,124,146\} & \{0,11,79,143\} & \{0,33,109,136\} & \{0,6,55,98\} &
  \{0,3,53,91\} & \{0,51,135\} \\ \{0,18,122,130\} & \{0,29,106,129\} & \{0,9,70,142\} &
  \{0,16,117,141\} & \{0,1,46,94\} & \{0,2,58,115\} & \{0,39,126\} \\
\end{array}$$

$$\begin{array}{lllllll}

\bullet\ \    v=365:\\ \{0,29,93,151,182\} & \{0,37,134,167\} & \{0,43,94,174\} &  \{0,11,114,169\} & \{0,15,150,172\} & \{0,38,140\}
    &\{0,3,108\} \\ \{0,14,152,170\} & \{0,17,96,117\} & \{0,25,85,166\}  &
  \{0,49,111,178\} &  \{0,34,149,162\}&\{0,2,74,84\} &\{0,5,104\}\\  \{0,20,143,175\} &
  \{0,52,118,179\} & \{0,48,98,173\} & \{0,9,142,168\} & \{0,30,154,177\}&\{0,16,164\} &
  \{0,7,139\} \\\{0,36,90,137,181\} & \{0,78,119,161\} & \{0,71,77,163\} & \{0,56,113,176\}&
  \{0,68,121,180\} &\{0,1,46,70\} & \{0,4,110\}\\ \{0,76,88,116\} &\{0,35,144,171\}   &
   \{0,65,73,160\} & \{0,39,146,165\} \\

\bullet\ \    v=371: \\
  \{0,14,61,94,184\} & \{0,30,147,171\} &  \{0,8,67,132\}& \{0,35,88,173\} &\{0,12,139,146\} & \{0,1,55,96\}
   &  \{0,71,176\} \\
\end{array}$$}

\section{Small $(4t,4,\{3,5\},1)$-PDFs for $t\equiv 0\ ({\rm mod}\ 6)$}

Here we list all the base blocks of a $(4t,4,\{3,5\},1)$-PDF for $t\equiv 0\ ({\rm mod}\ 6)$ and $6\leq t\leq 54$.

{\scriptsize
$$\begin{array}{lllllll}

\bullet\ \ t=6: &
  \{0,1,4,9,11\}\\

\bullet\ \ t=12: &
 \{0,1,3,9,20\} & \{0,4,18\} & \{0,5,21\} & \{0,7,22\} & \{0,10,23\}\\

\bullet\ \  t=18: \\
  \{0,1,3,7,24\} & \{0,8,35\} & \{0,10,29\} & \{0,12,32\} & \{0,14,30\} &
  \{0,5,31\} &\{0,9,34\} \\ \{0,11,33\}&  \{0,13,28\} \\

\bullet\ \  t=24: \\
  \{0,1,7,20,41\} & \{0,4,30\} & \{0,9,47\}&\{0,12,44\}&\{0,16,39\} &
  \{0,2,17\} & \{0,5,36\} \\\{0,10,43\}& \{0,14,42\}& \{0,18,45\}&
  \{0,3,25\} &\{0,8,37\} &\{0,11,46\} \\

\bullet\ \  t=30: \\
  \{0,1,3,32,56\} & \{0,13,20,41,59\} &\{0,23,50\} & \{0,8,52\} &\{0,12,48\} &\{0,22,47\} &\{0,15,34\}\\\{0,4,9,42,58\}& \{0,6,17,43,57\}
   & \{0,10,45\} \\

\bullet\ \  t=36: \\
  \{0,1,3,38,68\} &  \{0,4,9,43,70\}& \{0,12,64\} &\{0,19,60\} & \{0,22,55\} &
 \{0,20,51\} & \{0,8,62\}\\ \{0,14,25,42,71\}& \{0,6,13,53,69\}& \{0,15,59\} &\{0,23,49\}
   & \{0,10,58\}& \{0,18,50\}&\{0,21,45\} \\

\bullet\ \  t=42: \\
  \{0,1,3,44,77\} &\{0,10,21,58,81\}&\{0,13,51,67,82\} &
  \{0,4,9,49,79\}&\{0,8,65\} &\{0,26,73\}  &\{0,14,66\}\\ \{0,12,34,62,80\} &\{0,20,27,59,83\} &
  \{0,6,25,61,78\}&\{0,29,64\}  \\

\bullet\ \  t=48: \\
  \{0,1,3,50,89\} & \{0,8,18,71,93\}& \{0,17,33,74,95\} &\{0,23,79\} &\{0,28,65\} &\{0,14,81\} &\{0,25,68\}
  \\ \{0,4,9,55,91\} & \{0,11,30,70,94\}& \{0,12,38,72,92\} &\{0,29,73\} &\{0,31,66\}
   & \{0,15,76\} &\{0,27,69\} \\\{0,6,13,58,90\} \\

\bullet\ \  t=54: \\
  \{0,1,3,56,99\} & \{0,26,67,89,101\}& \{0,24,72,92,103\} &\{0,18,64\} &\{0,32,81\} &\{0,10,94\} &\{0,21,80\}
  \\ \{0,4,9,66,104\} & \{0,16,58,87,102\}& \{0,23,50,83,97\} &\{0,36,76\} &
  \{0,39,91\} & \{0,13,106\} &\{0,28,105\} \\\{0,19,25,70,107\} &
  \{0,8,69\} & \{0,7,85\}& \{0,17,90\} &\{0,30,65\}   \\

\end{array}$$}

\section{Small $(4t,4,\{3,5\},1)$-PDFs for $t\equiv 2\ ({\rm mod}\ 6)$}

Here we list all the base blocks of a $(4t,4,\{3,5\},1)$-PDF for $t\equiv 2\ ({\rm mod}\ 6)$ and $14\leq t\leq 110$.

{\scriptsize
$$\begin{array}{lllllll}

\bullet\ \ t=14: &
 \{0,1,9,22,25\} & \{0,4,10,15,27\} & \{0,2,20\} & \{0,7,26\}\\

\bullet\ \  t=20: \\
  \{0,1,3,7,32\} &
  \{0,5,13,28,39\}& \{0,9,30\} & \{0,12,36\} & \{0,16,38\} & \{0,17,35\}  & \{0,10,37\} \\ \{0,14,33\} \\

\bullet\ \ t=26: \\
  \{0,1,4,25,38\} &
  \{0,2,7,29,43\}& \{0,8,47\} &\{0,12,44\} &\{0,17,50\}& \{0,19,42\}  & \{0,10,45\} \\\{0,16,46\} &\{0,18,49\}&\{0,20,48\} &
  \{0,6,15\} & \{0,11,51\} \\

\bullet\ \  t=32: \\
  \{0,1,3,34,59\} & \{0,9,20,46,63\}& \{0,15,42,50,55\} &\{0,14,21,44,62\} &\{0,24,53\} &\{0,16,52\} &\{0,22,60\}
  \\ \{0,4,10,49,61\} &  \{0,19,47\}\\

\bullet\ \  t=38: \\
  \{0,1,3,40,71\} & \{0,8,19,65,75\}& \{0,6,13,55,72\} &\{0,22,54\} &\{0,25,52\} &\{0,18,61\} &\{0,24,50\}
  \\ \{0,4,9,45,73\} & \{0,16,30,51,74\} &\{0,29,63\} &
  \{0,15,62\} & \{0,12,60\}& \{0,20,53\} \\

\bullet\ \  t=44: \\
\{0,1,3,46,83\} & \{0,20,32,60,86\}& \{0,16,78\} &\{0,22,51\} &\{0,25,63\} &\{0,19,77\}
 & \{0,7,56\}\\ \{0,5,9,74,84\} &
\{0,6,14,48,87\}&\{0,23,64\} &\{0,27,57\}  & \{0,11,61\}& \{0,21,76\} &\{0,24,71\} \\\{0,15,33,68,85\} &\{0,31,67\}
 & \{0,13,72\}  \\

\bullet\ \  t=50: \\
\{0,1,3,49,92\} & \{0,15,26,70,97\}& \{0,22,67,84,98\}& \{0,18,37,75,96\} &\{0,12,81\} &\{0,29,68\} & \{0,8,73\}
 \\\{0,4,9,60,94\}&
\{0,6,58,86,93\} & \{0,20,33,74,99\}& \{0,23,53,63,95\}&\{0,16,77\} &\{0,36,83\}  &\{0,24,88\} \\

\bullet\ \  t=56: \\
\{0,1,3,58,101\} & \{0,16,64,92,111\}& \{0,22,73,93,110\} &\{0,27,109\} &\{0,33,77\} &
\{0,21,106\}\\
\{0,4,9,63,103\}& \{0,6,13,66,102\} & \{0,25,75,90,104\}&\{0,11,83\} &\{0,31,80\} &\{0,34,108\} \\ \{0,8,18,70,105\} & \{0,26,38,68,107\}&\{0,23,84\} &\{0,32,78\} &\{0,41,86\}& \{0,24,91\}\\
\end{array}$$}

{\scriptsize
$$\begin{array}{llllll}
\bullet\ \  t=62: \\
  \{0,1,3,64,113\} & \{0,25,47,100,121\}& \{0,26,43,83,116\}& \{0,15,95\} &\{0,29,122\} &\{0,35,86\} \\\{0,4,9,69,115\}  &\{0,6,13,72,114\} & \{0,27,81,105,119\}&\{0,16,118\}
   &\{0,30,85\} &\{0,37,104\}
  \\ \{0,8,18,76,117\} &
  \{0,20,39,91,123\}&\{0,23,120\} &\{0,31,87\} &
  \{0,44,89\} & \{0,11,88\}\\ \{0,28,107\} &\{0,34,70\} &\{0,48,98\}  & \{0,12,94\}\\

\bullet\ \  t=68: \\
  \{0,1,3,70,125\} & \{0,29,56,91,133\}& \{0,23,38,82,130\} &\{0,11,100,116,128\} &\{0,40,113\} &\{0,13,112\} \\\{0,46,111\}
  & \{0,4,9,75,127\} &\{0,34,51,87,132\}& \{0,14,74,115,135\} &\{0,37,95\}
   & \{0,50,114\}\\
\end{array}$$}

\section{Partition the set $[1,2u+1]\setminus\{(u-14)/2\}$ into pairs}

Here for each $u\in[36,131]$ and $u\equiv 2\ ({\rm mod}\ 4)$, we give a partition of the set $[1,2u+1]\setminus\{(u-14)/2\}$ into pairs $\{x_i,y_i\}$, $1\leq i\leq u$, such that $\{y_i-x_i:1\leq i\leq u\}=[17,16+u]$.

{\scriptsize
$$\begin{array}{lllllllllll}
\bullet\ \ $u=38:$\\
$(13,41)$&
$(14,43)$&
$(15,45)$&
$(16,47)$&
$(17,58)$&
$(19,61)$&
$(21,60)$&
$(23,63)$&
$(25,62)$&
$(27,65)$&
$(29,64)$,\\
$(10,35)$&
$(11,59)$&
$(31,75)$&
$(32,77)$&
$(33,76)$&
$(36,68)$&
$(22,73)$&
$(24,44)$&
$(28,50)$&
$(34,53)$&
$(1,37)$,\\
$(2,49)$&
$(3,56)$&
$(4,38)$&
$(5,54)$&
$(6,39)$&
$(7,57)$&
$(8,26)$&
$(9,55)$&
$(18,70)$&
$(20,74)$&
$(30,51)$,\\
$(40,67)$&
$(42,66)$&
$(46,72)$&
$(48,71)$&
$(52,69)$\\
\end{array}$$}

\newpage
{\scriptsize
$$\begin{array}{lllllllllll}
\bullet\ \ $u=42:$\\
$(10,36)$&
$(21,55)$&
$(23,58)$&
$(25,61)$&
$(27,64)$&
$(29,74)$&
$(31,77)$&
$(33,80)$&
$(35,78)$&
$(37,79)$&
$(39,83)$,\\
$(40,81)$&
$(43,82)$&
$(45,85)$&
$(46,84)$&
$(11,59)$&
$(12,62)$&
$(13,65)$&
$(15,66)$&
$(16,69)$&
$(17,50)$&
$(19,68)$,\\
$(38,70)$&
$(41,71)$&
$(42,73)$&
$(1,56)$&
$(2,20)$&
$(3,57)$&
$(4,60)$&
$(5,22)$&
$(6,63)$&
$(7,26)$&
$(8,30)$,\\
$(9,67)$&
$(18,47)$&
$(24,49)$&
$(28,51)$&
$(32,53)$&
$(34,54)$&
$(44,72)$&
$(48,75)$&
$(52,76)$\\

\bullet\ \ $u=46:$\\
$(21,55)$&
$(23,58)$&
$(25,61)$&
$(27,64)$&
$(29,67)$&
$(31,70)$&
$(33,79)$&
$(35,84)$&
$(37,85)$&
$(39,86)$&
$(41,91)$\\
$(43,87)$&
$(45,88)$&
$(47,89)$&
$(48,93)$&
$(50,90)$&
$(51,92)$&
$(13,65)$&
$(14,68)$&
$(15,71)$&
$(17,72)$&
$(19,76)$\\
\end{array}$$

$$\begin{array}{llllllllll}

\bullet\ \ $u=82:$\\
$(85,141)$&
$(87,140)$&
$(33,106)$&
$(41,107)$&
$(69,109)$&
$(71,100)$&
$(88,147)$&
$(90,148)$&
$(122,145)$&
$(57,133)$\\
$(63,128)$&
$(70,132)$&
$(53,134)$&
$(55,137)$&
$(27,110)$&
$(44,78)$&
$(46,116)$&
$(48,115)$&
$(50,80)$&
$(23,95)$\\
$(21,49)$&
$(25,114)$&
$(54,139)$&
$(68,142)$&
$(79,157)$&
$(29,113)$&
$(37,127)$&
$(39,108)$&
$(11,102)$&
$(12,67)$\\
$(13,105)$&
$(83,124)$&
$(56,117)$&
$(59,130)$&
$(35,98)$&
$(62,99)$&
$(64,103)$&
$(65,101)$&
$(26,119)$&
$(97,165)$\\
$(104,151)$&
$(112,163)$&
$(38,125)$&
$(45,131)$&
$(47,135)$&
$(86,138)$&
$(89,153)$&
$(92,152)$&
$(96,144)$&
$(91,136)$\\
$(94,143)$&
$(118,164)$&
$(123,158)$&
$(1,18)$&
$(2,20)$&
$(3,22)$&
$(4,24)$&
$(5,30)$&
$(6,28)$&
$(7,31)$\\
$(8,40)$&
$(9,84)$&
$(10,60)$&
$(14,111)$&
$(15,36)$&
$(16,42)$&
$(17,74)$&
$(19,73)$&
$(32,75)$&
$(43,76)$\\
\end{array}$$

$$\begin{array}{llllllllll}

\bullet\ \ $u=106:$\\
$(85,141)$&
$(111,174)$&
$(112,172)$&
$(119,192)$&
$(105,197)$&
$(106,189)$&
$(110,186)$&
$(120,185)$&
$(122,176)$&
$(127,182)$\\
$(60,163)$&
$(65,146)$&
$(97,165)$&
$(41,107)$&
$(115,203)$&
$(133,184)$&
$(134,198)$&
$(135,205)$&
$(98,173)$&
$(121,152)$\\
$(36,93)$&
$(143,180)$&
$(73,171)$&
$(70,132)$&
$(45,131)$&
$(101,195)$&
$(103,183)$&
$(118,200)$&
$(125,204)$&
$(61,161)$\\
$(62,175)$&
$(63,160)$&
$(155,207)$&
$(25,114)$&
$(87,208)$&
$(90,209)$&
$(100,212)$&
$(109,211)$&
$(39,108)$&
$(11,102)$\\
$(83,124)$&
$(56,117)$&
$(59,130)$&
$(84,188)$&
$(104,151)$&
$(72,178)$&
$(82,177)$&
$(89,179)$&
$(94,181)$&
$(23,128)$\\
$(95,167)$&
$(99,126)$&
$(123,162)$&
$(91,129)$&
$(92,136)$&
$(81,158)$&
$(113,159)$&
$(116,201)$&
$(137,196)$&
$(33,140)$\\
$(35,144)$&
$(37,145)$&
$(1,18)$&
$(2,20)$&
$(3,22)$&
$(4,24)$&
$(5,26)$&
$(6,28)$&
$(7,30)$&
$(8,32)$\\
\end{array}$$

$$\begin{array}{llllllllll}

\bullet\ \ $u=126:$\\
$(124,219)$&
$(111,174)$&
$(112,172)$&
$(119,192)$&
$(105,197)$&
$(89,169)$&
$(90,180)$&
$(92,178)$&
$(120,185)$&
$(122,176)$\\
$(129,222)$&
$(154,223)$&
$(115,203)$&
$(133,184)$&
$(134,198)$&
$(135,205)$&
$(113,242)$&
$(128,241)$&
$(138,244)$&
$(144,243)$\\
$(91,212)$&
$(97,236)$&
$(110,194)$&
$(140,249)$&
$(102,214)$&
$(101,225)$&
$(103,238)$&
$(136,193)$&
$(152,246)$&
$(153,228)$\\
$(159,240)$&
$(168,239)$&
$(121,247)$&
$(137,252)$&
$(145,250)$&
$(146,208)$&
$(195,251)$&
$(125,204)$&
$(83,215)$&
$(130,231)$\\
$(141,237)$&
$(127,182)$&
$(63,160)$&
$(155,207)$&
$(79,190)$&
$(100,227)$&
$(109,226)$&
$(87,201)$&
$(98,235)$&
$(163,245)$\\
$(114,248)$&
$(187,253)$&
$(117,166)$&
$(142,220)$&
$(149,202)$&
$(84,188)$&
$(104,151)$&
$(139,213)$&
$(150,217)$&
$(156,206)$\\
$(94,181)$&
$(77,162)$&
$(80,211)$&
$(93,161)$&
$(95,167)$&
$(116,175)$&
$(123,230)$&
$(108,218)$&
$(107,232)$&
$(126,229)$\\
\end{array}$$
}

\section{Partition the set $[1,2u+1]\setminus\{(u-30)/2\}$ into pairs}

Here for each $u\in[68,259]$ and $u\equiv 2\ ({\rm mod}\ 4)$, we give a partition of the set $[1,2u+1]\setminus\{(u-30)/2\}$ into pairs $\{x_i,y_i\}$, $1\leq i\leq u$, such that $\{y_i-x_i:1\leq i\leq u\}=[33,32+u]$.

{\scriptsize
$$\begin{array}{llllllllll}
\bullet\ \ $u=70:$&\\
$(38,81)$&
$(44,131)$&
$(45,135)$&
$(50,96)$&
$(52,92)$&
$(54,116)$&
$(32,95)$&
$(33,114)$&
$(35,117)$&
$(47,130)$\\
$(55,132)$&
$(57,136)$&
$(51,124)$&
$(18,68)$&
$(19,70)$&
$(63,98)$&
$(65,120)$&
$(13,58)$&
$(56,110)$&
$(53,94)$\\
$(15,99)$&
$(16,101)$&
$(17,105)$&
$(48,137)$&
$(14,106)$&
$(21,87)$&
$(22,83)$&
$(26,91)$&
$(27,118)$&
$(30,97)$\\
$(31,90)$&
$(37,107)$&
$(41,109)$&
$(43,80)$&
$(46,85)$&
$(60,113)$&
$(77,115)$&
$(59,134)$&
$(61,139)$&
$(28,86)$\\
$(64,133)$&
$(72,119)$&
$(73,121)$&
$(74,123)$&
$(67,111)$&
$(69,126)$&
$(1,34)$&
$(2,62)$&
$(3,89)$&
$(4,84)$\\
$(5,100)$&
$(6,102)$&
$(7,108)$&
$(8,79)$&
$(9,103)$&
$(10,82)$&
$(11,104)$&
$(12,112)$&
$(23,75)$&
$(24,88)$\\
$(25,122)$&
$(29,128)$&
$(36,78)$&
$(39,141)$&
$(40,138)$&
$(42,76)$&
$(49,125)$&
$(66,140)$&
$(71,127)$&
$(93,129)$\\
\end{array}$$}

{\scriptsize
$$\begin{array}{llllllllll}

\bullet\ \ $u=82:$\\
$(39,104)$&
$(51,117)$&
$(38,143)$&
$(54,130)$&
$(56,113)$&
$(60,108)$&
$(63,127)$&
$(65,121)$&
$(31,125)$&
$(33,136)$\\
$(35,135)$&
$(55,141)$&
$(17,101)$&
$(18,124)$&
$(64,145)$&
$(19,116)$&
$(72,131)$&
$(77,139)$&
$(83,144)$&
$(84,138)$\\
$(91,134)$&
$(29,142)$&
$(73,158)$&
$(85,152)$&
$(45,133)$&
$(46,106)$&
$(53,122)$&
$(58,105)$&
$(59,103)$&
$(67,107)$\\
$(68,109)$&
$(93,165)$&
$(22,120)$&
$(15,126)$&
$(16,115)$&
$(32,146)$&
$(47,155)$&
$(69,118)$&
$(28,123)$&
$(49,156)$\\
$(66,157)$&
$(36,129)$&
$(42,100)$&
$(95,150)$&
$(75,164)$&
$(79,159)$&
$(44,148)$&
$(81,119)$&
$(41,128)$&
$(43,114)$\\
$(86,161)$&
$(88,151)$&
$(90,160)$&
$(94,147)$&
$(23,132)$&
$(70,153)$&
$(1,34)$&
$(2,37)$&
$(3,40)$&
$(4,96)$\\
$(5,82)$&
$(6,48)$&
$(7,89)$&
$(8,87)$&
$(9,111)$&
$(10,78)$&
$(11,57)$&
$(12,102)$&
$(13,52)$&
$(14,110)$\\
\end{array}$$}

{\scriptsize
$$\begin{array}{llllllllll}

\bullet\ \ $u=106:$\\
$(39,104)$&
$(51,117)$&
$(54,130)$&
$(56,113)$&
$(60,108)$&
$(63,127)$&
$(65,121)$&
$(31,125)$&
$(33,136)$&
$(35,135)$\\
$(61,171)$&
$(82,172)$&
$(97,195)$&
$(110,196)$&
$(102,186)$&
$(114,167)$&
$(77,139)$&
$(83,144)$&
$(99,178)$&
$(30,154)$\\
$(132,190)$&
$(148,197)$&
$(95,150)$&
$(115,183)$&
$(116,209)$&
$(23,152)$&
$(24,131)$&
$(58,105)$&
$(78,194)$&
$(89,166)$\\
$(68,109)$&
$(93,165)$&
$(66,192)$&
$(67,185)$&
$(111,189)$&
$(34,168)$&
$(62,147)$&
$(32,160)$&
$(47,155)$&
$(100,204)$\\
$(119,201)$&
$(87,184)$&
$(94,193)$&
$(19,134)$&
$(107,188)$&
$(81,177)$&
$(84,205)$&
$(79,159)$&
$(80,181)$&
$(85,202)$\\
$(101,175)$&
$(126,199)$&
$(27,158)$&
$(28,141)$&
$(88,151)$&
$(71,173)$&
$(72,191)$&
$(70,153)$&
$(50,180)$&
$(120,174)$\\
$(123,156)$&
$(96,210)$&
$(103,208)$&
$(128,203)$&
$(91,182)$&
$(98,157)$&
$(112,164)$&
$(129,179)$&
$(142,213)$&
$(133,200)$\\
\end{array}$$}

{\scriptsize
$$\begin{array}{llllllllll}

$(89,166)$&
$(68,109)$&
$(93,165)$&
$(66,192)$&
$(34,147)$&
$(133,200)$&
$(136,244)$&
$(149,252)$&
$(95,220)$&
$(98,221)$\\
$(54,206)$&
$(59,210)$&
$(49,177)$&
$(119,201)$&
$(87,184)$&
$(94,193)$&
$(19,134)$&
$(97,242)$&
$(74,170)$&
$(88,190)$\\
$(84,205)$&
$(79,159)$&
$(80,181)$&
$(85,202)$&
$(101,175)$&
$(126,199)$&
$(142,197)$&
$(92,219)$&
$(31,155)$&
$(32,164)$\\
$(33,179)$&
$(115,227)$&
$(146,238)$&
$(112,234)$&
$(100,250)$&
$(120,213)$&
$(118,236)$&
$(141,209)$&
$(130,241)$&
$(73,187)$\\
$(123,156)$&
$(171,237)$&
$(104,154)$&
$(116,225)$&
$(111,253)$&
$(151,249)$&
$(103,208)$&
$(128,203)$&
$(131,218)$&
$(51,198)$\\
$(53,173)$&
$(1,35)$&
$(2,37)$&
$(3,39)$&
$(4,41)$&
$(5,43)$&
$(6,45)$&
$(7,47)$&
$(8,50)$&
$(9,52)$\\
$(10,55)$&
$(11,57)$&
$(12,56)$&
$(13,61)$&
$(14,71)$&
$(15,64)$&
$(16,67)$&
$(17,69)$&
$(18,72)$&
$(20,90)$\\
$(21,125)$&
$(22,82)$&
$(23,152)$&
$(24,158)$&
$(25,180)$&
$(26,145)$&
$(27,162)$&
$(28,161)$&
$(30,140)$&
$(36,172)$\\
\end{array}$$}

{\scriptsize
$$\begin{array}{llllllllll}

\bullet\ \ $u=146:$\\
$(155,237)$&
$(185,248)$&
$(108,251)$&
$(64,234)$&
$(71,215)$&
$(65,121)$&
$(73,242)$&
$(78,231)$&
$(91,174)$&
$(157,245)$\\
$(143,212)$&
$(29,160)$&
$(87,260)$&
$(112,238)$&
$(114,167)$&
$(94,184)$&
$(126,283)$&
$(142,289)$&
$(117,195)$&
$(124,189)$\\
$(129,229)$&
$(81,235)$&
$(86,224)$&
$(113,243)$&
$(151,246)$&
$(177,275)$&
$(187,263)$&
$(89,166)$&
$(68,109)$&
$(93,165)$\\
$(193,282)$&
$(204,271)$&
$(156,226)$&
$(136,244)$&
$(149,252)$&
$(95,220)$&
$(98,221)$&
$(59,210)$&
$(139,197)$&
$(152,225)$\\
$(96,264)$&
$(138,287)$&
$(123,257)$&
$(76,216)$&
$(104,214)$&
$(88,190)$&
$(84,205)$&
$(79,159)$&
$(80,181)$&
$(194,281)$\\
$(92,219)$&
$(168,239)$&
$(173,230)$&
$(74,240)$&
$(83,180)$&
$(99,276)$&
$(49,223)$&
$(144,286)$&
$(100,250)$&
$(120,213)$\\
$(118,236)$&
$(141,209)$&
$(130,241)$&
$(103,262)$&
$(56,201)$&
$(115,227)$&
$(131,266)$&
$(153,285)$&
$(171,278)$&
$(179,253)$\\
\end{array}$$}

{\scriptsize
$$\begin{array}{llllllllll}
\bullet\ \ $u=162:$\\
$(185,248)$&
$(193,282)$&
$(164,292)$&
$(176,321)$&
$(73,242)$&
$(78,231)$&
$(91,174)$&
$(20,75)$&
$(21,67)$&
$(22,82)$\\
$(29,160)$&
$(87,260)$&
$(137,299)$&
$(153,305)$&
$(94,184)$&
$(126,283)$&
$(142,289)$&
$(117,195)$&
$(124,189)$&
$(129,229)$\\
$(81,235)$&
$(83,269)$&
$(155,249)$&
$(163,276)$&
$(162,297)$&
$(187,263)$&
$(89,166)$&
$(68,109)$&
$(93,165)$&
$(54,236)$\\
$(181,303)$&
$(156,226)$&
$(136,244)$&
$(143,308)$&
$(95,220)$&
$(98,221)$&
$(59,210)$&
$(139,197)$&
$(152,225)$&
$(96,264)$\\
$(138,287)$&
$(123,257)$&
$(169,271)$&
$(77,255)$&
$(146,279)$&
$(128,312)$&
$(110,214)$&
$(141,209)$&
$(154,293)$&
$(168,239)$\\
$(173,230)$&
$(74,240)$&
$(111,267)$&
$(122,237)$&
$(172,301)$&
$(144,286)$&
$(100,250)$&
$(25,104)$&
$(116,178)$&
$(132,280)$\\
$(130,241)$&
$(103,262)$&
$(175,302)$&
$(127,186)$&
$(131,314)$&
$(134,311)$&
$(171,278)$&
$(179,253)$&
$(158,274)$&
$(147,284)$\\
\end{array}$$}

\hspace{1cm}
{\scriptsize
$$\begin{array}{llllllllll}
$(25,209)$&
$(26,230)$&
$(27,85)$&
$(28,106)$&
$(30,82)$&
$(35,88)$&
$(45,102)$&
$(70,133)$&
$(89,272)$&
$(100,256)$\\
$(111,176)$&
$(117,227)$&
$(157,328)$&
$(168,224)$&
$(173,258)$&
$(185,321)$&
$(192,337)$&
$(195,300)$&
$(222,334)$&
$(239,318)$\\
$(260,335)$&
$(261,332)$&
$(267,327)$&
$(277,341)$\\

\bullet\ \ $u=178:$\\
$(170,349)$&
$(204,322)$&
$(212,298)$&
$(225,339)$&
$(234,315)$&
$(290,344)$&
$(111,318)$&
$(112,314)$&
$(118,200)$&
$(124,215)$\\
$(73,242)$&
$(78,231)$&
$(91,174)$&
$(20,75)$&
$(120,313)$&
$(164,292)$&
$(29,160)$&
$(143,340)$&
$(266,346)$&
$(153,305)$\\
$(94,184)$&
$(126,283)$&
$(142,289)$&
$(137,299)$&
$(166,338)$&
$(129,229)$&
$(81,235)$&
$(83,269)$&
$(155,249)$&
$(163,276)$\\
$(162,297)$&
$(187,263)$&
$(99,183)$&
$(93,165)$&
$(54,236)$&
$(181,303)$&
$(156,226)$&
$(136,244)$&
$(217,312)$&
$(95,220)$\\
\end{array}$$}

{\scriptsize
$$\begin{array}{llllllllll}
$(98,221)$&
$(59,210)$&
$(256,360)$&
$(113,205)$&
$(96,264)$&
$(78,244)$&
$(81,235)$&
$(169,271)$&
$(77,255)$&
$(146,279)$\\
$(228,295)$&
$(105,285)$&
$(109,296)$&
$(154,293)$&
$(47,115)$&
$(84,245)$&
$(86,250)$&
$(25,67)$&
$(149,254)$&
$(123,272)$\\
$(144,286)$&
$(193,282)$&
$(153,370)$&
$(192,359)$&
$(209,361)$&
$(130,241)$&
$(261,332)$&
$(267,352)$&
$(127,186)$&
$(63,238)$\\
$(19,65)$&
$(21,74)$&
$(11,55)$&
$(158,274)$&
$(147,284)$&
$(103,175)$&
$(106,171)$&
$(49,190)$&
$(110,188)$&
$(152,343)$\\
$(107,265)$&
$(150,320)$&
$(167,288)$&
$(33,207)$&
$(281,371)$&
$(287,364)$&
$(108,251)$&
$(182,291)$&
$(101,247)$&
$(97,273)$\\
$(69,268)$&
$(37,240)$&
$(41,196)$&
$(43,243)$&
$(248,345)$&
$(202,275)$&
$(203,333)$&
$(206,309)$&
$(39,189)$&
$(159,252)$\\
$(177,294)$&
$(52,233)$&
$(121,319)$&
$(145,330)$&
$(56,180)$&
$(61,148)$&
$(31,232)$&
$(32,197)$&
$(201,270)$&
$(246,342)$\\
$(208,307)$&
$(71,259)$&
$(219,325)$&
$(51,139)$&
$(58,198)$&
$(191,329)$&
$(194,326)$&
$(229,374)$&
$(253,379)$&
$(131,323)$\\
\end{array}$$}

{\scriptsize
$$\begin{array}{llllllllll}
$(237,373)$&
$(257,365)$&
$(195,351)$&
$(214,390)$&
$(302,358)$&
$(225,339)$&
$(245,391)$&
$(178,341)$&
$(161,355)$&
$(116,223)$\\
$(23,85)$&
$(218,378)$&
$(224,322)$&
$(232,348)$&
$(247,362)$&
$(33,134)$&
$(41,182)$&
$(37,222)$&
$(199,388)$&
$(254,372)$\\
$(163,298)$&
$(211,385)$&
$(140,370)$&
$(149,331)$&
$(162,366)$&
$(158,381)$&
$(179,334)$&
$(306,363)$&
$(143,356)$&
$(111,336)$\\
$(123,319)$&
$(173,399)$&
$(191,329)$&
$(1,34)$&
$(2,36)$&
$(3,38)$&
$(4,40)$&
$(5,42)$&
$(6,44)$&
$(7,46)$\\
$(8,48)$&
$(9,50)$&
$(10,53)$&
$(12,54)$&
$(14,62)$&
$(15,64)$&
$(20,75)$&
$(22,88)$&
$(24,126)$&
$(25,196)$\\
$(26,84)$&
$(28,138)$&
$(30,207)$&
$(31,76)$&
$(43,271)$&
$(45,272)$&
$(49,112)$&
$(67,300)$&
$(69,268)$&
$(80,212)$\\
$(82,168)$&
$(100,249)$&
$(104,335)$&
$(114,314)$&
$(118,313)$&
$(120,273)$&
$(128,292)$&
$(131,328)$&
$(135,349)$&
$(141,347)$\\
$(164,357)$&
$(169,403)$&
$(184,404)$&
$(190,401)$&
$(194,396)$&
$(200,392)$&
$(204,395)$&
$(229,374)$&
$(236,380)$&
$(277,405)$\\
\end{array}$$}

\hspace{1cm}
{\scriptsize
$$\begin{array}{llllllllll}
$(218,378)$&
$(104,249)$&
$(271,403)$&
$(247,362)$&
$(33,134)$&
$(41,182)$&
$(37,222)$&
$(199,388)$&
$(254,372)$&
$(224,327)$\\
$(273,401)$&
$(73,203)$&
$(84,277)$&
$(162,366)$&
$(69,289)$&
$(43,168)$&
$(306,363)$&
$(116,357)$&
$(118,292)$&
$(309,426)$\\
$(120,356)$&
$(125,304)$&
$(314,380)$&
$(349,412)$&
$(229,334)$&
$(220,419)$&
$(221,374)$&
$(230,414)$&
$(204,361)$&
$(240,395)$\\
$(245,391)$&
$(212,375)$&
$(301,424)$&
$(341,425)$&
$(1,34)$&
$(2,36)$&
$(3,38)$&
$(4,40)$&
$(5,42)$&
$(6,44)$\\
$(7,46)$&
$(8,48)$&
$(9,50)$&
$(10,53)$&
$(12,54)$&
$(14,62)$&
$(15,64)$&
$(20,75)$&
$(22,67)$&
$(24,191)$\\
$(31,169)$&
$(45,161)$&
$(76,319)$&
$(80,272)$&
$(82,140)$&
$(87,283)$&
$(90,190)$&
$(99,322)$&
$(108,294)$&
$(119,313)$\\
$(123,337)$&
$(126,351)$&
$(129,369)$&
$(131,370)$&
$(133,297)$&
$(138,333)$&
$(150,387)$&
$(153,355)$&
$(164,409)$&
$(170,405)$\\
$(176,383)$&
$(177,421)$&
$(183,416)$&
$(192,402)$&
$(195,392)$&
$(200,348)$&
$(206,423)$&
$(207,418)$&
$(214,390)$&
$(232,415)$\\
\end{array}$$}

\hspace{1cm}
{\scriptsize
$$\begin{array}{llllllllll}
$(247,362)$&
$(33,134)$&
$(41,182)$&
$(37,222)$&
$(199,388)$&
$(300,429)$&
$(224,327)$&
$(168,284)$&
$(94,333)$&
$(121,359)$\\
$(283,390)$&
$(69,289)$&
$(56,144)$&
$(90,254)$&
$(306,363)$&
$(116,357)$&
$(118,292)$&
$(309,426)$&
$(258,444)$&
$(264,401)$\\
$(314,380)$&
$(349,412)$&
$(169,422)$&
$(170,427)$&
$(130,366)$&
$(139,331)$&
$(204,361)$&
$(214,318)$&
$(245,391)$&
$(212,375)$\\
$(301,424)$&
$(205,393)$&
$(208,445)$&
$(193,282)$&
$(207,432)$&
$(219,325)$&
$(246,372)$&
$(203,448)$&
$(221,364)$&
$(232,449)$\\
$(234,446)$&
$(287,383)$&
$(261,428)$&
$(302,358)$&
$(304,453)$&
$(153,335)$&
$(166,417)$&
$(173,369)$&
$(1,34)$&
$(2,36)$\\
$(3,38)$&
$(4,40)$&
$(5,42)$&
$(6,44)$&
$(7,46)$&
$(8,48)$&
$(9,50)$&
$(10,53)$&
$(12,54)$&
$(13,60)$\\
$(14,62)$&
$(15,64)$&
$(17,68)$&
$(18,70)$&
$(20,225)$&
$(31,96)$&
$(59,210)$&
$(73,328)$&
$(75,213)$&
$(80,180)$\\
$(95,351)$&
$(100,155)$&
$(101,348)$&
$(105,298)$&
$(106,339)$&
$(108,220)$&
$(109,344)$&
$(110,188)$&
$(113,296)$&
$(119,368)$\\
\end{array}$$}

\hspace{1cm}
{\scriptsize
$$\begin{array}{llllllllll}
$(167,288)$&
$(49,184)$&
$(72,151)$&
$(308,410)$&
$(97,305)$&
$(218,441)$&
$(229,400)$&
$(18,70)$&
$(20,109)$&
$(93,311)$\\
$(39,189)$&
$(122,277)$&
$(125,341)$&
$(52,233)$&
$(80,283)$&
$(280,340)$&
$(58,200)$&
$(76,286)$&
$(243,377)$&
$(32,197)$\\
$(107,265)$&
$(60,164)$&
$(71,118)$&
$(414,475)$&
$(142,408)$&
$(84,338)$&
$(87,158)$&
$(195,404)$&
$(198,405)$&
$(251,398)$\\
$(206,420)$&
$(35,117)$&
$(185,386)$&
$(62,136)$&
$(85,258)$&
$(145,373)$&
$(146,330)$&
$(190,370)$&
$(202,275)$&
$(262,407)$\\
$(271,434)$&
$(236,435)$&
$(318,372)$&
$(331,393)$&
$(15,64)$&
$(133,384)$&
$(139,285)$&
$(91,174)$&
$(247,362)$&
$(33,134)$\\
$(41,182)$&
$(37,222)$&
$(199,388)$&
$(300,429)$&
$(224,327)$&
$(168,284)$&
$(94,333)$&
$(220,469)$&
$(245,381)$&
$(69,289)$\\
$(56,144)$&
$(73,150)$&
$(306,363)$&
$(116,357)$&
$(152,402)$&
$(214,409)$&
$(309,426)$&
$(113,169)$&
$(264,401)$&
$(314,380)$\\
$(349,412)$&
$(252,392)$&
$(149,328)$&
$(171,367)$&
$(221,466)$&
$(228,465)$&
$(231,423)$&
$(105,213)$&
$(106,359)$&
$(261,467)$\\
\end{array}$$}

\hspace{1cm}
{\scriptsize
$$\begin{array}{llllllllll}
$(140,256)$&
$(162,228)$&
$(168,445)$&
$(169,394)$&
$(180,444)$&
$(199,388)$&
$(204,428)$&
$(225,463)$&
$(241,439)$&
$(264,432)$\\
$(309,380)$&
$(320,480)$&
$(365,401)$&
$(416,469)$&
$(422,468)$&
$(424,474)$\\

\bullet\ \ $u=250:$\\
$(7,113)$&
$(8,48)$&
$(302,485)$&
$(388,444)$&
$(424,495)$&
$(57,209)$&
$(86,316)$&
$(124,215)$&
$(315,440)$&
$(323,415)$\\
$(157,376)$&
$(165,397)$&
$(414,475)$&
$(266,346)$&
$(118,290)$&
$(237,448)$&
$(255,494)$&
$(264,461)$&
$(22,67)$&
$(24,82)$\\
$(137,299)$&
$(178,378)$&
$(243,468)$&
$(249,487)$&
$(102,350)$&
$(111,326)$&
$(114,336)$&
$(187,263)$&
$(260,324)$&
$(278,353)$\\
$(79,269)$&
$(181,303)$&
$(156,226)$&
$(227,354)$&
$(217,312)$&
$(155,417)$&
$(320,498)$&
$(141,343)$&
$(298,492)$&
$(317,488)$\\
$(206,433)$&
$(220,396)$&
$(228,420)$&
$(250,464)$&
$(84,360)$&
$(106,273)$&
$(96,270)$&
$(100,368)$&
$(47,115)$&
$(361,473)$\\
\end{array}$$}

\hspace{1cm}
{\scriptsize
$$\begin{array}{llllllllll}
$(74,179)$&
$(76,297)$&
$(116,357)$&
$(191,462)$&
$(194,364)$&
$(71,325)$&
$(83,239)$&
$(160,400)$&
$(161,419)$&
$(349,412)$\\
$(252,392)$&
$(230,517)$&
$(231,466)$&
$(238,401)$&
$(450,505)$&
$(463,509)$&
$(105,213)$&
$(261,467)$&
$(232,351)$&
$(234,460)$\\
$(141,343)$&
$(142,403)$&
$(149,328)$&
$(150,387)$&
$(366,472)$&
$(415,477)$&
$(245,381)$&
$(251,398)$&
$(259,436)$&
$(287,383)$\\
$(313,456)$&
$(327,459)$&
$(304,453)$&
$(153,335)$&
$(312,474)$&
$(113,397)$&
$(126,279)$&
$(127,272)$&
$(129,402)$&
$(154,293)$\\
$(240,394)$&
$(323,512)$&
$(330,513)$&
$(288,391)$&
$(68,331)$&
$(282,470)$&
$(307,437)$&
$(339,465)$&
$(166,423)$&
$(167,377)$\\
$(17,56)$&
$(19,73)$&
$(347,395)$&
$(318,478)$&
$(358,476)$&
$(152,356)$&
$(23,301)$&
$(27,158)$&
$(337,446)$&
$(332,413)$\\
$(218,441)$&
$(90,359)$&
$(241,439)$&
$(172,416)$&
$(348,471)$&
$(365,499)$&
$(299,504)$&
$(320,498)$&
$(369,507)$&
$(375,482)$\\
$(376,496)$&
$(406,503)$&
$(205,455)$&
$(208,483)$&
$(212,454)$&
$(1,34)$&
$(2,36)$&
$(3,38)$&
$(4,40)$&
$(5,42)$\\
$(6,44)$&
$(7,48)$&
$(8,195)$&
$(9,78)$&
$(10,53)$&
$(12,54)$&
$(13,221)$&
$(14,243)$&
$(16,66)$&
$(21,61)$\\
$(26,271)$&
$(28,295)$&
$(31,188)$&
$(59,225)$&
$(65,344)$&
$(75,322)$&
$(97,211)$&
$(98,262)$&
$(119,338)$&
$(120,190)$\\
$(140,372)$&
$(156,445)$&
$(157,256)$&
$(162,434)$&
$(183,417)$&
$(219,484)$&
$(224,510)$&
$(226,422)$&
$(242,506)$&
$(286,486)$\\
$(292,514)$&
$(302,418)$&
$(305,485)$&
$(336,511)$&
$(367,500)$&
$(370,458)$&
$(378,429)$&
$(449,516)$\\
\end{array}$$}

\section{$(4t,4,\{3,5\},1)$-CDFs for $t\equiv 1\ ({\rm mod}\ 2)$, $t\geq 7$ and $t\neq11$.}

Here we list all base blocks of a $(4t,4,\{3,5\},1)$-CDF for $t\equiv 1\ ({\rm mod}\ 2)$, $t\geq 7$ and $t\neq11$.

$(1)$ When $t\equiv 1\ ({\rm mod}\ 6)$ and $t\geq 7$, the conclusion follows immediately from Lemma $3.5$ in $[33]$.

$(2)$ When $t=29$, a $(4t,4,\{3,5\},1)$-CDF is listed below:

\vspace{0.1cm}
\begin{tabular}{llllll}
$\{0,3,69\},$ &$\{0,11,86\},$ & $\{0,16,83\},$ & $\{0,18,82\},$ & $\{0,24,77\},$ & $\{0,1,6,20,27\},$ \\
$\{0,4,74\},$ &$\{0,12,85\},$ & $\{0,17,79\},$ & $\{0,22,78\},$ & $\{0,28,76\},$ & $\{0,2,10,25,61\},$ \\
$\{0,9,81\},$ &$\{0,13,84\}.$ &&&& \\
\end{tabular}

\vspace{0.1cm}
\noindent When $t=35$, a $(4t,4,\{3,5\},1)$-CDF is listed below:

\vspace{0.1cm}
\begin{tabular}{llllll}
$\{0,1,32,78,104\},$ &$\{0,4,20\},$ & $\{0,7,34\},$ & $\{0,10,96\},$ & $\{0,13,100\},$ & $\{0,25,81\},$ \\
$\{0,3,33,76,98\},$ &$\{0,5,23\},$ & $\{0,8,82\},$ & $\{0,11,102\},$ & $\{0,14,99\},$ & $\{0,28,80\},$\\
$\{0,2,19\},$ &$\{0,6,21\},$ & $\{0,9,92\},$ & $\{0,12,101\},$ & $\{0,24,71\},$ & $\{0,29,79\}.$\\
\end{tabular}

\vspace{0.1cm}
\noindent When $t\equiv 5,11\ ({\rm mod}\ 24)$ and $t\geq 53$, take $A_1=\{0,1,t-3,2t+8,3t-1\}$ and  $A_2=\{0,3,t-2,2t+6,3t-7\}$ as two base blocks with block size five. Let $S=([1,t-1]\cup[2t+1,3t-1])$ and $T=\{1,3,t-13,t-9,t-5,t-4,t-3,t-2,2t+2,2t+3,2t+5,2t+6,2t+7,2t+8,3t-11,3t-10,3t-8,3t-7,3t-2,
3t-1\}$. Then we shall show that $S\setminus T$ can be partitioned into triples $\{a_i,b_i,c_i\},$ $1\leq i\leq (2t-22)/3,$ such that $a_i+b_i\equiv c_i\ ({\rm mod}\ 4t)$.

By Lemma $4.12$ with $(d,u,k)=(3,(t-17)/6,(t-17)/3)$, the set $[3,(t-11)/2]\setminus \{(t-13)/2\}$ can be partitioned into triples $\{a'_i,b'_i,c'_i\},$ $1\leq i\leq (t-17)/6,$ such that $a'_i+b'_i=c'_i.$ Taking $a_i=2a'_i,$ $b_i=2b'_i$ and $c_i=2c'_i$ for $1\leq i\leq (t-17)/6,$ we know that $[6,t-11]_{e}\setminus \{t-13\}$ can be partitioned into triples $\{a_i,b_i,c_i\},$ $1\leq i\leq (t-17)/6,$ such that $a_i+b_i=c_i.$

Furthermore, $(S\setminus T)\setminus ([6,t-11]_{e}\setminus \{t-13\})$ can be partitioned into triples $\{a_i,b_i,c_i\},$ $(t-11)/6\leq
i\leq (2t-22)/3$, such that $a_i+b_i\equiv c_i\ ({\rm mod}\ 4t)$ as follows:

\vspace{0.1cm}
\begin{tabular}{ll}
$\{13+2j,(5t-15)/2-j,(5t+11)/2+j\},$ & $j\in[0,(t-35)/2]\setminus\{(t-37)/2\};$
\end{tabular}

\begin{tabular}{lll}
$\{4,t-24,t-20\},$ & $\{t-12,2t+9,3t-3\},$ & $\{2,(5t+1)/2,(5t+5)/2\},$\\
$\{11,t-18,t-7\},$ & $\{t-14,2t+1,3t-13\},$ & $\{(5t-13)/2,(5t-3)/2,t-8\},$\\
$\{5,3t-9,3t-4\},$ & $\{9,(5t-9)/2,(5t+9)/2\},$ & $\{(5t-11)/2,(5t-1)/2,t-6\},$\\
$\{t-16,2t+11,3t-5\},$ & $\{7,(5t-7)/2,(5t+7)/2\},$ & $\{(5t-5)/2,(5t+3)/2,t-1\},$\\
$\{t-10,2t+4,3t-6\}.$ &  & \\
\end{tabular}

\vspace{0.1cm}
\noindent Then $\{\{0,a_i,c_i\}:1\leq i\leq (2t-22)/3\}\cup\{A_1,A_2\}$ forms a $(4t,4,\{3,5\},1)$-CDF.

\vspace{0.1cm}
$(3)$ When $t=17$, a $(4t,4,\{3,5\},1)$-CDF is listed below:

\vspace{0.1cm}
\begin{tabular}{llllll}
$\{0,1,5,11,48\},$ &$\{0,3,16,35,44\},$ & $\{0,2,14\},$ & $\{0,7,46\},$ & $\{0,8,50\},$ & $\{0,15,38\}.$
\end{tabular}

\vspace{0.1cm}
\noindent When $t=23$, a $(4t,4,\{3,5\},1)$-CDF is listed below:

\vspace{0.1cm}
\begin{tabular}{llllll}
$\{0,1,5,11,20\},$ &$\{0,7,55\},$ & $\{0,12,64\},$ & $\{0,17,67\},$ & $\{0,21,56\},$ & $\{0,22,53\},$ \\
$\{0,2,49,62,65\},$ &$\{0,8,66\},$ & $\{0,14,68\},$ & $\{0,18,51\}.$ &  & \\
\end{tabular}

\vspace{0.1cm}
\noindent When $t=41$, a $(4t,4,\{3,5\},1)$-CDF is listed below:

{\tabcolsep 0.07in
\begin{tabular}{llllll}
$\{0,1,38,90,122\},$ &$\{0,2,25\},$ & $\{0,7,35\},$ & $\{0,11,110\},$ & $\{0,15,107\},$ & $\{0,18,119\},$ \\
$\{0,3,39,86,116\},$ &$\{0,4,26\},$ & $\{0,8,93\},$ & $\{0,12,115\},$ & $\{0,16,118\},$ & $\{0,20,111\},$\\
$\{0,19,114\},$ &$\{0,5,29\},$ & $\{0,9,40\},$ & $\{0,13,109\},$ & $\{0,17,105\},$ & $\{0,34,94\},$\\
$\{0,33,100\},$ &$\{0,6,27\},$ & $\{0,10,108\},$ & $\{0,14,120\}.$ &  & \\
\end{tabular}}

\vspace{0.1cm}
\noindent When $t\equiv 17,23\ ({\rm mod}\ 24)$ and $t\geq47$, take $A_1=\{0,1,t-3,2t+8,3t-1\}$ and $A_2=\{0,3,t-2,2t+4,3t-7\}$ as two base blocks with block size five. Let $S=([1,t-1]\cup[2t+1,3t-1])$ and $T=\{1,3,t-11,t-9,t-5,t-4,t-3,t-2,2t+1,2t+2,2t+4,2t+5,2t+7,2t+8,3t-11,3t-10,3t-7,3t-6,3t-2,
3t-1\}$.

By Lemma $4.11$ with $(d,u)=(3,(t-17)/6)$, the set $[3,(t-13)/2]$ can be partitioned into triples $\{a'_i,b'_i,c'_i\},$ $1\leq i\leq (t-17)/6,$ such that $a'_i+b'_i=c'_i.$ Taking $a_i=2a'_i,$ $b_i=2b'_i$ and $c_i=2c'_i$ for $1\leq i\leq (t-17)/6,$ we know that $[6,t-13]_{e}$ can be partitioned into triples $\{a_i,b_i,c_i\},$ $1\leq i\leq (t-17)/6,$ such that $a_i+b_i=c_i.$

Furthermore, $(S\setminus T)\setminus [6,t-13]_{e}$ can be partitioned into triples $\{a_i,b_i,c_i\},$ $(t-11)/6\leq
i\leq (2t-22)/3$, such that $a_i+b_i\equiv c_i\ ({\rm mod}\ 4t)$ as follows:

\vspace{0.1cm}
\begin{tabular}{ll}
$\{11+2j,(5t-13)/2-j,(5t+9)/2+j\},$ & $j\in[0,(t-37)/2];$
\end{tabular}

\begin{tabular}{lll}
$\{7,t-14,t-7\},$ & $\{t-12,2t+3,3t-9\},$ & $\{5,(5t-5)/2,(5t+5)/2\},$\\
$\{4,t-24,t-20\},$ & $\{t-22,2t+9,3t-13\},$ & $\{(5t-9)/2,(5t-7)/2,t-8\},$\\
$\{9,3t-12,3t-3\},$ & $\{t-18,2t+10,3t-8\},$ & $\{(5t-11)/2,(5t-1)/2,t-6\},$\\
$\{t-16,2t+11,3t-5\},$ & $\{2,(5t+3)/2,(5t+7)/2\},$ & $\{(5t-3)/2,(5t+1)/2,t-1\},$\\
$\{t-10,2t+6,3t-4\}.$ &  & \\
\end{tabular}

\vspace{0.1cm}
\noindent Then $\{\{0,a_i,c_i\}:1\leq i\leq (2t-22)/3\}\cup\{A_1,A_2\}$ forms a $(4t,4,\{3,5\},1)$-CDF.

\vspace{0.1cm}
$(4)$ When $t=15$, a $(4t,4,\{3,5\},1)$-CDF is listed below:

\vspace{0.1cm}{\tabcolsep 0.07in
\begin{tabular}{lllllll}
$\{0,1,3,12,20\}$,&
$\{0,4,25\}$,&
$\{0,5,28\}$,&
$\{0,6,24\}$,&
$\{0,7,29\}$,&
$\{0,10,26\}$,
$\{0,13,27\}$.
\end{tabular}}

\vspace{0.1cm}
\noindent When $t\equiv 15,21\ ({\rm mod}\ 24)$ and $t\geq21$, take $A=\{0,1,t-2,2t+3,3t-1\}$ as the base block with block size five. Let $S=([1,t-1]\cup[2t+1,3t-1])$ and $T=\{1,t-4,t-3,t-2,2t+1,2t+2,2t+3,3t-5,3t-2,3t-1\}$.

By Lemma $4.12$ with $(d,u,k)=(2,(t-9)/6,t/3-3)$, the set $[2,(t-5)/2]\setminus \{(t-7)/2\}$ can be partitioned into triples $\{a'_i,b'_i,c'_i\},$ $1\leq i\leq (t-9)/6,$ such that $a'_i+b'_i=c'_i.$ Taking $a_i=2a'_i,$ $b_i=2b'_i$ and $c_i=2c'_i$ for $1\leq i\leq (t-9)/6,$ we know that $[4,t-5]_{e}\setminus \{t-7\}$ can be partitioned into triples $\{a_i,b_i,c_i\},$ $1\leq i\leq (t-9)/6,$ such that $a_i+b_i=c_i.$

Furthermore, $(S\setminus T)\setminus ([4,t-5]_{e}\setminus \{t-7\})$ can be partitioned into triples $\{a_i,b_i,c_i\},$ $(t-3)/6\leq
i\leq 2t/3-4$, such that $a_i+b_i\equiv c_i\ ({\rm mod}\ 4t)$ as follows:

\vspace{0.1cm}
\begin{tabular}{ll}
$\{7+2j,(5t-9)/2-j,(5t+5)/2+j\},$ & $j\in[0,(t-17)/2]\setminus\{(t-19)/2\};$
\end{tabular}

\begin{tabular}{lll}
$\{5,t-12,t-7\},$ & $\{t-8,2t+5,3t-3\},$ & $\{(5t-7)/2,(5t-5)/2,t-6\},$\\
$\{3,3t-7,3t-4\},$ & $\{2,(5t-1)/2,(5t+3)/2\},$ & $\{(5t-3)/2,(5t+1)/2,t-1\}.$\\
\end{tabular}

\vspace{0.1cm}
\noindent Then $\{\{0,a_i,c_i\}:1\leq i\leq 2t/3-4\}\cup\{A\}$ forms a $(4t,4,\{3,5\},1)$-CDF.

\vspace{0.1cm}
$(5)$ When $t=9$, a $(4t,4,\{3,5\},1)$-CDF is listed below:

\vspace{0.1cm}
\begin{tabular}{lllllll}
$\{0,1,3,13,17\}$, &$\{0,5,11\}$,& $\{0,7,15\}$.
\end{tabular}

\vspace{0.1cm}
\noindent When $t\equiv 3,9\ ({\rm mod}\ 24)$ and $t\geq27$, take $A=\{0,1,t-2,2t+3,3t-1\}$ as the base block with block size five. Let $S=([1,t-1]\cup[2t+1,3t-1])$ and $T=\{1,t-4,t-3,t-2,2t+1,2t+2,2t+3,3t-5,3t-2,3t-1\}$.

By Lemma $4.11$ with $(d,u)=(2,(t-9)/6)$, the set $[2,(t-7)/2]$ can be partitioned into triples $\{a'_i,b'_i,c'_i\},$ $1\leq i\leq (t-9)/6,$ such that $a'_i+b'_i=c'_i.$ Taking $a_i=2a'_i,$ $b_i=2b'_i$ and $c_i=2c'_i$ for $1\leq i\leq (t-9)/6,$ we know that $[4,t-7]_{e}$ can be partitioned into triples $\{a_i,b_i,c_i\},$ $1\leq i\leq (t-9)/6,$ such that $a_i+b_i=c_i.$

Furthermore, $(S\setminus T)\setminus [4,t-7]_{e}$ can be partitioned into triples $\{a_i,b_i,c_i\},$ $(t-3)/6\leq
i\leq 2t/3-4$, such that $a_i+b_i\equiv c_i\ ({\rm mod}\ 4t)$ as follows:

\vspace{0.1cm}
\begin{tabular}{ll}
$\{7+2j,(5t-9)/2-j,(5t+5)/2+j\},$ & $j\in[0,(t-19)/2];$
\end{tabular}

\begin{tabular}{lll}
$\{5,t-10,t-5\},$ & $\{t-8,2t+4,3t-4\},$ & $\{(5t-7)/2,(5t-5)/2,t-6\},$\\
$\{3,3t-6,3t-3\},$ & $\{2,(5t-1)/2,(5t+3)/2\},$ & $\{(5t-3)/2,(5t+1)/2,t-1\}.$\\
\end{tabular}

\vspace{0.1cm}
\noindent Then $\{\{0,a_i,c_i\}:1\leq i\leq 2t/3-4\}\cup\{A\}$ forms a $(4t,4,\{3,5\},1)$-CDF.

\section{$(16t,16,\{3,5\},1)$-CDFs for $t\equiv 0\ ({\rm mod}\ 2)$ and $t\geq 4$}

Here we list all base blocks of a $(16t,16,\{3,5\},1)$-CDF for $t\equiv 0\ ({\rm mod}\ 2)$ and $t\geq 4$.

$(1)$ When $t\equiv 4\ ({\rm mod}\ 6)$ and $t\geq 4$, the conclusion follows immediately from Lemma $3.8$ in $[33]$.

$(2)$ When $t=12$, a $(16t,16,\{3,5\},1)$-CDF is listed below:

\vspace{0.1cm}{\tabcolsep 0.06in
\begin{tabular}{lllllll}
$\{0,1,34\},$ &$\{0,6,128\},$ &$\{0,11,126\},$ &$\{0,19,142\},$ &$\{0,27,136\},$&$\{0,30,119\},$\\
$\{0,2,37\},$ &$\{0,7,131\},$ &$\{0,15,140\},$ &$\{0,20,137\},$ &$\{0,28,134\},$&$\{0,31,44\},$\\
$\{0,3,41\},$ &$\{0,8,105\},$ &$\{0,16,127\},$ &$\{0,21,135\},$ &$\{0,18,130\},$&$\{0,32,46\},$\\
$\{0,4,43\},$ &$\{0,9,138\},$ &$\{0,17,102\},$ &$\{0,26,139\},$ &$\{0,29,133\},$&$\{0,25,47,118,141\},$\\
$\{0,5,45\},$ &$\{0,10,110\},$ &$\{0,42,143\}.$ & &&\\
\end{tabular}}

\vspace{0.1cm}
\noindent When $t=24$, a $(16t,16,\{3,5\},1)$-CDF is listed below:

\vspace{0.1cm}
\begin{tabular}{ll}
$\{0,11+2j,244+j\},$& $j\in[0,16];$ \\
$\{0,53+2j,266+j\},$& $j\in[0,17];$ \\
\end{tabular}

\begin{tabular}{lllllll}
$\{0,1,45\},$ &$\{0,5,89\},$ &$\{0,9,91\},$ &$\{0,16,78\},$ &$\{0,26,287\},$&$\{0,74,265\},$\\
$\{0,2,30\},$ &$\{0,6,40\},$ &$\{0,10,66\},$ &$\{0,18,88\},$ &$\{0,42,284\},$&$\{0,90,239\},$\\
$\{0,3,54\},$ &$\{0,7,93\},$ &$\{0,12,76\},$ &$\{0,20,263\},$ &$\{0,52,286\},$&$\{0,92,214\},$\\
$\{0,4,36\},$ &$\{0,8,58\},$ &$\{0,14,94\},$ &$\{0,22,60\},$ &$\{0,68,237\},$&$\{0,49,95,238,285\}.$\\
\end{tabular}

\vspace{0.1cm}
\noindent When $t\equiv 0\ ({\rm mod}\ 12)$ and $t\geq36$, take $A=\{0,2t+1,4t-1,10t-2,12t-3\}$ as the base block with block size five. Let $S=([1,4t-1]\cup[8t+1,12t-1])\setminus \{t,2t,3t,9t,10t,11t\}$ and $T_1=\{2t-2,2t-1,2t+1,4t-1,8t+2,8t+3,10t-4,10t-2,10t+1,12t-3\}$.

By Lemma $4.12$ with $(d,u,k)=(3,t/6-2,t/12)$, the set $[3,t/2-3]\setminus \{t/4\}$ can be partitioned into triples $\{a'_i,b'_i,c'_i\},$ $1\leq i\leq t/6-2,$ such that $a'_i+b'_i=c'_i.$ Taking $a_i=4a'_i,$ $b_i=4b'_i$ and $c_i=4c'_i$ for $1\leq i\leq t/6-2,$ we have that $T_2=\{4r:\ r=3,4,\ldots,t/2-3\}\setminus \{t\}$ can be partitioned into triples $\{a_i,b_i,c_i\},$ $1\leq i\leq t/6-2,$ such that $a_i+b_i=c_i.$

Furthermore, $(S\setminus T_1)\setminus T_2$ can be partitioned into triples $\{a_i,b_i,c_i\},$ $t/6-1\leq
i\leq 8t/3-6$, such that $a_i+b_i\equiv c_i\ ({\rm mod}\ 16t)$ as follows:

\vspace{0.1cm}
\begin{tabular}{ll}
$\{11+2j,10t-7-j,10t+4+j\},$ & $j\in[0,t-8]\setminus\{2,t/2-5\};$ \\
$\{2t+5+2j,9t-3-j,11t+2+j\},$ &$j\in[0,t-7];$ \\
$\{6+4j,3t-2-2j,3t+4+2j\},$ &$j\in[0,t/2-3];$ \\
\end{tabular}

\begin{tabular}{lll}
$\{7,8,15\},$ & $\{t+1,11t-3,12t-2\},$ & $\{9t-2,9t-1,2t-3\},$\\
$\{2,3,5\},$ & $\{2t-8,2t+3,4t-5\},$ & $\{10t-9,10t+2,4t-7\},$\\
$\{4,10t-5,10t-1\},$ & $\{2t+2,10t-6,12t-4\},$ & $\{19t/2-2,21t/2-1,4t-3\},$\\
$\{9,10t-3,10t+6\},$ & $\{2t-4,10t+3,12t-1\},$ & $\{8t+1,11t+1,3t+2\},$\\
$\{1,11t-2,11t-1\}.$ &  & \\
\end{tabular}

\vspace{0.1cm}
\noindent Then $\{\{0,a_i,c_i\}:1\leq i\leq 8t/3-6\}\cup\{A\}$ forms a $(16t,16,\{3,5\},1)$-CDF.

$(3)$ When $t=14$, a $(16t,16,\{3,5\},1)$-CDF is listed below:

\vspace{0.1cm}{\tabcolsep 0.06in
\begin{tabular}{lllllll}
$\{0,1,39\},$ &$\{0,6,52\},$ &$\{0,11,163\},$ &$\{0,17,155\},$ &$\{0,26,156\},$&$\{0,33,53\},$\\
$\{0,2,43\},$ &$\{0,7,123\},$ &$\{0,12,47\},$ &$\{0,19,147\},$ &$\{0,27,151\},$&$\{0,34,165\},$\\
$\{0,3,40\},$ &$\{0,8,127\},$ &$\{0,13,159\},$ &$\{0,21,162\},$ &$\{0,30,166\},$&$\{0,49,148\},$\\
$\{0,4,48\},$ &$\{0,9,129\},$ &$\{0,15,160\},$ &$\{0,36,158\},$ &$\{0,32,153\},$&$\{0,29,54,143,161\},$\\
$\{0,5,50\},$ &$\{0,10,149\},$ &$\{0,16,150\},$ &$\{0,51,118\},$ &$\{0,23,167\},$&$\{0,31,55,142,164\}.$\\
\end{tabular}}

\vspace{0.1cm}
\noindent When $t=26$, a $(16t,16,\{3,5\},1)$-CDF is listed below:

\vspace{0.1cm}
\begin{tabular}{ll}
$\{0,13+2j,264+j\},$& $j\in[0,16];$ \\
$\{0,57+2j,287+j\},$& $j\in[0,14];$ \\
\end{tabular}

{\tabcolsep 0.06in
\begin{tabular}{lllllll}
$\{0,1,51\},$ &$\{0,7,54\},$ &$\{0,14,70\},$ &$\{0,20,303\},$ &$\{0,74,306\},$&$\{0,96,258\},$\\
$\{0,2,38\},$ &$\{0,8,94\},$ &$\{0,16,60\},$ &$\{0,24,309\},$ &$\{0,76,231\},$&$\{0,98,281\},$\\
$\{0,3,93\},$ &$\{0,9,97\},$ &$\{0,22,84\},$ &$\{0,34,66\},$ &$\{0,89,101\},$&$\{0,99,259\},$\\
$\{0,5,87\},$ &$\{0,10,40\},$ &$\{0,4,311\},$ &$\{0,68,282\},$ &$\{0,92,304\},$&$\{0,53,102,263,305\},$\\
$\{0,6,64\},$ &$\{0,11,91\},$ &$\{0,18,302\},$ &$\{0,72,100\},$ &$\{0,95,310\},$&$\{0,55,103,262,308\}.$\\
\end{tabular}}

\vspace{0.1cm}
\noindent When $t\equiv 2\ ({\rm mod}\ 12)$ and $t\geq38$, take $A_1=\{0,2t+1,4t-2,10t+3,12t-7\}$ and $A_2=\{0,2t+3,4t-1,10t+2,12t-4\}$ as two base blocks with block size five. Let $S=([1,4t-1]\cup[8t+1,12t-1])\setminus \{t,2t,3t,9t,10t,11t\}$ and $T_1=\{2t-10,2t-6,2t-4,2t-3,2t+1,2t+3,4t-2,4t-1,8t+1,8t+2,8t+3,8t+5,10t-8,10t-7,10t-5,10t-3,
10t+2,10t+3,12t-7,12t-4\}$.

By Lemma $4.12$ with $(d,u,k)=(2,(t-8)/6,(t+10)/12)$, the set $[2,t/2-2]\setminus \{(t+2)/4\}$ can be partitioned into triples $\{a'_i,b'_i,c'_i\},$ $1\leq i\leq (t-8)/6,$ such that $a'_i+b'_i=c'_i.$ Taking $a_i=4a'_i,$ $b_i=4b'_i$ and $c_i=4c'_i$ for $1\leq i\leq (t-8)/6,$ we have that $T_2=\{4r:\ r=2,3,\ldots,t/2-2\}\setminus \{t+2\}$ can be partitioned into triples $\{a_i,b_i,c_i\},$ $1\leq i\leq (t-8)/6,$ such that $a_i+b_i=c_i.$

Furthermore, $(S\setminus T_1)\setminus T_2$ can be partitioned into triples $\{a_i,b_i,c_i\},$ $(t-2)/6\leq
i\leq (8t-28)/3$, such that $a_i+b_i\equiv c_i\ ({\rm mod}\ 16t)$ as follows:

\vspace{0.1cm}
\begin{tabular}{ll}
$\{13+2j,10t-9-j,10t+4+j\},$& $j\in[0,t-10]\setminus\{t-11\};$ \\
$\{2t+5+2j,9t-4-j,11t+1+j\},$& $j\in[0,t-10]\setminus\{t/2+4\};$ \\
$\{6+4j,3t-2-2j,3t+4+2j\},$& $j\in[0,t/2-5]\setminus\{(t-6)/4\};$ \\
\end{tabular}

\begin{tabular}{lll}
$\{3,12t-8,12t-5\},$ & $\{2t-1,10t-2,12t-3\},$ & $\{8t+4,7t/2+1,23t/2+5\},$\\
$\{5,12t-6,12t-1\},$ & $\{2t+4,10t-6,12t-2\},$ & $\{10t-4,10t+1,4t-3\},$\\
$\{2,4t-11,4t-9\},$ & $\{2t-5,9t+2,11t-3\},$ & $\{9t-2,11t-5,4t-7\},$\\
$\{9,4t-13,4t-4\},$ & $\{t+2,9t-3,10t-1\},$ & $\{9t-1,11t-4,4t-5\},$\\
$\{4,2t+2,2t+6\},$ & $\{1,11t-2,11t-1\},$ & $\{5t/2+1,17t/2-8,11t-7\},$\\
$\{7,2t-9,2t-2\},$ & $\{11,3t+2,3t+13\}.$ & \\
\end{tabular}

\vspace{0.1cm}
\noindent Then $\{\{0,a_i,c_i\}:1\leq i\leq (8t-28)/3\}\cup\{A_1,A_2\}$ forms a $(16t,16,\{3,5\},1)$-CDF.

$(4)$ When $t=6$, the conclusion follows from Example $2.9$. When $t\equiv 6\ ({\rm mod}\ 12)$ and $t\geq 18$, take $A=\{0,2t+1,4t-1,10t-2,12t-3\}$ as the base block with block size five. Let $S=([1,4t-1]\cup[8t+1,12t-1])\setminus \{t,2t,3t,9t,10t,11t\}$ and $T_1=\{2t-2,2t-1,2t+1,4t-1,8t+2,8t+3,10t-4,10t-2,10t+1,12t-3\}$.

By Lemma $4.12$ with $(d,u,k)=(2,t/6-1,(t+6)/12)$, the set $[2,t/2-1]\setminus \{(t+2)/4\}$ can be partitioned into triples $\{a'_i,b'_i,c'_i\},$ $1\leq i\leq t/6-1,$ such that $a'_i+b'_i=c'_i.$ Taking $a_i=4a'_i,$ $b_i=4b'_i$ and $c_i=4c'_i$ for $1\leq i\leq t/6-1,$ we have that $T_2=\{4r:\ r=2,3,\ldots,t/2-1\}\setminus \{t+2\}$ can be partitioned into triples $\{a_i,b_i,c_i\},$ $1\leq i\leq t/6-1,$ such that $a_i+b_i=c_i.$

Furthermore, $(S\setminus T_1)\setminus T_2$ can be partitioned into triples $\{a_i,b_i,c_i\},$ $t/6\leq
i\leq 8t/3-6$, such that $a_i+b_i\equiv c_i\ ({\rm mod}\ 16t)$ as follows:

\vspace{0.1cm}
\begin{tabular}{ll}
$\{9+2j,10t-6-j,10t+3+j\},$& $j\in[0,t-7]\setminus\{2\};$ \\
$\{2t+5+2j,9t-3-j,11t+2+j\},$& $j\in[0,t-7]\setminus\{t/2\};$ \\
$\{10+4j,3t-4-2j,3t+6+2j\},$ &$j\in[0,t/2-4]\setminus\{(t-10)/4\};$ \\
\end{tabular}

\begin{tabular}{lll}
$\{6,7,13\},$ & $\{2,4t-7,4t-5\},$ & $\{2t+2,9t-1,11t+1\},$\\
$\{1,3t+4,3t+5\},$ & $\{t+2,11t-3,12t-1\},$ & $\{8t+1,7t/2+1,23t/2+2\},$\\
$\{4,3t-2,3t+2\},$ & $\{2t+3,10t-5,12t-2\},$ & $\{5t/2+1,17t/2-3,11t-2\},$\\
$\{3,10t+2,10t+5\},$ & $\{2t-3,10t-1,12t-4\},$ & $\{9t-2,11t-1,4t-3\},$\\
$\{5,10t-8,10t-3\}.$ &  & \\
\end{tabular}

\vspace{0.1cm}
\noindent Then $\{\{0,a_i,c_i\}:1\leq i\leq 8t/3-6\}\cup\{A\}$ forms a $(16t,16,\{3,5\},1)$-CDF.

$(5)$ When $t=8$, the conclusion follows from Example $2.10$. When $t\equiv 8\ ({\rm mod}\ 12)$ and $t\geq20$, take $A_1=\{0,2t+1,4t-2,10t+3,12t-7\}$ and $A_2=\{0,2t+3,4t-1,10t+2,12t-4\}$ as two base blocks with block size five. Let $S=([1,4t-1]\cup[8t+1,12t-1])\setminus \{t,2t,3t,9t,10t,11t\}$ and $T_1=\{2t-10,2t-6,2t-4,2t-3,2t+1,2t+3,4t-2,4t-1,8t+1,8t+2,8t+3,8t+5,10t-8,10t-7,10t-5,
10t-3,10t+2,10t+3,12t-7,12t-4\}$.

By Lemma $4.12$ with $(d,u,k)=(2,(t-8)/6,(t+4)/12)$, the set $[2,t/2-2]\setminus \{t/4\}$ can be partitioned into triples $\{a'_i,b'_i,c'_i\},$ $1\leq i\leq (t-8)/6,$ such that $a'_i+b'_i=c'_i.$ Taking $a_i=4a'_i,$ $b_i=4b'_i$ and $c_i=4c'_i$ for $1\leq i\leq (t-8)/6,$ we have that $T_2=\{4r:\ r=2,3,\ldots,t/2-2\}\setminus \{t\}$ can be partitioned into triples $\{a_i,b_i,c_i\},$ $1\leq i\leq (t-8)/6,$ such that $a_i+b_i=c_i.$

Furthermore, $(S\setminus T_1)\setminus T_2$ can be partitioned into triples $\{a_i,b_i,c_i\},$ $(t-2)/6\leq
i\leq (8t-28)/3$, such that $a_i+b_i\equiv c_i\ ({\rm mod}\ 16t)$ as follows:

\vspace{0.1cm}
\begin{tabular}{ll}
$\{13+2j,10t-9-j,10t+4+j\},$& $j\in[0,t-10]\setminus\{1,t/2-5\};$ \\
$\{2t+5+2j,9t-4-j,11t+1+j\},$& $j\in[0,t-10];$ \\
$\{6+4j,3t-2-2j,3t+4+2j\},$& $j\in[0,t/2-5];$ \\
\end{tabular}

\begin{tabular}{lll}
$\{4,11,15\},$ & $\{2,4t-13,4t-11\},$ & $\{10t-2,10t-1,4t-3\},$\\
$\{3,2t-5,2t-2\},$ & $\{2t+4,10t-6,12t-2\},$ & $\{9t-2,11t-5,4t-7\},$\\
$\{7,2t-1,2t+6\},$ & $\{2t+2,9t-3,11t-1\},$ & $\{9t-1,11t-3,4t-4\},$\\
$\{1,12t-6,12t-5\},$ & $\{t+3,11t-4,12t-1\},$ & $\{19t/2-4,21t/2-1,4t-5\},$\\
$\{5,12t-8,12t-3\},$ & $\{10t-10,10t+1,4t-9\},$ & $\{8t+4,11t-2,3t+2\},$\\
$\{9,10t-4,10t+5\}.$ &  & \\
\end{tabular}

\vspace{0.1cm}
\noindent Then $\{\{0,a_i,c_i\}:1\leq i\leq (8t-28)/3\}\cup\{A_1,A_2\}$ forms a $(16t,16,\{3,5\},1)$-CDF.

\end{document}